\documentclass[12pt]{amsart}

\usepackage{amsfonts}
\usepackage{latexsym}
\usepackage{amssymb}
\usepackage{amscd}
\usepackage{graphics}
\usepackage{amsmath}
\usepackage{mathrsfs}

\input xy
\xyoption{all}

\newcommand{\cA}{\mathcal{A}}
\newcommand{\cB}{\mathcal{B}}

\newcommand{\cG}{\mathcal{G}}

\newcommand{\cL}{\mathcal{L}}
\newcommand{\cN}{\mathcal{N}}

\newcommand{\cP}{\mathcal{P}}

\newcommand{\cS}{\mathcal{S}}

\newcommand{\mcH}{\mathcal{H}}
\newcommand{\cM}{\mathcal{M}}

\newcommand{\ca}{\mathfrak{a}}
\newcommand{\cb}{\mathfrak{b}}
\newcommand{\cg}{\mathfrak{g}}
\newcommand{\ch}{\mathfrak{h}}

\newcommand{\cp}{\mathfrak{p}}
\newcommand{\cn}{\mathfrak{n}}
\newcommand{\cs}{\mathfrak{s}}
\newcommand{\cz}{\mathfrak{z}}

\newcommand{\cu}{\mathfrak{u}}
\newcommand{\cw}{\mathfrak{w}}
\newcommand{\cv}{\mathfrak{v}}

\newcommand{\G}{\mathbb{G}}

\renewcommand{\H}{\mathbb{H}}

\newcommand{\M}{\mathbb{M}}

\newcommand{\N}{\mathbb{N}}

\newcommand{\R}{\mathbb{R}}
\newcommand{\U}{\mathbb{U}}
\renewcommand{\P}{\mathbb{P}}

\renewcommand{\H}{\mathbb{H}}

\newcommand{\ep}{\varepsilon}

\newcommand{\ph}{\varphi}

\newcommand{\sm}{\setminus}
\newcommand{\dist}{\mbox{dist}}
\newcommand{\gope}{\circledcirc}

\newcommand{\lan}{\langle}
\newcommand{\ran}{\rangle}
\newcommand{\lls}{\mbox{\large $($}}
\newcommand{\rls}{\mbox{\large $)$}}
\newcommand{\lLs}{\mbox{\Large $($}}
\newcommand{\rLs}{\mbox{\Large $)$}}

\newcommand{\ra}{\rightarrow}
\newcommand{\lra}{\longrightarrow}

\newcommand{\medint}{\hbox{\vrule height3.5pt depth-2.8pt
width4pt}\mkern-14mu\int\nolimits}

\newcommand{\diam}{\mbox{diam}}
\newcommand{\Exp}{\mbox{\rm Exp}}

\newcommand{\Lip}{\mbox{Lip}}

\renewcommand{\span}{\mbox{span}}

\renewcommand{\dist}{\mbox{\rm dist}}

\newcommand{\Id}{\mbox{Id}}

\newtheorem{The}{Theorem}[section]
\newtheorem{Lem}[The]{Lemma}
\newtheorem{Def}{Definition}[section]
\newtheorem{Rem}[The]{Remark}
\newtheorem{Pro}[The]{Proposition}
\newtheorem{Cor}[The]{Corollary}
\newtheorem{Exa}[The]{Example}

\begin{document}

\title[Towards Differential Calculus in stratified groups]
{{\bf Towards Differential Calculus in stratified groups}}
\author{Valentino Magnani}
\address{Valentino Magnani: Dipartimento di Matematica
\\ Largo Bruno Pontecorvo 5 \\ $\;\;$ I-56127 Pisa}
\email{magnani@dm.unipi.it}
\maketitle

\vskip-8truecm
\tableofcontents

\pagebreak

\section{Introduction}

The relationship between the geometry of stratified groups
and several branches of Mathematics, as PDEs, Differential Geometry,
Complex Analysis, Control Theory and Geometric Measure Theory
has known an increasing interest in the last decade. 

The initial motivation of the present work was the study of an
implicit function theorem between stratified groups
with respect to the corresponding notion differentiability.
Then we realized that to formulate and prove this theorem
a series of basic results were needed. These ones naturally lead us to the 
development of an Intrinsic Differential Calculus for group-valued mappings
defined on stratified groups. 

This project belongs to the general program of developing Analysis in
metric spaces, as explained in several books \cite{AmbTil}, \cite{BBI}, \cite{DavSem}, \cite{HajKos}, \cite{Hein1}, \cite{Semmes}.
Other monographs are focused on Lie groups and sub-Riemannian manifolds, \cite{BelRis}, \cite{CDPT}, \cite{FS}, \cite{Gr2}, \cite{Montgom}, \cite{VSC}.
Calculus and Analysis in stratified groups have been developed by several 
authors with a vast literature, here we mention a list of references for a rough overview that is certainly far from being complete,
\cite{ASCV}, \cite{ArcMor}, \cite{Bal}, \cite{BHIT}, \cite{BSCV}, \cite{CapCow},
\cite{CG},  \cite{CHMY}, \cite{CitMan}, \cite{DGN1}, \cite{DGN3}, \cite{FSSC4}, \cite{FSSC5}, \cite{FSSC6}, \cite{GN}, \cite{Hein}, \cite{KirSer}, 
\cite{LeoMas}, \cite{LeoRig}, \cite{Mag5}, \cite{Mag8A}, \cite{MagVit}, 
\cite{Mat05}, \cite{MonRick}, \cite{MonSer}, \cite{Pau2}, \cite{RitRos}.

In the present work the notion of differentiability plays a central role
and essentially enters every result.
This notion in the setting of stratified groups has been introduced by P. Pansu
to obtain a rigidity theorem for quasi-isometries of the quaternionic
hyperbolic space, where he proved a Rademacher type theorem
for Lipschitz mappings with respect to his notion of differentiability, \cite{Pan2}.

Our {\em intrinsically $C^1$-smooth Calculus} starts from the following 
\begin{The}[Characterization of continuously P-differentiable mappings]\label{PdifContact} 
Let $\Omega\subset\G$ be an open set and let $f:\Omega\lra\M$.
Then the following statements are equivalent
\begin{itemize}
	\item $f$ is continuously P-differentiable\,,
	\item letting $f(x)=\exp\circ F(x)=\exp\sum_{j=1}^\upsilon F_j(x)$, 
the mappings $F_j:\Omega\lra W_j$,
equal to $\pi_j\circ F$, are continuously P-differentiable and the formulae
\begin{eqnarray}
&&\pi_1\circ df(x)=dF_1(x) \label{eqF_11} \\
&&dF_j(x)(h)=\sum_{n=2}^\upsilon\frac{(-1)^n}{n!}\;\pi_j
\left([F(x),dF(x)(h)]_{n-1}\right)\,. \label{eqDF_j1} \label{eqdF}
\end{eqnarray}
hold for every $j=2,\ldots,\upsilon$ and every $h\in\G$.
 	\item $f$ is continuously h-differentiable contact\,.
\end{itemize}	
The integer $\upsilon$ denotes the step of $\M$ and
the symbol $d$ denotes the P-differential read between Lie algebras.
\end{The}
We have denoted by $\G$ and $\M$ a stratified group and a graded group, 
respectively. The corresponding Lie algebras $\cG$ of $\G$ 
and $\cM$ of $\M$ are decomposed into the direct sums of layers
$V_i$ and $W_i$, respectively. The mappings $\pi_j$ indicate the canonical
projection onto the $j$-th layer of a graded algebra.

Notice that h-differentiability, introduced in Definition~\ref{hdiff},
is a much weaker notion than P-differentiability. According to
Theorem~\ref{PdifContact}, we have that continuous h-differentiability
and contact property given by \eqref{eqdF} yield continuous P-differentiability.
The main point of Theorem~\ref{PdifContact} is that it allows for the study of
P-differentiable mappings regardless of the notion of P-differentiability
and using instead the following system of first order nonlinear PDEs
\begin{equation}\label{hcondF}
\left\{\begin{array}{c}
X_1F_j=\sum_{n=2}^\upsilon\frac{(-1)^n}{n!}\;\pi_j\left([F,X_1F]_{n-1}\right)\\
X_2F_j=\sum_{n=2}^\upsilon\frac{(-1)^n}{n!}\;\pi_j\left([F,X_2F]_{n-1}\right)\\
\!\!\!\!\!\vdots\qquad\vdots\qquad\vdots\qquad\vdots\qquad\vdots\qquad\vdots\qquad\vdots \\
X_mF_j=\sum_{n=2}^\upsilon\frac{(-1)^n}{n!}\;\pi_j\left([F,X_mF]_{n-1}\right)
\end{array}\right.
\end{equation}
for every $j=2,\ldots,\upsilon$, where $(X_1,\ldots,X_m)$ is a basis
of the first layer of $\cG$. In fact, formulae \eqref{eqDF_j1} recursively define the
P-differential $dF_j$ for every $j=2,\ldots,\upsilon$, starting from $dF_1$
and in particular yield the system \eqref{hcondF}.
These $m\times(\upsilon-1)$ equations exactly characterize the contact
property of $f$ and can be used to study existence of Lipschitz extensions
for mappings between stratified groups, \cite{Mag10}.
Another feature of  these equations is that they allow to transmit
the regularity of $F_1$ to the remaining components $F_j$ of higher
order. The same phenomenon occurs in the ``algebraic regularity''
proved in Theorem~6.1 of \cite{CapCow}.

The crucial point in the proof of Theorem~\ref{PdifContact} is the
quantitative estimate \eqref{pansuestim} for the P-difference quotient
of a horizontal curve, that is obtained in Theorem~\ref{pansu}.
This estimate plays a key role also in the proof of the mean value inequality and it allows for the extension of Pansu result, \cite{Pan2},
on a.e. P-differentiability of Lipschitz mappings to the case when the
target is a graded group, see Corollary~\ref{Pansugrad}.
On the other hand, this key estimate is obtained in a rather
``simple'' way, once the technical Lemma~\ref{keylemma} is adopted,
and it boils down to integrate recursively the differential equations
\begin{eqnarray}\label{cntctcurve}
	\dot\gamma_i(t)=\sum_{n=2}^\upsilon\frac{(-1)^n}{n!}\;
	\pi_i\left([\gamma(t),\dot\gamma(t)]_{n-1}\right), \qquad i=2,\ldots,\upsilon,
\end{eqnarray}
that must be satisfied by every horizontal curve $\Gamma=\exp\gamma\in\M$,
where $\gamma_i=\pi_i\circ \gamma$ and $i=2,\ldots,\upsilon$,
according to Proposition~\ref{horprop}.
Lemma~\ref{keylemma} is one of the major technical results of this paper
and will appear several times throughout the present work. It
corresponds to the linearization of the homogeneous addends $c_n$
appearing in the Baker-Campbell-Hausdorff formula \eqref{absBCH}.

If one assumes apriori $C^1$ smoothness of $f$, then the
characterization of contact mappings through P-differentiabiltiy
is already shown in \cite{War}. 
The intriguing feature of Theorem~\ref{PdifContact}, along with its technical 
difficulty, stems from the fact that we do not assume any ``extrinsic'' regularity of $f$, making this characterization sharp.
Continuously h-differentiable mappings are continuous, since their components are locally Lipschitz with respect to the homogeneous distance of the domain. Furthermore, adding the contact property to the continuous h-differentiability,
then Theorem~\ref{PdifContact} and Corollary~\ref{corlip} yield the local Lipschitz property of $f$ with respect to homogeneous distances.
However, let us point out that one can find continuously h-differentiable contact mappings that are not differentiable on sets of positive measure with respect to the Euclidean notion of differentiability, \cite{Mag3}. 

Our next result is the mean value inequality
for P-differentiable mappings, that represents a crucial tool
in the proof of our implicit function theorem.
Even in this case, the proof relies on a preliminary
study of horizontal curves. In Corollary~\ref{hdifflip} we characterize 
Lipschitz curves in graded groups with respect to a homogeneous distance
as Lipschitz curves in the Euclidean sense that also
satisfy the differential equations \eqref{cntctcurve}.
In the case of stratified groups, where the Lie bracket condition
is satisfied, this is a well known fact that can be found in
Section~11.1 of \cite{HajKos} for the more general Carnot-Carath\'eodory
spaces.

\begin{The}[Mean Value Inequality]\label{unifestim}
Let $\Omega\subset\G$ be an open subset and consider
a continuously P-differentiable mapping $f:\Omega\lra\M$.
Let $\Omega_1,\Omega_2\subset\G$ be open subsets such that 
\begin{equation}\label{inclomega1}
\Big\{x\in\G\,\Big|\,d(x,\Omega_1)\leq c(\G,d)\,
N\,\mbox{\rm diam}(\Omega_1)\Big\}\subset\Omega_2\,,
\end{equation}
where $\Omega_2$ is compactly contained in $\Omega$.
Here $c(\G,d)$ and $N$ are geometric constants only depending on
the metric space $(\G,d)$, see both Lemma~\ref{140FS} and Definition~\ref{defcaDef}.
Then there exist a constant $C$, only depending on $\G$, $\max_{x\in\overline{\Omega}_2}\|dF_1(x)\|$ and on 
the modulus of continuity $\omega_{\overline{\Omega}_2,dF_1}$
of $x\lra df(x)$, defined in \eqref{modcont}, such that
\begin{equation}\label{keyunifest}
\frac{\rho\Big(f(x)^{-1}f(y),Df(x)(x^{-1}y)\Big)}{d(x,y)}\leq 
C\;\;\big[\omega_{\overline{\Omega}_2,df}\big(N\,c(\G,d)\,d(x,y)\big)
\big]^{1/\iota^2}\,.
\end{equation}
for every $x,y\in\overline{\Omega}_1$, with $x\neq y$.
The integer $\iota$ denotes the step of $\G$.
\end{The}
Several difficulties are hidden in this estimate. 
First, P-differentiable mappings in general are not $C^1$ smooth
in the classical sense. Second, working by single components
does not suffice, since this would lead to estimates on the
Euclidean norm of the difference quotient, that does not fit into
the notion of P-differentiability.
Third, the family of horizontal curves along which one integrates
the P-differential is not manageable.
In fact, condition \eqref{inclomega1} has only the technical motivation
to make sure that the special family of piecewise horizontal lines
that connect points of $\Omega_1$ and along which we apply \eqref{pansuestim}
are all contained in $\Omega_2$. 

An immediate corollary of the mean value inequality is a short proof of
the inverse mapping theorem in stratified groups.
On the other hand, the difficulty in applying Theorem~\ref{unifestim}
to prove an intrinsic implicit function theorem arises from the algebraic
problems related to a proper factorization of the group $\G$.
In the commutative case, that is the classical viewpoint,
the implicit mapping is defined by decomposing the space into a product
of two linear subspaces, naturally given by the tangent space
to the level set and one of its complementary subspaces,
as for instance the orthogonal subspace, if a scalar product is fixed.
Extending this argument to stratified groups does not work when
we seek the second subspace.
In fact, the first one is automatically defined as the kernel of
the P-differential of the defining mapping, that is also 
a normal subgroup, but in general we cannot claim a group
structure in the ``complementary subspace''. 
For this reason, according to our terminology, we consider the special
class of {\em h-epimorphisms}, as the ``regular surjective h-homomorphisms", that yield the natural splitting of $\G$, when they are the P-differential
of the defining mapping. Notice that this splitting of the group is a
necessary condition in order to state an implicit function theorem.
In fact, the existence of the complementary subgroup $H$ in Theorem~\ref{implth} is a consequence of Proposition~\ref{GlinH}, since $Df(\overline x)$ is an
h-epimorphism. To emphasize this special class of surjective h-homomorphisms,
in Example~\ref{exahr2}, we provide a surjective h-homomorphism that is not an h-epimorphism.

Notice that not all Lie subgroups are considered, but only the special class of 
homogeneous subgroups, that are closed under dilations.
This turns out to be rather natural, in the perspective to
interpret a subgroup as blow-up of a suitable intrinsically regular
set, according to Theorem~\ref{blwlevelimage}. 
Using the terminology of Section~\ref{homsubgr},
two homogeneous subgroups $P$ and $H$ are {\em complementary}
if they satisfy the conditions $\G=PH$ and $P\cap H=\{e\}$.

Now we are in the position to establish one of the central results of this paper.
\begin{The}[Implicit Function Theorem]\label{implth}
Let $\Omega\subset\G$ be an open set and
let $f:\Omega\lra\M$ be continuously P-differentiable,
where $\overline{x}\in\Omega$ and the P-differential $Df(\overline{x}):\G\lra\M$
is an h-epimorphism.
Let $N$ be the kernel of $Df(\overline{x})$ and let $H$
be a complementary subgroup. 
Then there exist $r,s>0$, with
$D^N_{\overline{n},r}D^H_{\overline{h},s}\subset\Omega$,
along with a unique mapping
$\ph:D^N_{\overline{n},r}\lra D^H_{\overline{h},s}$, such that
\begin{eqnarray}\label{impleq}
f^{-1}\lls f(\overline{x})\rls\cap D^N_{\overline{n},r}D^H_{\overline{h},s}
=\{n \ph(n)\mid n\in D^N_{\overline{n},r}\}.
\end{eqnarray}
Furthermore, there exists a constant $\kappa>0$ such that the Lipschitz-type estimate
\begin{eqnarray}\label{intrlip}
d\big(\ph(n),\ph(n')\big)\leq \kappa \;d\big(\ph(n')^{-1}n^{-1}n'\ph(n')\big)
\end{eqnarray}
holds. In particular, the mapping $\ph$ is $1/\iota$-H\"older continuous
with respect to the metrics $d$ in $D^H_{\overline{h},r}$
and $\|\cdot\|$ in $D^N_{\overline{n},s}$.
\end{The}
The induced splitting of $\G$ is precisely an inner semidirect
product, that is in general not direct. This should explain 
why the classical contraction-mapping principle here
seems to be not applicable.
More precisely, to use this argument, in the case $\overline{x}$ equals
the unit element $e\in\G$, one should consider the mapping $F(h)=L(h)^{-1}f(nh)$, 
where $L$ is the restriction of the P-differential $Df(e)$ to
the complementary subgroup $H$.
In the classical case, one shows that the mapping $F$
is a contraction for every $n$ sufficiently close to $e$. 
In stratified groups, using the notion of P-differentiability,
the analogous argument does not work, due to noncommutativity.
To overcome this point, we show that the mapping $F_n(h)=f(nh)$ is uniformly
biLipschitz with respect to $n$ and has constant nonvanishing 
topological degree as $n$ varies in a compact neighbourhood
of the unit element of $N$. This gives existence and uniqueness
of the implicit mapping.
We also wish to point out how inequality \eqref{intrlip} surprisingly
fits into the intrinsic notion of Lipschitz mapping given in Definition~3.1 of \cite{FSSC7}, when the ambient space is an Heisenberg group.

It is now natural to investigate the natural counterpart of the implicit function theorem, namely, the rank theorem.
To translate this theorem in terms of continuously P-differentiable mappings,
we first single out the class injective h-homomorphisms that provide the splitting of the target group, according to Proposition~\ref{GlinHmono}. 
\begin{The}[Rank Theorem]\label{Pembed}
Let $f:\Omega\lra\M$ be a continuously P-differentiable mapping,
where $\Omega$ is an open subset of $\G$.
Let $\overline x\in\Omega$ and let $Df(\overline{x}):\G\lra\M$
be the P-differential. Let us assume that $Df(\overline x)$ is an h-monomorphism
of image $H$ and let $N$ be a normal complementary subgroup.
Let $p:\M\lra\ H$ be the associated canonical projection.
Then there exist neighbourhoods $V\subset\Omega$ of $\overline{x}$ and
$W\subset H$ of $p(f(\overline{x}))$ along with mappings $\ph:W\lra N$,
$J:H\lra\G$ and $\Psi:f(V)\lra\M$ such that we have
\begin{eqnarray}\label{imagerep}
f(V)=\{h\ph(h)\mid h\in W\}\quad\mbox{and}\quad \Psi\circ f\circ J_{|V'}=I_{|V'},
\end{eqnarray}
where $I:H\hookrightarrow\M$ is the restriction of the identity mapping $\mbox{Id}_\M$
to $H\subset\M$, the open subset $V'\subset H$ is equal to $J^{-1}(V)$ and
$J=\big(p\circ Df(\overline{x})\big)^{-1}$.
Furthermore, setting $F=\exp^{-1}\circ\ph$, then there exists
$C>0$ such that for every $h,h'\in W$ we have
\begin{eqnarray}\label{lipparam}
\|F(h)-F(h')\|\leq C\, d(h,h')\,.
\end{eqnarray}
\end{The}
The proof of this result uses completely different tools with respect
to the implicit function theorem. Here the key observation is that
the projection $p:\M\lra H$ is an h-epimorphism, hence it is Lipschitz.
Then we apply our inverse mapping theorem. 

Implicit function theorem and rank theorem are the standard tools
to define differentiable manifolds. 
Analogously, in an obvious way Theorem~\ref{implth} and Theorem~\ref{Pembed}
define subsets that are expected to possess some intrinsic regularity.
These are the $(\G,\M)$-regular sets, that we distinguish into those contained
in $\G$, that are suitable level sets and those of $\M$, 
that are suitable images, see Section~\ref{GMregsets} for more details.
By estimate \eqref{holdergraph}, that immediately follows
from \eqref{intrlip}, $(\G,\M)$-regular sets of $\G$ can be
locally parametrized by $1/\iota$-H\"older mappings.
In the case of $(\H^1,\R)$-regular sets of $\H^1$, hence $\iota=2$,
B. Kirchheim and F. Serra Cassano have proved that
the embedding H\"older exponent $1/2$ cannot be improved in general.
In the same work \cite{KirSer}, the authors provide through a nontrivial
construction an example of $(\H^1,\R)$-regular set of $\H^1$
of Euclidean Hausdorff dimension $5/2$, that is clearly not
rectifiable in the Federer sense, \cite{Fed}.
Understanding the intrinsic regularity of these sets is certainly hard
and very far from being completely understood.
For instance, a fine characterization of $(\H^n,\R)$-regular sets as 
suitable 1-codimensional intrinsic graphs has been recently
established in \cite{ASCV}.

As an immediate consequence of both rank theorem and implicit function theorem, we have the following
\begin{Cor}\label{GMreggraph}
Every $(\G,\M)$-regular set is locally an intrinsic graph.
\end{Cor}
Intrinsic graphs and $(\G,\M)$-regular sets of $\G$ first appeared in
the works of Franchi, Serapioni and Serra Cassano, \cite{FSSC3}, \cite{FSSC4}
in the case $\M=\R$ and more recently in Heisenberg groups with
$\M=\R^k$, \cite{FSSC6}.
In general stratified groups, these sets have been considered in
\cite{Mag5}, where the problem of studying their metric and
topological properties has been raised.
For instance, Corollary~\ref{GMreggraph} shows that
the topological codimension of $(\G,\M)$-regular sets corresponds to the topological dimension of $\M$. 

A detailed study of $(\H^n,\R^k)$-regular sets, $1\leq k\leq n$,
has been accomplished in \cite{FSSC6}, 
where a corresponding implicit function theorem and an area-type formula have 
been established. In the same work, following the Rumin complex \cite{Rum90},
the authors lay the foundations for an intrinsic theory of currents,
where ``rectifiable'' sets correspond to both $(\H^n,\R^k)$-regular sets
and $(\R^k,\H^n)$-regular sets.
In the terminology of \cite{FSSC6}, these two class of sets
correspond to {\em low codimensional $\H$-regular surfaces}
and {\em low dimensional $\H$-regular surfaces}, respectively.
As an example of another geometry to be considered,
Theorem~\ref{intregh21} provides the complete list of all possible 
intrinsically regular sets of the 6-dimensional complexified
Heisenberg group $\H_2^1$.
To obtain this classification we need both algebraic and analytical
tools. In fact, we first find all possible factorizations
of $\H_2^1$ into an inner semidirect product of two complementary subgroups,
then we apply our implicit function theorem joined with
Theorem~\ref{legalg}.

Low dimensional $\H$-regular surfaces of $\H^n$ are nontrivial examples
of $(\G,\M)$-regular sets of $\M$, where $\G=\R^k$ and $\M=\H^n$.
More generally, Theorem~\ref{legalg} characterizes $(\R^k,\M)$-submanifolds
as $k$-dimensional $C^1$-Legendrian submanifolds of $\M$.
This theorem is a straightforward consequence of Theorem~\ref{PdifContact}.
As another consequence of Theorem~\ref{legalg}, Legendrian submanifolds are
intrinsic graphs. If we consider general couples of stratified groups,
then $(\G,\M)$-regular sets of $\M$ can be thought of as a smooth version of
``$N$-rectifiable sets'' studied in \cite{Pau2}, where $\G=N$.

From the metric viewpoint, estimate \eqref{lipparam} shows a
stronger regularity of image sets with respect to level sets,
where the ``nonlinear estimate'' \eqref{intrlip} holds.
However, in the direction of intrinsic regularity we have the following
\begin{The}[Intrinsic Blow-up]\label{blwlevelimage}
Under hypotheses of Theorem~\ref{implth}, 
we consider the set $S=f^{-1}\big(f(\overline{x})\big)$. 
Then for every $R>0$ we have
\[
D_R\cap \delta_{1/\lambda}\big((\overline{x})^{-1}S\big)
\lra D_R\cap N\quad\mbox{as}\quad\lambda\ra0^+
\]
with respect to the Hausdorff convergence of sets. In particular,
$\mbox{\rm Tan}(S,\overline{x})=N$.

Under hypotheses of Theorem~\ref{Pembed}, setting $S=f(V)$,
for every $R>0$, we have
\[
D_R\cap \delta_{1/\lambda}\big((\overline{x})^{-1}S\big)\lra D_R\cap H
\quad\mbox{as}\quad\lambda\ra0^+
\]
with respect to the Hausdorff convergence of sets. In particular,
$\mbox{\rm Tan}(S,\overline{x})=H$.
\end{The}
As a consequence of this theorem, we notice another strong difference
between level sets and image sets in general. In fact, all
homogeneous tangent cones to a $(\G,\M)$-regular set of $\M$ are clearly
h-isomophic to $\M$ and in particular have all the same
Hausdorff dimension. Furthermore, by the area formula of \cite{Mag},
their Hausdorff dimension coincides with that of $\G$ and 
we have an integral formula for their Hausdorff measure.
On the other hand, homogeneous tangent cones to a $(\G,\M)$-regular set of $\G$ are not necessarily h-isomorphic to each other, as Example~\ref{nnhisotan} shows.
By the way, Corollary~\ref{GMGTan} shows that all of the homogeneous
tangent cones to a $(\G,\M)$-regular set of $\G$ have the same Hausdorff dimension and this one equals to
\[
\mbox{$\mcH$-dim$(\G)-\mcH$-dim($\M$)}.
\]
This suggests that the Hausdorff dimension of a $(\G,\M)$-regular
set of $\G$ should coincide with that of its homogeneous tangent cones,
but this problem still claims to be investigated. 

From the previous results, it is clear that the richness of $(\G,\M)$-regular sets is connected to the richness of h-homomorphisms between
$\G$ and $\M$ along with their factorizing properties.
To advertize novel applications of both the rank theorem and
the implicit function theorem, in Section~\ref{sectfactor}
we consider some special couples of groups where all injective h-homomorphisms
are h-monomorphisms and all surjective h-homomorphisms are h-epimorphisms.
We believe that possibly many more ``geometries'' can be discovered
through a deeper algebraic investigation.

Finally, we give a concise overview of the paper.
In Section~\ref{preldef}, we recall the main definitions
that will be used throughout. Section~\ref{techlem} develops the
technical machinery of the paper, that plays a key role in the proofs
of our main results.
In Section~\ref{estimhcurve}, we establish quantitative estimates
on the P-difference quotient of horizontal curves and we apply them
to characterize continuously P-differentiable mappings.
In Section~\ref{ABSC} we characterize absolutely continuous curves
with respect to a homogeneous distance and through this characterization
we prove the Lipschitz property of continuously P-differentiable
mappings with values in a graded group.
Section~\ref{MVI} contains the proof of the mean value inequality
for continuously P-differentiable mappings.
As an immediate application, in Subsection~\ref{invmapthe} we give 
a proof of the inverse mapping theorem in stratified groups.
Section~\ref{homsubgr} recalls the notion of homogeneous subgroup,
introduces complementary subgroups and gives
simple characterizations of both h-epimorphisms and h-monomorphisms.
In Section~\ref{hquot} we show that the quotient of a graded group
by a normal homogeneous subgroup is still a graded group
and the analogous statement holds for stratified groups.
Section~\ref{MProofs} is devoted to the proof of both the implicit function 
theorem and the rank theorem.
In Section~\ref{GMregsets} we prove the everywhere existence of the
homogeneous tangent cone to $(\G,\M)$-regular sets and we show by
an example that they might not be h-isomorphic to each other in the case
of $(\G,\M)$-regular sets of $\G$.
In Section~\ref{sectfactor} we give the notions
of h-quotients and of h-embeddings.
We introduce the notion of factorizing group as a quotient
and factorizing group as a subgroup, providing corresponding examples.
In Section~\ref{SecEx}, we characterize all $C^1$ smooth Legendrian
submanifolds in graded groups as $(\R^k,\M)$-regular sets and we 
find all intrinsically regular sets of Heisenberg groups and
of the complexified Heisenberg group. 

\vskip.2truecm
{\bf Acknowledgements.}
I am grateful to Fulvio Ricci for his fruitful comments and
pleasant discussions. I thank Yu. L. Sachkov for the interesting
discussions we had on Lie groups theory when I was at ISAS of Trieste.
I thank Raul Serapioni for his subtle observations on the issue
of higher codimensional intrinsically regular level sets.
I wish to thank Alessandro Ottazzi for some inspiring discussions
on factorizations of stratified groups.
\section{Preliminaries and definitions}\label{preldef}
All Lie groups we consider in this paper are real, connected,
simply connected and finite dimensional.
A {\em graded group} is a Lie group
$\G$, whose Lie algebra $\cG$ can be written as the direct sum
of subspaces $V_i$, called layers, such that
\begin{eqnarray}\label{grading}
[V_i,V_j]\subset V_{i+j}
\end{eqnarray}
and $\cG=V_1\oplus\cdots\oplus V_\iota$.
The integer $\iota$ is the step of nilpotence of $\G$, \cite{FS}.
A graded group $\G$ is {\em stratified} if its layers satisfy
the stronger condition $[V_i,V_j]=V_{i+j}$.

The grading of $\cG$ allows us to introduce a one-parameter group
of Lie algebra automorphisms $\delta_r:\cG\lra\cG$, defined as
$\delta_r(X)=r^i$ if $X\in V_i$, where $r>0$.
These mappings are called {\em dilations}.
Taking into account that the exponential mapping $\exp:\cG\lra\G$
is a diffeomorphism for simply connected nilpotent Lie groups,
we can read dilations in the group $\G$ through the mapping $\exp$
and mantain the same notation. Recall from Theorem~2.14.3 of \cite{Vara}
that the differential of the exponential mapping is given by the followig formula
\begin{equation}\label{dexp}
d\exp\,(X)=\mbox{Id}-\sum_{n=2}^\infty\frac{(-1)^n}{n!}\,\mbox{ad}(X)^{n-1}\,,
\end{equation}
\begin{Exa}\label{versub}{\rm Let $\H^1$ denote the Heisenberg group
and let $\ch^1=\span\{X,Y,Z\}$ be its Lie algebra with $[X,Y]=Z$.
Then the ``vertical'' subgroup $\exp\big(\span\{X,Z\}\big)$ is a graded group,
but it is not stratified.}
\end{Exa} 
The group operation can be read in the algebra as follows
\begin{eqnarray}\label{absBCH}
X\gope Y=\sum_{j=1}^\iota	c_n(X,Y)
\end{eqnarray}
where  $c_1(X,Y)=X+Y$ and the addends $c_n$ are given by induction through
the Baker-Campbell-Hausdorff formula
\begin{eqnarray}\label{kBCH}
&&(n+1)\,c_{n+1}(X,Y)=\frac{1}{2}\;[X-Y,c_n(X,Y)] \\
&&+\sum_{\substack{p\geq 1\\ 2p\leq n}}K_{2p}
\sum_{\substack{k_1,\ldots,k_{2p}>0 \\
k_1+\cdots k_{2p}=n}}[c_{k_1}(X,Y),[\cdots,[c_{k_{2p}}(X,Y),X+Y],],\ldots,],
\nonumber
\end{eqnarray}
see Lemma~2.15.3 of \cite{Vara}.
Analyzing \eqref{kBCH}, one easily notices that
\begin{eqnarray}
c_n(\lambda X,\lambda Y)=\lambda^n\,c_n(X,Y)	
\end{eqnarray}
for every $X,Y\in\cG$ and $\lambda\in\R$. These formulae will be
important in the next section.

The metric structure of a graded group is given by a
continuous, left invariant distance $d:\G\times\G\lra\R$ 
such that $d(\delta_rx,\delta_ry)=r\,d(x,y)$ for every
$x,y\in\G$ and $r>0$. Every distance satisfying these properties is a {\em homogeneous distance}.

The Carnot-Carath\'eodory distance is an important example of homogeneous distance that can be defined in stratified groups, since they
satisfy the Lie bracket generating condition, see for instance \cite{Gr1}.

Notice that graded groups may not satisfy this condition.
On the other hand, according to Example~\ref{versub},
when a graded group is contained in a stratified group,
the restriction of the Carnot-Caratheodory distance
to the graded subgorup provides an example of homogeneous distance.
In general, it is possible to construct homogeneous distances in every
graded group, \cite{FSSC5}, \cite{RicPC}.

We denote by $e$ the unit element and to simplify notation
we set $d(x)=d(x,e)$. Notice that left invariance and symmetry
of $d$ imply the equality $d(x^{-1})=d(x)$. 
An open ball of center $x$ and radius $r$ with respect to
a homogeneous distance will be denoted by $B_{x,r}$.
The corresponding closed ball will be denoted by $D_{x,r}$

\subsection{h-homomorphisms and notions of differentiability}\label{homodiff}
\begin{Def}[h-homomorphism]{\rm
Let $\G$ and $\M$ be graded groups with dilations
$\delta_r^\G$ and $\delta_r^\M$, respectively.
We say that a group homomorphism $L:\G\lra\M$
such that $L(\delta_r^\G x)=\delta_r^\M L(x)$ for every $x\in\G$ and $r>0$
is a {\em homogeneous homomorphism},
in short {\em h-homomorphism}. 
Invertible h-homomorphisms will be called {\em h-isomorphisms}.
}\end{Def}
\begin{Rem}{\rm
Analogous terminology will be used for the corresponding
Lie algebra homomorphisms of graded algebras that commute
with dilations.}
\end{Rem}
\begin{Rem}\label{disthomo}
{\rm Concerning injective and surjective h-homomorphisms, we will
use the classical terminology of h-epimorphism and h-monomorphisms
to indicate special classes of surjective and injective h-homomorphisms. 
In fact, a surjective linear mapping of vector spaces 
is characterized by the existence of a right inverse that is also
linear. Analogously, injective linear mappings are characterized by the
existence of a linear left inverse mapping.
The analogous characterization for either surjective or
injective algebra homomorphisms does not work, as we will see in the
next example.}
\end{Rem}
\begin{Exa}\label{exahr2}{\rm
Let $l:\ch^1\lra\R^2$ be a surjective h-homomorphism.
Clearly $l(\cv)=\R^2$ and $l(\cz)=\{0\}$, where $\ch^1=\cv\oplus\cz$.
We show that there is no right inverse that is also an
h-homomorphism. Assume by contradiction that there exists an
h-homomorphism $\tau:\R^2\lra\ch^1$ that is a right inverse.
Then the property $l\circ\tau=\Id_{\R^2}$ and the fact that
$\tau$ is an h-homomorphism imply that $\tau(\R^2)$ is a 2-dimensional
homogeneous subalgebra of $\ch^1$.
Clearly $\tau(\R^2)$ cannot intersect $\ker l=\cz$, but this conflicts
with Example~\ref{nnexcmpl}, where we show that any 2-dimensional homogeneous
subalgebra of $\ch^1$ contains $\cz$.
}\end{Exa}
As a consequence of the previous example,
in the category of graded algebras and
h-homomorphisms, requiring the existence of a right inverse homomorphism
is a stronger condition than surjectivity.
This motivates the following
\begin{Def}\label{hepihmono}{\rm
We say that an h-homomorphism is an {\em h-epimorphism} if it has
a right inverse that is also an h-homomorphism. We say
that an h-homomorphism is an {\em h-monomorphism} if it has a left
inverse that is also an h-homomorphism.}
\end{Def}
In Subsection~\ref{hepimono}, we will see how either h-epimorphisms
or h-monomorphisms can be characterized by their property of factorizing
either the domain or the codomain.
\begin{Def}[P-differentiability]{\rm
Let $\G$ and $\M$ be graded groups with homogeneous distances
$d$ and $\rho$, respectively. Let $\Omega$ be an open subset of $\G$
and consider $f:\Omega\lra\M$.
We say that $f$ is {\em P-differentiable} at $x\in\Omega$ if
there exists an h-homomorphism $L:\G\lra\M$ such that
$$
\frac{\rho\left(f(x)^{-1}f(xh),L(h)\right)}{d(h)}\lra0
\quad\mbox{as}\quad h\ra e\,.
$$
The h-homomorphism $L$ satisfying this limit is unique
and it is called {\em P-differential} of $f$ at $x$.
We denote $L$ by $Df(x)$. When we read the P-differential
between Lie algebras, we will denote it by $df(x)$.
}\end{Def}
\begin{Def}[Horizontal differentiability]\label{hdiff}{\rm
Let $\Omega\subset\G$ be an open set and let $M$ be a smooth manifold.
We consider a mapping $f:\Omega\lra M$ and $x\in\Omega$.
We say that $f$ is {\em horizontally differentiable}
at $x$, in short {\em h-differentiable} at $x$, 
if there exists a neighbourhood $U$ of the origin in
the first layer $V_1$ such that the restriction  
\[
U\ni X\lra f(x\exp X)\in M
\]
is differentiable at the origin. The differential $L:H_x\G\lra T_{f(x)}M$
will be denoted by $d_Hf(x)$.
We say that $f$ is {\em continuously} h-differentiable
in the case $x\lra d_Hf(x)$ is a continuous mapping.
If $M$ is a graded group $\M$, then we also use the notation
$D_Hf(x):\exp V_1\lra\M$ to denote the h-differential between
the corresponding groups.
}\end{Def}
\begin{Rem}\label{Pdifhdif}{\rm
We notice that P-differentiability of $\R^k$-valued mappings
on an open subset $\Omega\subset\G$ implies h-differentiability.
It suffices to restrict P-differentiability to horizontal directions
and to observe that $d(\exp X)/\|X\|$ is
bounded away from zero and from above with constants independent
from $X$ as it varies in $V_1$. This yields
\begin{eqnarray}\label{dhfP}
D_Hf(x)(\exp X)=Df(x)(\exp X)\qquad \mbox{for every}\qquad X\in V_1.
\end{eqnarray}
Then h-differentiability in general is a weaker notion.
On the other hand, the little regularity of h-differentiable mappings
suffices to introduce contact mappings.}
\end{Rem}
\begin{Def}[Contact mapping]{\rm
Let $f:\Omega\lra\M$ be an h-differentiable mapping.
We say that $f$ is a {\em contact mapping} if the inclusion
$d_Hf(x)(H_x\G)\subset H_{f(x)}\M$ holds for every $x\in\Omega$.
}\end{Def}
\begin{Rem}{\rm
Notice that $\R^k$-valued continuously h-differentiable mappings
are automatically contact. This family of mappings have been
already considered in \cite{FSSC4} as mappings of ${\bf C}^1_\G$.}
\end{Rem}
\begin{Pro}[Chain rule]\label{chain}
Let $f:\Omega\lra\U$ be P-differentiable at $x\in\Omega$
and let $g:\Upsilon\lra\M$ be P-differentiable at $f(x)\in\Upsilon$,
where $\Omega\subset\G$ and $\Upsilon\subset\U$ are open subsets
and $f(\Omega)\subset\Upsilon$.
Then $g\circ f:\Omega\lra\M$ is P-differentiable at $x$ and
$D(g\circ f)(x)=Dg\big(f(x)\big)\circ Df(x)$.
\end{Pro}
The proof of this result is straightforward, see for instance
Proposition~3.2.5 of \cite{MagPhD}.
\section{Technical lemmata}\label{techlem}

Throughout this section, we denote by $\G$ a graded group,
equipped with a homogeneous distance $d$.
Its Lie algebra $\cG$ of layers $V_i$ has step $\iota.$
On $\cG$, seen as a finite dimensional real vector space
we fix a norm $\|\cdot\|$.
Bilinearity of brackets yields a constant $\beta>0$, such that
for every $X,Y\in\cG$ we have
\begin{equation}\label{lcst}
\|[X,Y]\|\leq\beta\, \|X\|\;\|Y\|.
\end{equation}
\begin{Lem}
Let $\nu>0$ and let $n=2,\ldots,\iota$. Then there exists
a constant $\alpha_n(\nu)$ only depending on $n$ and $\nu$ such that
\begin{eqnarray}\label{bilestim}
\|c_n(X,Y)\|\leq \alpha_n(\nu)\;\|[X,Y]\|
\end{eqnarray}
whenever $\|X\|,\|Y\|\leq\nu$.
\end{Lem}
{\sc Proof.}
Our statement is trivial for $n=2$, being $c_2(X,Y)=[X,Y]/2$.
Assume that it is true for every $j=2,\ldots,n$, with $n\geq2$.
It suffices to observe that $[c_{k_{2p}}(X,Y),X+Y]\neq0$ in \eqref{kBCH} implies $k_{2p}>1$, then inductive hypothesis yields
$$
\|c_{k_{2p}}(X,Y)\|\leq \alpha_{k_{2p}}(\nu)\; \|[X,Y]\|\,.
$$
Using this estimate in \eqref{kBCH} and observing that
$\|c_{k_i}(X,Y)\|\leq2\nu$, whenever $k_i=1$,
the thesis follows. $\Box$
\begin{Lem}
Let $c_n(X,Y)$ be as in \eqref{absBCH}.
Then for each $n=2,\ldots,\iota$ there exists
a set of real numbers 
$\big\{e_{n,\alpha}\mid \alpha\in \{1,2\}^{n-1}\big\}$ only
depending on $\cG$ such that for every $A_1,A_2\in\cG$ we have
\begin{eqnarray}\label{mlinc}
c_n(A_1,A_2)=\sum_{\alpha\in\{1,2\}^{n-1}}e_{n,\alpha}\;
L_n(A_\alpha,A_1+A_2)\,,
\end{eqnarray}
where $A_\alpha=(A_{\alpha_1},\ldots,A_{\alpha_{n-1}})$,
$L_1=\mbox{Id}_\cG$ 
and for $n\geq2$ the $n$-linear mapping $L_n:\cG^n\lra\cG$ is defined by
\begin{equation}\label{Ln}
L_n(X_1,X_2,\ldots,X_n)=[X_1,[X_2,[\ldots,[X_{n-1},X_n]]\ldots]\,.
\end{equation}
\end{Lem}
{\sc Proof.}
Let $L_1:\cG\lra\cG$ be the identity mapping and let
$k_1,\ldots,k_p$ be positive integers with
$n=k_1+k_2\cdots+k_p$ where $p\in\N$.
Iterating Jacoby identity of the Lie product, it is not
difficult to check that
\begin{eqnarray}\label{nnassp}
&&[L_{k_1}(X^1_1,\ldots,X^1_{k_1}),[L_{k_2},(X^2_1,\ldots,X^2_{k_2}),[\cdots,\\
&&[L_{k_{p-1}}(X^{p-1}_1,\ldots,X^{p-1}_{k_{p-1}}), 
L_{k_p},(X^p_1,\ldots,X^p_{k_p})]]\ldots] \nonumber \\
&&=\sum_{\sigma_1\in S_{k_1},\ldots,\sigma_{p-1}\in S_{k_{p-1}}}
s_{\sigma_1,\ldots,\sigma_{p-1}}
\;\;L_n\big(X_{\sigma_1}^1,X^2_{\sigma_2},\ldots,X^{p-1}_{\sigma_1},
X^p\big)\,,\nonumber
\end{eqnarray}
where $S_j$ is the set of all permutations on $j$ elements,
$s_{\sigma_1,\ldots,\sigma_{p-1}}\in\{-1,0,1\}$,
$$
X^p=(X^p_1,\ldots,X^p_{k_p})\quad\mbox{ and}\quad
X^j_{\sigma_j}=(X^j_{\sigma_j(1)},X^j_{\sigma_j(2)},\ldots,
X^j_{\sigma_j(k_j)})\,.
$$
Our statement can be proved by induction.
It is clearly true for $n=2$ taking $e_1=1/2$ and $e_2=0$,
due to the formula $c_2(A_1,A_2)=[A_1,A_2]/2$.
Let us assume that \eqref{mlinc} holds for all $c_j(A_1,A_2)$
with $j\leq n$. By recurrence equation \eqref{kBCH}, we have
\begin{eqnarray}
&&(n+1)\,c_{n+1}(A_1,A_2)=\frac{1}{2}\sum_{\alpha\in\{1,2\}^{n-1}}
e_{n,\alpha}\;\left[A_1-A_2,
L_n(A_\alpha,A_1+A_2)\right] \\
&&+\sum_{\substack{p\geq 1\\ 2p\leq n}}K_{2p}
\sum_{\substack{k_1,\ldots,k_{2p}>0 \\
k_1+\cdots k_{2p}=n}}
\sum_{\substack{\;\alpha_i\in\{1,2\}^{k_i-1}\\i=1,\ldots,2p}}
e_{k_1,\alpha_1}\;e_{k_2,\alpha_2}\;\cdots\; e_{k_{2p},\alpha_{2p}}\nonumber \\
&&[L_{k_1}(A_{\alpha_1},A_1+A_2),[L_{k_2}(A_{\alpha_2},A_1+A_2,),
[\cdots,[L_{k_{2p}}(A_{\alpha_{2p}},A_1+A_2),A_1+A_2]]\ldots]
\nonumber
\end{eqnarray}
Then applying \eqref{nnassp} we have proved that $c_{n+1}(A_1,A_2)$
can be represented as a linear combination of terms
$L_{n+1}(A_\alpha,A_1+A_2)$ where $\alpha\in\{1,2\}^n$.
This concludes our proof. $\Box$
\vskip.2truecm
Following the notation of \cite{Pan2}, Section~4.5, in the next definition
we introduce the iterated Lie bracket.
\begin{Def}{\rm
Let $X,Y\in\cG$. The $k$-th bracket is defined by
\begin{eqnarray}
[X,Y]_k=\underbrace{[X,[X,[\cdots,[X}_{\mbox{\tiny $k$ times}},Y],],\ldots,]
\quad\mbox{and}\quad[X,Y]_0=Y\,.
\end{eqnarray}
}\end{Def}
\begin{Lem}\label{keylemma}
Let $c_n(X,Y)$ be as in \eqref{absBCH}, where $n=2,\ldots,\iota$. Then we have
\begin{eqnarray}\label{cnrest}
c_n(X,Y)=\frac{(-1)^{n-1}}{n!}\;
\left[\frac{Y-X}{2},X+Y\right]_{n-1}+R_n(X,Y)
\end{eqnarray}
and for every  $\nu>0$ there exists a positive
nondecreasing function $C(n,\cdot)$ such that
\begin{eqnarray}\label{estimsmall}
\|R_n(X,Y)\|\leq C(n,\nu)\, \|X+Y\|^3
\end{eqnarray}
whenever $\|X\|,\|Y\|\leq\nu$.
\end{Lem}
{\sc Proof.}
Our statement is trivial for $n=2$, being $R_2(X,Y)$ vanishing
for any $X,Y\in\cG$. Let us consider $n\geq3$.
We apply formula \eqref{mlinc} to \eqref{kBCH}, getting
\begin{eqnarray}\label{cnrestpar}
c_n(X,Y)=-\frac{1}{n}\;\left[\frac{Y-X}{2},c_{n-1}(X,Y)\right]+E_n(X,Y)
\end{eqnarray}
where for every $A_1,A_2\in\cG$ we have
\begin{eqnarray*}
&&E_n(A_1,A_2)=\sum_{\substack{p\geq 1\\ 2p\leq n-1}}K_{2p}
\sum_{\substack{k_1,\ldots,k_{2p}>0 \\k_1+\cdots k_{2p}=n-1}}
\sum_{\substack{\;\alpha_i\in\{1,2\}^{k_i-1}\\i=1,\ldots,2p}}
e_{k_1,\alpha_1}\;e_{k_2,\alpha_2}\;\cdots\; e_{k_{2p},\alpha_{2p}}\\
&&[L_{k_1}(A_{\alpha_1},A_1+A_2),[L_{k_2}(A_{\alpha_2},A_1+A_2,), 
[\cdots,[L_{k_{2p}}(A_{\alpha_{2p}},A_1+A_2),A_1+A_2]]\ldots]\nonumber\,.
\end{eqnarray*}
As a consequence, there exist constants $\widetilde C_p>0$ such that
\begin{eqnarray}\label{enest}
&&\|E_n(X,Y)\|\leq\sum_{1\leq p\leq(n-1)/2} \widetilde C_p
\;\; \nu^{n-2p-1}\;\|X+Y\|^{2p+1}\\
&&\leq\sum_{1\leq p\leq(n-1)/2} 4^p\;\widetilde C_p
\;\; \nu^{n-3}\;\|X+Y\|^3\,.\nonumber
\end{eqnarray}
Iterating \eqref{cnrestpar}, we get
\begin{eqnarray}
&& c_n(X,Y)=\frac{(-1)^2}{n(n-1)}\left[\frac{Y-X}{2},c_{n-2}(X,Y)\right]_2
+E_n(X,Y) \\
&& -\frac{1}{n}\left[\frac{Y-X}{2},E_{n-1}(X,Y)\right]
=\frac{(-1)^{n-1}}{n!}\left[\frac{Y-X}{2},c_1(X,Y)\right]_{n-1}
+E_n(X,Y)\nonumber \\
&&+\sum_{k=1}^{n-2}\frac{(-1)^k}{n(n-1)\cdots(n-k+1)}
\left[\frac{Y-X}{2},E_{n-k}(X,Y)\right]_k\nonumber
\end{eqnarray}
therefore \eqref{enest} concludes our proof. $\Box$
\begin{Rem}{\rm
As an immediate consequence of Lemma~\ref{keylemma},
we achieve the estimate
\begin{eqnarray}\label{leftinveucl}
\|(-\xi)\gope\eta\|\leq C(\nu) \|\xi-\eta\|
\end{eqnarray}
for every $\xi,\eta\in\cG$ satisfying $\|\xi\|,\|\eta\|\leq\nu$.
To see this, it suffices to apply \eqref{cnrest} and \eqref{estimsmall}
to the Baker-Campbell-Hausdorff formula
$(-\xi)\gope\eta=-\xi+\eta+\sum_{n=2}^\iota c_n(-\xi,\eta)$.
}\end{Rem}
As an immediate extension of inequality (+) at p.13 of \cite{Pan2},
see also \cite{Pan1}, we have the following lemma.
\begin{Lem}\label{estimi}
Let $\pi^i:\cG\lra V_i\oplus\cdots\oplus V_\iota$
be the natural projection and let $U$ be
a bounded open neighbourhood of the unit element
$e\in\G$. Then there exists a constant $K_U>0$, depending on $U$,
such that
\begin{eqnarray}\label{estimil}
\|\pi^i\left(\exp^{-1}(x)\right)\|\leq K_U\;d(x)^i	
\end{eqnarray}
holds for every $x\in U$.
\end{Lem}
{\sc Proof.}
We first define $\cS=\{v\in\cG\mid d(\exp v)=1\}$.
Let us fix $v\in\cS$ and $s>0$ satisfying the
condition $s\leq M$, with $M=\max_{x\in U} d(x)$.
We have
$$
\|\pi^i(\delta_sv)\|=s^i\;\Big\|\sum_{j=i}^\iota s^{j-i}\,\pi_j(v)\Big\|
\leq C s^i\,,
$$
where $C=\sum_{j=i}^\iota M^{j-i}\,\max_{u\in\cS}\|\pi_j(u)\|$\,.
This concludes the proof. $\Box$
\begin{Rem}{\rm
From Proposition~1.5 of \cite{FS}, it is easy to check that
there exists a constant $\kappa(\nu)$ such that whenever $\|\xi\|\leq\nu$
there holds 
\begin{eqnarray}\label{rhonormiota}
d(\exp\xi)\leq\kappa(\nu)\,\|\xi\|^{1/\iota}\,.
\end{eqnarray} 
}\end{Rem}
\begin{Rem}{\rm 
As a byproduct of \eqref{leftinveucl} and \eqref{rhonormiota} we obtain
the well known estimate
\begin{eqnarray}\label{rhoestnorm}
d\big(\exp\xi,\exp\eta\big)\leq C(\nu)\; \|\xi-\eta\|^{1/\iota}
\end{eqnarray}
for every $\xi,\eta\in\cG$ satisfying $\|\xi\|,\|\eta\|\leq\nu$,
see \cite{RotStein} for its proof in the more general setting
of vector fields satisfying the H\"ormander condition.
}\end{Rem}
\begin{Rem}{\rm
Notice that for every homogeneous distance $d$, there
exists a constant $C_{d,\|\cdot\|}$ such that
\[
C_{d,\|\cdot\|}^{-1}\;\sum_{j=1}^\iota\|\pi_j(\xi)\|^{1/j}\leq
d(\exp\xi)\leq C_{d,\|\cdot\|}\;\sum_{j=1}^\iota\|\pi_j(\xi)\|^{1/j}\,.
\]
In particular, we have
\begin{equation}\label{lip1}
d(\exp\xi,\exp\eta)\geq C_{d,\|\cdot\|}^{-1}\;\|\pi_1(\xi-\eta)\|\,.
\end{equation}
for every $\xi,\eta\in\cG$.
}\end{Rem}
\begin{Lem}\label{lemconjest}
Let $x,y\in\G$ and let $\nu>0$ be such that
$d(x),d(y)\leq\nu$.
Then there exists a constant $C(\nu)$ only depending on
$\G$ and on $\nu$ such that
\begin{eqnarray}\label{conjest1}
d(y^{-1}xy)\leq C(\nu)\;\|\exp^{-1}(x)\|^{1/\iota}\,.
\end{eqnarray}
\end{Lem}
{\sc Proof.}
Let us fix $x=\exp\xi$ and $y=\exp h$. By Lemma~\ref{estimi}
we find a constant $k_\nu>0$ such that $\|\xi\|,\|h\|\leq k_\nu$.
The Baker-Campbell-Hausdorff fomuls yields
\begin{eqnarray}\label{bchassoc}
(-h)\gope\big(\xi\gope h\big)
=\xi+\sum_{n=2}^\iota c_n(\xi,h)+\sum_{n=2}^\iota c_n(-h,\xi\gope h).
\end{eqnarray}
The same formula along with \eqref{bilestim},
also gives the estimate
$$
\|\xi\gope h\|\leq k_\nu\,\Big(2+\sum_{n=2}^\iota
\alpha_n(k_\nu)\,\beta\,k_\nu\Big)
=A(\nu).
$$
The bilinear estimate \eqref{bilestim} also yields
\begin{eqnarray}\label{iterest}
\|c_n(-h,\xi\gope h)\|\leq\alpha_n\big(A(\nu)\big)\,
\|[h,\xi\gope h]\|.
\end{eqnarray}
Observing that 
$$
[h,\xi\gope h]=\big[h,\xi+\sum_{n=2}^\iota c_n(\xi,h)\big]
$$
we achieve the estimate
$$
\|[h,\xi\gope h]\|\leq\beta\,\|h\|\,\|\xi\|\,
\big(1+\sum_{n=2}^\iota\alpha_n(k_\nu)\,\|h\|\big).
$$
Joining this estimate with \eqref{iterest} we achieve
\begin{eqnarray}\label{iterest1}
\|c_n(-h,\xi\gope h)\|\leq
\beta\,\alpha_n\big(A(\nu)\big)\,\Big(1+\sum_{n=2}^\iota\nu\,\alpha_n(\nu)\Big)\,\|h\|\,\|\xi\|
\end{eqnarray}
Thus, formula \eqref{bchassoc}, estimates \eqref{bilestim} and
\eqref{iterest1} give a constant $B(\nu)>0$ such that
\begin{eqnarray}\label{fingope}
\|(-h)\gope\xi\gope h\|\leq \|\xi\|\;\big(1+B(\nu)\|h\|\big).
\end{eqnarray}
The previous inequality yields
$$
\|(-h)\gope\xi\gope h\|\leq k_\nu\,\big(1+k_\nu\;B(\nu)\big),
$$
then \eqref{rhonormiota} gives a constant $B_1(\nu)>0$ such that
$$
d(y^{-1}xy)\leq B_1(\nu)\;\|(-h)\gope\xi\gope h\|^{1/\iota}
$$
Thus, applying \eqref{fingope} we achieve a constant $B_2(\nu)>0$
such that 
$$
d(y^{-1}xy)\leq B_2(\nu)\;\|\xi\|^{1/\iota}\,.
$$
This finishes the proof. $\Box$
\begin{Rem}{\rm
The previous lemma also provides another variant of its
estimate. In fact, the condition $d(x)\leq\nu$
implies that $\|\exp^{-1}(x)\|$ is bounded by $d(x)$
up to a factor only depending on $\nu$,
due to Lemma~\ref{estimi}.
As a result, the assumption $d(x),d(y)\leq\nu$
gives a constant $C(\nu)$ such that the following estimate holds
\begin{eqnarray}\label{conjest}
d(y^{-1}xy)\leq C(\nu)\;d(x)^{1/\iota}.
\end{eqnarray}
}\end{Rem}
\begin{Lem}\label{lestimprod}
Let $N$ be a positive integer and let $A_j,B_j\in\G$
with $j=1,\ldots,N$. Let $\nu>0$ be such that 
$d(B_jB_{j+1}\cdots B_N)\leq\nu$ and $d(A_j,B_j)\leq\nu$
for every $j=1,\ldots,N$. Then there exists $K_\nu>0$ such that
\begin{eqnarray}\label{estimprod}
d\big(A_1A_2\cdots A_N,B_1B_2\cdots B_N\big)\leq K_\nu\;
\sum_{j=1}^Nd(A_j,B_j)^{1/\iota}.
\end{eqnarray}
\end{Lem}
{\sc Proof.}
We define
\[
\hat{B}_j=B_jB_{j+1}\cdots B_N,\qquad \hat{A}_j=A_jA_{j+1}\cdots A_N
\]
and use left invariance of $d$ to obtain
\[
d\big(\hat{A}_1,\hat{B}_1\big)\leq d(A_N,B_N)
+\sum_{j=1}^{N-1}d\big(\hat{B}_{j+1}A_j^{-1}B_j\hat{B}_{j+1}\big)\,.
\]
By our hypothesis and \eqref{conjest} we get
\begin{eqnarray*}
&&d\big(\hat{A}_1,\hat{B}_1\big)\leq d(A_N,B_N)
+C(\nu)\;\sum_{j=1}^{N-1}d\big(A_j,B_j\big)^{1/\iota}
\leq d(A_N,B_N)^{1/\iota}\nu^{1-1/\iota}\\
&&+C(\nu)\;\sum_{j=1}^{N-1}d\big(A_j,B_j\big)^{1/\iota}
\leq\max\{C(\nu),\nu^{1-1/\iota}\}\sum_{j=1}^Nd(A_j,B_j)^{1/\iota}.
\end{eqnarray*}
This shows the validity of \eqref{estimprod}. $\Box$
\begin{Lem}\label{bchdifflem}
Let $X,Y,D_1,D_2\in\cG$ and let
$\nu\geq\max\{\|X\|,\|Y\|,\|D_1\|,\|D_2\|\}$. Then there exists
a constant $\kappa_n>0$ such that
\begin{eqnarray}\label{bchdiffprod}
\|c_n(X+D_1,Y+D_2)-c_n(X,Y)\|\leq \gamma_n\; \nu^{n-1}\; \max\{\|D_1\|,\|D_2\|\}
\end{eqnarray}
for every $n=2,\ldots,\iota$.
\end{Lem}
{\sc Proof.}
We set $X=A_1^1$, $D_1=A_1^2$, $Y=A_2^1$ and $D_2=A_2^2$ and
apply \eqref{mlinc}.
Taking into account that $L_n$ is a multilinear mapping defined
in \eqref{Ln}, we obtain
\begin{eqnarray*}
&&c_n(A_1^1+A_1^2,A_2^1+A_2^2)\\
&&=\sum_{\alpha\in\{1,2\}^{n-1}}e_{n,\alpha}\;
L_n(A^1_{\alpha_1}+A^2_{\alpha_1},\ldots,A^1_{\alpha_{n-1}}+A^2_{\alpha_{n-1}},A_1^1+A_1^2+A_2^1+A_2^2) \\
&&=\sum_{\alpha\in\{1,2\}^{n-1}}e_{n,\alpha}
\sum_{\gamma\in\{1,2\}^{n-1}}
L_n(A^{\gamma_1}_{\alpha_1},\ldots,A^{\gamma_{n-1}}_{\alpha_{n-1}},A_1^1+A_1^2+A_2^1+A_2^2)\\
&&=\sum_{\alpha\in\{1,2\}^{n-1}}e_{n,\alpha}
\sum_{\gamma\in\{1,2\}^{n-1}}
L_n(A^{\gamma_1}_{\alpha_1},\ldots,A^{\gamma_{n-1}}
_{\alpha_{n-1}},A_1^1+A_2^1)\\
&&+\sum_{i=1,2}\sum_{\alpha\in\{1,2\}^{n-1}}e_{n,\alpha}
\sum_{\gamma\in\{1,2\}^{n-1}}
L_n(A^{\gamma_1}_{\alpha_1},\ldots,A^{\gamma_{n-1}}_{\alpha_{n-1}},A_i^2)\\
&&=c_n(A_1^1,A_2^1)+\sum_{\alpha\in\{1,2\}^{n-1}}e_{n,\alpha}
\sum_{\substack{\gamma\in\{1,2\}^{n-1}\\\gamma\neq(1,1,\ldots,1)}}
L_n(A^{\gamma_1}_{\alpha_1},\ldots,A^{\gamma_{n-1}}
_{\alpha_{n-1}},A_1^1+A_2^1)\\
&&+\sum_{i=1,2}\sum_{\alpha\in\{1,2\}^{n-1}}e_{n,\alpha}
\sum_{\gamma\in\{1,2\}^{n-1}}
L_n(A^{\gamma_1}_{\alpha_1},\ldots,A^{\gamma_{n-1}}_{\alpha_{n-1}},A_i^2)
\end{eqnarray*}
where the equality in the line before last we have applied
again formula \eqref{mlinc}.
Taking into account the multilinear estimate
\[
\|L_n(X_1,X_2,\ldots,X_n)\|\leq\beta^n\;\prod_{j=1}^n\|X_j\|,
\]
where $\beta$ is given in \eqref{lcst}.
As a consequence, a short calculation yields
\[
\|c_n(A_1^1+A_1^2,A_2^1+A_2^2)-c_n(A_1^1,A_2^1)\|\leq
2^{n+1}\,\beta^n\,\nu^{n-1}\,\max\{\|A_1^2\|,\|A_2^2\|\}\,
\sum_{\alpha\in\{1,2\}^{n-1}}|e_{n,\alpha}|\,,
\]
leading us to the conclusion. $\Box$
%
%
%
%
%
%
%
%
%
%
%
%
%
%
%
%
\section{P-differentiability of curves in graded groups}\label{estimhcurve}
In this section we study differentiability properties of curves $\Gamma:[a,b]\lra\M$, where $\M$ is a graded group.
We denote by $\gamma$ the corresponding curve
$\exp^{-1}\circ\Gamma$ with values in the Lie algebra $\cM$ of $\M$.
The components of $\gamma$ taking values in the layers $W_i$
are denoted by $\gamma_i$, hence $\gamma=\sum_{i=1}^\upsilon\gamma_i$.
We have $\gamma_i=\pi_i\circ \gamma$, where 
$\pi_i:\cM\lra W_i$ is the canonical projection
on layers of degree $i$ and $\cM=W_1\oplus\cdots\oplus W_\upsilon$.
We denote by $\|\cdot\|$ a fixed norm on $\cM$.
The {\em horizontal subspace} of $\M$ at $x\in\M$ is
defined as follows
\[
H_x\M=\{U(x)\mid U\in W_1\}.
\]
\begin{Def}{\rm
We say that a curve $\Gamma:[a,b]\lra\M$ is {\em horizontal} if
$\gamma:[a,b]\lra\cM$ is absolutely continuous and for a.e. 
differentiability point $t\in[a,b]$ the inclusion $\dot\gamma(t)\in H_{\gamma(t)}\M$ holds.
}\end{Def}
\noindent
Identifying any $T_x\M$ with $\cM$ and applying formula \eqref{dexp}
for the differential of the exponential mapping,
we have
\begin{equation}\label{dexpgamma}
\dot\Gamma(t)=
\dot\gamma(t)-\sum_{n=2}^\upsilon\frac{(-1)^n}{n!}\,\mbox{ad}(\gamma(t))^{n-1}(\dot\gamma(t))=
\dot\gamma(t)-\sum_{n=2}^\upsilon\frac{(-1)^n}{n!}\,
[\gamma(t),\dot\gamma(t)]_{n-1}\,.
\end{equation}
Then $\dot\Gamma(t)\in W_1$ if and only if
$$
\pi_i\Big(\dot\gamma(t)-\sum_{n=2}^\upsilon\frac{(-1)^n}{n!}\,
[\gamma(t),\dot\gamma(t)]_{n-1}\Big)=0
$$
for all $i\geq2$. This immediately proves the following
\begin{Pro}\label{horprop}
Let $\gamma:[a,b]\lra\cM$ be an absolutely continuous curve.
Then $\Gamma$ is a horizontal curve if and only if the differential equation
\begin{eqnarray}\label{eqhdiffcurve}
	\dot\gamma_i(t)=\sum_{n=2}^\upsilon\frac{(-1)^n}{n!}\;
	\pi_i\left([\gamma(t),\dot\gamma(t)]_{n-1}\right)
	\end{eqnarray}
a.e. holds for each $i=2,\ldots,\upsilon$.
\end{Pro}
\begin{Def}\label{mcw}{\rm
Let $\gamma:[a,b]\lra\cM$ be a locally summable
curve and let $\lambda\neq0$.
The {\em sup-average} of $\gamma$ at $t$ is defined as follows
\begin{eqnarray}\label{cMt}
\cA_t^\lambda(\gamma)=\left\{\begin{array}{ll}
\displaystyle\sup_{0\leq\tau\leq \lambda}
\;\;\medint_t^{t+\tau}\|\gamma(l)\|\,dl &  \mbox{if $\;\lambda>0$} \\
\displaystyle\sup_{\lambda\leq\tau\leq 0}
\;\;\medint_{t+\tau}^t\|\gamma(l)\|\,dl &  \mbox{if $\;\lambda<0$}
\end{array}
\right.\,.
\end{eqnarray}
}\end{Def}
Notice that for some $t$ the sup-average $\cA_t^\lambda(\gamma)$
takes values in $[0,+\infty]$, as $\gamma$ is not necessarily
bounded in a neighbourhood of $t$.
\begin{Lem}\label{theoestim}
There exist strictly increasing functions
$\Upsilon_i:[0,+\infty[\lra[0,+\infty[$, with $i=2,\ldots,\upsilon$,
which only depend on $\cM$, are infinitesimal at zero
and satisfy the following properties.
For any horizontal curve $\Gamma:[-\alpha,\alpha]\lra\M$
such that $\Gamma(0)=e$, the estimates  
\begin{eqnarray}\label{estimhor}
\left|\int_0^\lambda\|\dot\gamma_i(t)\|\,dt\right|\leq
\Upsilon_i(L)\,
\cA_0^\lambda(\dot\gamma_1-X)\,|\lambda|^i
\end{eqnarray}
hold, where
$L=\max
\left\{\cA_0^{-\alpha}(\dot\gamma_1),\cA_0^\alpha(\dot\gamma_1)\right\}$,
$X\in\cM$, $\|X\|\leq L$
and the function $\cA_0^\lambda(\cdot)$ has been introduced
in Definition~\ref{mcw}.
\end{Lem}
{\sc Proof.}
First of all, we assume that $\cA_0^\lambda(\dot\gamma_1)<+\infty$,
otherwise the proof becomes trivial.
Due to Proposition~\ref{horprop}, the differential equations
\eqref{eqhdiffcurve} hold.
Applying \eqref{eqhdiffcurve} for $i=2$, we achieve
\begin{equation}\label{gmdot2}
\dot\gamma_2(t)=\frac{t}{2}\;\Big[\frac{\gamma_1(t)}{t},\dot\gamma_1(t)\Big]=
\frac{t}{2}
\left(\left[\medint_0^t\dot\gamma_1(\tau)\,d\tau,X\right]
+\left[\medint_0^t\dot\gamma_1(\tau)\,d\tau,
\dot\gamma_1(t)-X\right]\right)\,.
\end{equation}
For every $0<s<\alpha$, it follows that
\begin{eqnarray*}
&&\!\!\!\int_0^s\|\dot\gamma_2(t)\|\,dt 
\leq\frac{s}{2}\int_0^s\left\|\left[\medint_0^t\dot\gamma_1(\tau)\,d\tau,X\right]\right\|\,dt
+\frac{s}{2}\int_0^s\left\|\left[\medint_0^t\dot\gamma_1(\tau)\,
d\tau,\dot\gamma_1(t)-X\right]\right\|\,dt \\
&&\leq\frac{\beta Ls}{2}\left(\int_0^s\cA^t_0(\dot\gamma_1-X)\,dt
+s\,\cA_0^s(\dot\gamma_1-X)\right)\,.
\end{eqnarray*}
This shows that
\begin{eqnarray}\label{os2}
\int_0^s\|\dot\gamma_2(t)\|\,dt\leq
\beta\,L\,s^2\,\cA_0^s(\dot\gamma_1-X).
\end{eqnarray}
From the grading property \eqref{grading},
equations \eqref{eqhdiffcurve} can be stated as follows
\begin{eqnarray}\label{confshape}
\dot\gamma_j(t)=\sum_{n=2}^\upsilon\frac{(-1)^n}{n!}
\sum_{\substack{1\leq l_1,\ldots,l_n\leq\upsilon\\ l_1+\cdots +l_n=j}}
[\gamma_{l_1}(t),[\gamma_{l_2}(t),[\cdots,[\gamma_{l_{n-1}}(t),
\dot\gamma_{l_n}(t)],],\ldots,]
\end{eqnarray}
for every $j=2,\ldots,\upsilon$.
We will proceed by induction,
assuming that for every $j=2,\ldots,i$ and every $\upsilon\geq i$,
there exists a positive constant $\kappa_j(\beta)$ such that
\begin{eqnarray}\label{indsj}
\int_0^s\|\dot\gamma_j(t)\|\,dt\leq\kappa_j(\beta)\,F_j(L)
\;\cA^s_0(\dot\gamma_1-X)\,s^j,
\end{eqnarray}
where $F_j:[0,+\infty[\lra[0,+\infty[$ is a strictly increasing
function, that is infinitesimal at zero.
We have already proved this statement for $i=2$.
Using \eqref{confshape} for $j=i+1$ and applying the
inductive hypothesis \eqref{indsj}, we get
\begin{eqnarray*}
&&\int_0^s\|\dot\gamma_{i+1}(t)\|\,dt\leq\sum_{n=2}^\upsilon\frac{\beta^n}{n!}
\sum_{\substack{1\leq l_1,\ldots,l_n\leq\upsilon\\ l_1+\cdots +l_n=i+1}}
\int_0^s\|\gamma_{l_1}(t)\|\|\gamma_{l_2}(t)\|\cdots\|\gamma_{l_{n-1}}(t)\|
\|\dot\gamma_{l_n}(t)\|\,dt\\
&&\!\!\!\leq\sum_{n=2}^\upsilon\frac{\beta^n}{n!}
\sum_{\substack{1\leq l_1,\ldots,l_n\leq\upsilon\\ l_1+\cdots +l_n=i+1}}
\int_0^s\|\dot\gamma_{l_1}(t)\|\,dt\int_0^s\|\dot\gamma_{l_2}(t)\|\,dt\cdots
\int_0^s\|\dot\gamma_{l_{n-1}}(t)\|\,dt\int_0^s\|\dot\gamma_{l_n}(t)\|\,dt\\
&&\leq\sum_{n=2}^\upsilon\frac{\beta^n}{n!}
\sum_{\substack{1\leq l_1,\ldots,l_n\leq\upsilon\\ l_1+\cdots +l_n=i+1}}
\kappa_{l_1}\cdots\kappa_{l_n}\,F_{l_1}(L)\cdots F_{l_n}(L)\
\cA_0^s(\dot\gamma_1-X)^n\,s^{i+1}\\
&&\leq\sum_{n=2}^\upsilon\frac{\beta^n}{n!}
\sum_{\substack{1\leq l_1,\ldots,l_n\leq\upsilon\\ l_1+\cdots +l_n=i+1}}
\kappa_{l_1}\cdots\kappa_{l_n}\,F_{l_1}(L)\cdots F_{l_n}(L)\
(2L)^{n-1}\,\cA_0^s(\dot\gamma_1-X)\,s^{i+1}\,.
\end{eqnarray*}
This proves estimates \eqref{indsj} for every $j=2,\ldots,\upsilon$.
To complete the proof, it suffices to apply these estimates
to the Lipschitz curve $\tilde\Gamma(t)=\Gamma(-t)$
and replace $X$ by $-X$. $\Box$
\begin{Cor}\label{corlt}
There exists a strictly increasing function
$\Upsilon:[0,+\infty[\lra[0,+\infty[$, infinitesimal at zero
and only depending on $\cM$ that has the following property.
Let $\Gamma:[a,b]\lra\M$ be a horizontal curve
and assume that $t\in[a,b]$ is an approximate continuity
points $t$ of $\dot\gamma_1$. Then
\begin{eqnarray}\label{estimdif}
\;\;\left\|\pi_i\left(\exp^{-1}\Big(
\delta_{1/h}\big(\Gamma(t)^{-1}\Gamma(t+h)\big)\Big)\right)\right\|
\leq\Upsilon(L_t)\;\cA^h_t\big(\dot\gamma_1-\dot\gamma_1(t)\big)
\end{eqnarray}
for every $i=2,\ldots,\upsilon$, where
$L_t=\max\left\{\cA^{b-t}_t(\dot\gamma_1),\cA^{a-t}_t(\dot\gamma_1)\right\}$.
\end{Cor}
{\sc Proof.}
Clearly, $\gamma_1$ is absolutely continuous.
Let $t$ be an approximate continuity point
of $\dot\gamma_1$ and then also a differentiability point.
It suffices to apply Lemma~\ref{theoestim} to the curve
$h\lra\Gamma(t)^{-1}\Gamma(t+h)$ with $X=\dot\gamma_1(t)$.
In fact, we observe that
$$
\pi_1\left(\exp^{-1}\big(\Gamma(t)^{-1}\Gamma(t+h)\big)\right)
=\gamma_1(t+h)-\gamma_1(t)
$$
and that 
$$
\cA^\lambda_0\big(\dot\gamma_1(t+\cdot)-\dot\gamma_1(t)\big)
=\cA^\lambda_t\big(\dot\gamma_1(\cdot)-\dot\gamma_1(t)\big)\,.
$$
This implies that
\begin{eqnarray}
\;\;\left\|\pi_i\Big(\exp^{-1}\big(\Gamma(t)^{-1}\Gamma(t+h)\big)\Big)\right\|\leq\Upsilon_i(L_t)\,|h|^i\,
\cA^h_t\big(\dot\gamma_1-\dot\gamma_1(t)\big)
\end{eqnarray}
where $\Upsilon_i$ are given in Lemma~\ref{theoestim}
and 
$$
L_t=\max\left\{\cA^{b-t}_t(\dot\gamma_1),\cA^{a-t}_t(\dot\gamma_1)\right\}\,.
$$
Thus, setting $\Upsilon=\max_{i=2,\ldots,\upsilon}\Upsilon_i$,
and using the definition of dilations, our claim follows. $\Box$
\begin{The}\label{pansu}
There exists a nondecreasing function
$\Upsilon:[0,+\infty[\lra[0,+\infty[$ only depending on $\cM$
with the following property.
If $\Gamma:[a,b]\lra\M$ is a horizontal curve,
then for all approximate continuity points $t$ of $\dot\gamma_1$
the following estimate
\begin{eqnarray}\label{pansuestim}
\;\qquad\left\|\delta_{1/h}
\Big(-h\,\dot\gamma_1(t)\gope\big(-\gamma(t)\big)\gope\gamma(t+h)\Big)
\right\|\leq\Upsilon(L_t)\;\cA^h_t\big(\dot\gamma_1-\dot\gamma_1(t)\big)
\end{eqnarray}
holds, where $L_t=\max\left\{\cA^{b-t}_t(\dot\gamma_1),\cA^{a-t}_t(\dot\gamma_1)\right\}$.
In particular, $\Gamma$ is a.e.
P-differentiable.
\end{The}
{\sc Proof.}
Let $t$ be an approximate continuity point of $\dot\gamma_1$ and define
$$
h\lra\Gamma(t)^{-1}\Gamma(t+h)
=\exp\theta(h)
=\exp\big(\theta_1(h)+\cdots+\theta_\upsilon(h)\big)\,,
$$
with $\theta_i(h)\in V_i$. Notice that in particular
$\theta_1(h)=\gamma_1(t+h)-\gamma_1(t)$. In view of \eqref{estimdif},
there exists a constant $K>0$ depending on $\cM$ and
a strictly increasing function $\Upsilon:[0,+\infty[\lra[0,+\infty[$
infinitesimal at zero such that
\begin{eqnarray}\label{thetahi}
\left\|\frac{\theta_i(h)}{h^i}\right\|\leq \Upsilon(L_t)\,
\cA^h_t\big(\dot\gamma_1-\dot\gamma_1(t)\big).
\end{eqnarray}
for every $i=2,\ldots,\upsilon$, where we have set
$$
L_t=\max\left\{\cA^{b-t}_t(\dot\gamma_1),
\cA^{a-t}_t(\dot\gamma_1)\right\}.
$$
Thus, taking into account that $t$ is a differentiability point of
$\gamma_1$ and that $0$ is an approximate continuity point of
$\dot\theta_1$, we obtain that $\Gamma$ is P-differentiable at $t$ and
that
$$
D\Gamma(t)(\lambda)=\exp\big(\lambda\,\dot\gamma_1(t)\big)\,.
$$
Notice also that
$$
\left\|\frac{\theta_1(h)}{h}-\dot\gamma_1(t)\right\|
\leq \cA^h_t\big(\dot\gamma_1-\dot\gamma_1(t)\big),
$$
then \eqref{thetahi} yields
\begin{eqnarray}\label{thetamin}
\left\|\delta_{1/h}\theta(h)-\dot\gamma_1(t)\right\|\leq
\Big[1+(\upsilon-1)\Upsilon(L_t)\Big]
\;\cA^h_t\big(\dot\gamma_1-\dot\gamma_1(t)\big)\,.
\end{eqnarray}
Applying \eqref{cnrest} and \eqref{thetamin} to
$$
-\dot\gamma_1(t)\gope\delta_{1/h}\theta(h)
=\delta_{1/h}\theta(h)-\dot\gamma_1(t)+
\sum_{n=2}^\upsilon c_n\big(-\dot\gamma_1(t),\delta_{1/h}\theta(h)\big)\,,
$$
we get the estimate
\begin{eqnarray}\label{longestdot}
&&\left\|-\dot\gamma_1(t)\gope\delta_{1/h}\theta(h)\right\|\leq
\Big[1+(\upsilon-1)\Upsilon(L_t)\Big]
\;\cA^h_t\big(\dot\gamma_1-\dot\gamma_1(t)\big)\\
&&+\sum_{n=2}^\upsilon\frac{\beta^{n-1}}{n!}
\left\|\frac{\delta_{1/h}\theta(h)+\dot\gamma_1(t)}{2}\right\|^{n-1}
\Big[1+(\upsilon-1)\Upsilon(L_t)\Big]
\cA^h_t\big(\dot\gamma_1-\dot\gamma_1(t)\big)\nonumber\\
&&+\sum_{n=2}^\upsilon
\left\|R_n\left(-\dot\gamma_1(t),\delta_{1/h}\theta(h)\right)\right\|\,.\nonumber
\end{eqnarray}
Observing that $\|\dot\gamma_1(t)\|\leq L_t$
and $\cA^h_t\big(\dot\gamma_1-\dot\gamma_1(t)\big)\leq2L_t$,
inequality \eqref{thetamin} implies that
\begin{eqnarray*}
&&\left\|\delta_{1/h}\theta(h)\right\|\leq
3L_t\Big[1+(\upsilon-1)\Upsilon(L_t)\Big]=\Upsilon_1(L_t),
\end{eqnarray*}
then we apply \eqref{estimsmall} to the third line of
\eqref{longestdot}, obtaining
\begin{eqnarray}
&&\left\|-\dot\gamma_1(t)\gope\delta_{1/h}\theta(h)\right\|\leq
\bigg(1+\sum_{n=2}^\upsilon\frac{\beta^{n-1}}{n!}
\left(\frac{\Upsilon_1(L_t)+L_t}{2}\right)^{n-1}\bigg)\\
&&\,\Big[1+(\upsilon-1)\Upsilon(L_t)\Big]
\cA^h_t\big(\dot\gamma_1-\dot\gamma_1(t)\big)\nonumber \\
&& +\Big(\sum_{n=2}^\upsilon C\lLs n,\Upsilon_1(L_t)\rLs\Big)\,
\Big[1+(\upsilon-1)\Upsilon(L_t)\Big]^3\,4\,L_t^2\;
\cA^h_t\big(\dot\gamma_1-\dot\gamma_1(t)\big)\nonumber\,.
\end{eqnarray}
This concludes the proof. $\Box$
\begin{Cor}[Pansu]\label{Pansugrad}
Let $A$ be a measurable set of $\G$ and let 
$f:A\lra\M$ be a Lipschitz mapping. Then $f$ is a.e.
P-differentiable.
\end{Cor}
{\sc Proof.}
By Theorem~\ref{pansu}, rectifiable curves in graded groups
are in particular a.e. P-differentiable. This allows us to extend Proposition~3.6
of \cite{Mag} to the case of mappings with values in a graded group $\M$.
Thus, repeating the same arguments of the proof of Theorem~3.9
of \cite{Mag} the claim is achieved. $\Box$
%
%
%
%
%
%
%
%
%
%
%
%
%
%
%
%
%
%
%
%
%
%
%
%
%
%
%
\subsection{Characterization of P-differentiable mappings}\label{contactmap}
Here we denote by $\G$ and $\M$ two arbitrary graded groups,
where $\G$ is stratified and $\Omega\subset\G$ is an open set.
To obtain estimates on the P-difference quotient of
P-differentiable mappings, we will use the following family
of piecewise horizontal lines.
\begin{Def}\label{PN}{\rm
Let $N$ be a positive integer, $X_1,\ldots,X_m$ be a basis of the first
layer $V_1$, $i_1,\ldots,i_N\in\{1,\ldots,m\}$ be fixed integers and
for every $a=(a_1,\ldots,a_N)\in\R^N$ define
\begin{eqnarray}
P^s(a)=\left\{\begin{array}{ll}
e &  \mbox{if $s=0$} \\
\delta_{a_1}h_{i_1}\delta_{a_2}h_{i_2}\cdots\delta_{a_s}h_{i_s} &
\mbox{if $s=1,\ldots,N$}
\end{array}\right.\,,
\end{eqnarray}
where $h_1=\exp X_1,\ldots,h_m=\exp X_m$ and we have assumed that $d(h_i)=1$.
}\end{Def}
Lemma 1.40 of \cite{FS} yields the following
\begin{Lem}\label{140FS}
For every stratified group $\G$ there exists an integer $N$
and a family of integers $\{i_1,\ldots,i_N\}\subset\{1,\ldots,m\}$ 
depending on $\G$, such that the mapping $P^N$ of Definition~\ref{PN}
sends a neighbourhood of the origin in $\R^N$ onto a neighbourhood of the identity $e\in\G$.
\end{Lem}
\begin{Def}\label{defcaDef}{\rm
Under conditions of Lemma~\ref{140FS}, we define the number
\begin{eqnarray}\label{defca}
c(\G,d)=\max_{\substack{s=1,\ldots,N\\ a=(a_s)\in (P^N)^{-1}(D_1)}}|a_s|,
\end{eqnarray}
that only depends on the algebraic structure of $\G$ and on the
homogeneous distance $d$ used to define both $P^N$ and the closed unit ball
$D_1$.
}\end{Def}
\begin{Rem}\label{rmkPN}{\rm
Being the mapping $P^N$ homogeneous with respect to dilations,
namely, $P^N(ra)=\delta_r\big(P^N(a)\big)$, it follows that
$P^N$ is surjective onto $\G$. Homogeneity yields
\begin{eqnarray}\label{defcar}
\max_{\substack{s=1,\ldots,N\\ a=(a_s)\in (P^N)^{-1}(D_r)}}|a_s|
=r\;c(\G,d)
\end{eqnarray}
for every $r>0$ and implies that the preimage of compact
sets of $\G$ is compact in $\R^N$.
This shows that definition \eqref{defca} is well posed.
}\end{Rem}
\begin{Def}{\rm
Let $f:K\lra Y$ be a vector valued continuous function
on a compact metric space $(K,\rho)$.
Then we define the {\em modulus of continuity} of $f$ on $K$ as 
\begin{eqnarray}\label{modcont}
\omega_{K,f}(t)=\max_{\substack{x,y\in K\\ \rho(x,y)\leq t}}
\|f(x)-f(y)\|\,.
\end{eqnarray}
}\end{Def}
\begin{The}\label{PdifPdifj}
Let $f:\Omega\lra\M$ be P-differentiable at $x\in\Omega$,
where $f=\exp\circ F$ and $F=F_1+\cdots+F_\upsilon$
with $F_j:\Omega\lra W_j$.
Then every $F_j$ is P-differentiable at $x$ and we have the formulae
\begin{eqnarray}
&&\pi_1\circ df(x)=dF_1(x) \label{eqF_1}\\
&&dF_i(x)(h)=\sum_{n=2}^\upsilon\frac{(-1)^n}{n!}\;\pi_i
\left([F(x),dF(x)(h)]_{n-1}\right)\,. \label{eqDF_j}
\end{eqnarray}
for every $i=2,\ldots,\upsilon$ and every $h\in\G$ 
\end{The}
{\sc Proof.}
We will use the notation $|h|=d(h)$.
By definition of P-differentiability we have that
$$
\delta_{1/|h|}\left(f(x)^{-1}f(xh)\right)
$$
converges to $Df(x)(h)$ as $h\ra e$, uniformly with respect to the parameter
$\delta_{1/|h|}h$ varying in a compact set.
Then the difference quotients 
$$
\pi_i\circ\exp^{-1}\Big[\delta_{1/|h|}\left(f(x)^{-1}f(xh)\right)\Big]\,,
$$
uniformly converge for every $i=1,\ldots,\upsilon$.
In other words, the difference quotient
$$
\frac{F_1(xh)-F_1(x)}{|h|}
$$ 
along with 
\begin{equation}\label{convi}
\frac{F_i(xh)-F_i(x)+\sum_{n=2}^\upsilon
\pi_i\left(c_n\lls -F(x),F(xh)\rls\right)}{|h|^i}
\end{equation}
converge as $h\ra0$ whenever $i=2,\ldots,\upsilon$.
In particular, $F_1$ is $P$-differentiable at $x$
and
\[
\pi_1\circ\exp^{-1}\big(Df(x)(h)\big)=DF_1(x)(h)
\]
that implies \eqref{eqF_1}.
Writing \eqref{convi} for $i=2$, we notice that this expression
can be written as 
\begin{eqnarray}\label{conv2}
\frac{F_2(xh)-F_2(x)-
\frac{1}{2}[F_1(x),F_1(xh)]}{|h|^2}\,.
\end{eqnarray}
The convergence of the previous quotient and P-differentiability of $F_1$
at $x$ imply that  $F_2$ is also differentiable at $x$ and there holds
$$
DF_2(x)(h)=\frac{1}{2}\,[F_1(x),DF_1(x)(h)].
$$
By induction, we assume that for $j\geq2$, the vector-valued mapping $F_i$
is differentiable at $x$ and that 
\begin{eqnarray*}
	DF_i(x)(h)=\sum_{n=2}^\upsilon\frac{(-1)^n}{n!}\;
	\pi_i\left([\overline F_{i-1}(x),
	D\overline F_{i-1}(x)(h)]_{n-1}\right)
\end{eqnarray*}
holds for every $i=1,\ldots,j$, where we have set
\[
\overline F_i=F_1+F_2+\cdots+F_i.
\]
By P-differentiability \eqref{convi}, we have
\begin{equation}\label{convip1}
\frac{F_{i+1}(xh)-F_{i+1}(x)+\sum_{n=2}^\upsilon
\pi_{i+1}\left(c_n\lls -F(x),F(xh)\rls\right)}{|h|}\lra0\quad
\mbox{as}\quad h\ra0\,.
\end{equation}
Now we notice that
\begin{equation}\label{eqci}
\pi_{i+1}\left(c_n\lls -F(t),F(t+h)\rls\right)
=\pi_{i+1}\left(c_n\lls -\overline{F}_i(t),
\overline{F}_i(t+h)\rls\right)\,.
\end{equation}
Furthermore, inductive hypothesis and Lemma~\ref{keylemma}
ensure the existence of 
\begin{equation}\label{indhypl}
\lim_{|h|\ra0}
\frac{c_n(-\overline F_i(x),\overline F_i(x+h))}{|h|}=
\frac{(-1)^{n-1}}{n!}
\left[\overline{F}_i(x),D\overline{F}_i(x)(h)\right]_{n-1}.
\end{equation}
Thus, joining \eqref{convip1}, \eqref{eqci} and \eqref{indhypl},
it follows that there exists $DF_{i+1}(x)$ and
$$
DF_{i+1}(x)(h)=\sum_{n=2}^{i+1}\frac{(-1)^n}{n!}\;\pi_{i+1}
\left([\overline{F}_i(x),D\overline F_i(x)]_{n-1}\right)\,.
$$
This shows that $F:\Omega\lra\cM$ is $P$-differentiable at $t$,
hence in the previous formula we can replace $\overline F_i$
by $F$, achieving \eqref{eqDF_j}.
This ends the proof. $\Box$
\begin{Cor}\label{Pansuimpl}
Let $\Gamma:[a,b]\lra\M$ be P-differentiable at $t\in[a,b]$.
Then $\gamma$ is differentiable at t and satisfies
equations \eqref{eqhdiffcurve}.
\end{Cor}
This corollary is an immediate consequence of Theorem~\ref{PdifPdifj}.
Now we can prove the first result stated in the introduction.
\vskip.2truecm
{\sc Proof of Theorem~\ref{PdifContact}.}
If $f:\Omega\lra\M$ is continuously P-differentiable, then 
Theorem~\ref{PdifPdifj} implies the P-differentiability of all $F_j$'s and the 
validity of both formulae \eqref{eqF_11} and \eqref{eqDF_j1}.
These equations in turn give the continuoity of the P-differential
of every $F_j$. This shows that the first condition implies the second one.

Now, we assume that the second condition holds.
As observed in Remark~\ref{Pdifhdif}, P-differentiability of
vector-valued mappings implies h-differentiability, hence
all $F_j$'s are in particular h-differentiable.
In addition, formula \eqref{dhfP} yields the continuity of the horizontal differential. By our assumption \eqref{eqF_11} and
\eqref{eqDF_j1} hold, therefore combining these equations with
Proposition~\ref{horprop}, it follows that $f$ also has the
contact property. We have shown that $f$ is continuously h-differentiable contact, that is the third condition.

Finally, we assume that $f$ is a continuously h-differentiable
contact mapping and consider
\[
\Gamma_X(t)=\exp\big(\sum_{j=1}^\iota\gamma_{X,j}(t)\Big)
=f\big(c_X(t)\big)
=\exp\Big(\sum_{j=1}^\iota F_j\big(c_X(t)\big)\Big)
\]
where $c_X(t)=x\exp(tX)$ and
$\exp^{-1}\circ\Gamma_X=\gamma_X=\sum_{j=1}^\iota\gamma_{X,j}$ with
$t\lra\gamma_{X,j}(t)\in V_j$.
By definition of h-differentiability we immediately observe that
also $F_j$ are continuously h-differentiable
and we have 
\[
\dot\gamma_{X,j}(t)=XF_j(c_X(t))
\]
By our assumptions, we know that $\Gamma_X$ is a horizontal curve,
then we apply the key estimate \eqref{pansuestim} to $\Gamma_X$
at $t=0$, achieving
\[
\left\|\delta_{1/h}
\Big(-h\,XF_1(x)\gope\big(-\gamma_X(0)\big)\gope\gamma_X(h)\Big)
\right\|\leq\Upsilon(L)\;\cA^h_0\big((XF_1)\circ c_X-XF_1(x)\big)
\]
where $\Upsilon:[0,+\infty[\lra[0,+\infty[$ is a nondecreasing function
depending on $\M$, $\cA_0^h$ is defined in \eqref{cMt} and 
\begin{eqnarray*}
L=\max\left\{\cA^\ep_0(\dot\gamma_{X,1}),\cA^{-\ep}_0(\dot\gamma_{X,1})\right\}.
\end{eqnarray*}
Here $\ep>0$ is fixed such that
$c_X([-\ep,\ep])\subset U$ for every $X\in V_1$ with $\|X\|=1$,
where $U\subset\Omega$ is a fixed compact neighbourhood of $x$.
For every $|h|\leq\ep$, we have
\begin{eqnarray*}
	\cA^h_0\big((XF_1)\circ c_X-XF_1(x)\big)
	\leq \max_{\substack{y,z\in U\\	\|\exp^{-1}(y)-\exp^{-1}(z)\|\leq|h|}}
	\|XF_1(z)-XF_1(y)\|\,.
\end{eqnarray*}
We notice that $\dot\Gamma_X(0)=XF_1(x)=d_Hf(x)(X)$, then setting
\begin{eqnarray*}
\omega_{d_Hf,U}(|h|)
=\max_{\substack{y,z\in U\\	\|\exp^{-1}(y)-\exp^{-1}(z)\|\leq\ep}}
\|d_Hf(y)-d_Hf(z)\|
\end{eqnarray*}
and observing that
\[
L\leq\max_{y\in U}\|d_Hf(y)\|
\]
\begin{equation}\label{uniflipest}
\left\|\delta_{1/h}
\Big(\!-h\,XF_1(x)\gope\big(\!-\gamma_X(0)\big)\gope\gamma_X(h)\Big)
\right\|\leq\Upsilon\big(\max_{y\in U}\|d_Hf(y)\|\big)\;
\omega_{d_Hf,U}(|h|)\,.
\end{equation}
The right hand side is independent of $x$ and tends to zero
as $h\ra0$, then taking the exponential of elements in the
left hand side, we get the uniform convergence of
\[
\delta_{1/h}\bigg(\exp\big(-h XF_1(y)\big)
\Big(f(y)^{-1}f(y\exp(hX)\Big)\bigg)
\]
to the unit element $e$ as $h\ra0$, where
$y$ varies in the smaller neighbourhood 
$$
U_{-\ep}=\{y\in U\mid \mbox{dist}_{\|\cdot\|}(y,U^c)\geq\ep\}\,.
$$
of $x$. Thus, arguing as in Corollaire~3.3 of \cite{Pan2}, we achieve 
the P-differentiability of $f$ in a possibly smaller neighbourhood
of $x$. This shows the everywhere P-differentiability of $f$.
Observing that
\begin{eqnarray}\label{hpdif}
d_Hf(x)(X)=df(x)(X)
\end{eqnarray}
for every $x\in\Omega$ and $X\in V_1$, we get that 
$x\lra Df(x)_{|\exp V_1}$ is continuous. 
The fact that $df(x)_{|V_j}$ polynomially depends on 
$df(x)_{|V_1}$ concludes the proof. $\Box$
%
%
%
%
%
%
%
%
%
%
%
%
%
%
%
%
\section{Absolutely continuous curves in graded groups}\label{ABSC}
This section is devoted to the characterization of 
absolutely continuous curves in graded groups. $\M$ and $\cM$
denote a graded group along with its Lie algebra,
equipped with a homogeneous distance $\rho$ and 
a norm $\|\cdot\|$, respectively.
\begin{Def}{\rm
Let $\gamma:[0,s]\lra\cM$ be an absolutely continuonus curve in $\cM$
and let $\Sigma=(t_0,t_1,\ldots,t_N)$ be a partition
of $[0,s]$, where $t_0=0$, $t_i<t_{i+1}$ and $t_N=s$ for
some $N\in\mathbb N$. The associated ``sum'' with respect to
the group operation in $\cM$ is defined by
\begin{eqnarray}\label{ARiemsum}
\sigma_\Sigma(\gamma)(s)=\sum_{k=0}^{N-1}\;
\gamma(t_k)^{-1}\gope\gamma(t_{k+1})\,.
\end{eqnarray}
We also set $\|\Sigma\|=\max_{1\leq i\leq N}(t_i-t_{i-1})$.
}\end{Def}
\begin{Lem}\label{laireb}
Let $\gamma:[0,s]\lra\cM$ be an absolutely continuous curve. Then there exists 
\begin{eqnarray}\label{aireb}
\lim_{\|\Sigma\|\ra0^+}
\sigma_\Sigma(s)
=\gamma(s)-\gamma(0)+\sum_{n=2}^\iota\frac{(-1)^{n-1}}{n!}
\int_0^s[\gamma(l),\dot\gamma(l)]_{n-1}\,dl\,.
\end{eqnarray}
\end{Lem}
{\sc Proof.}
Applying \eqref{absBCH} and \eqref{ARiemsum}, we get
$$
\sigma_\Sigma(s)=\sum_{k=0}^{N-1} \sum_{n=1}^\iota c_n\big(-\gamma(t_k),\gamma(t_{k+1})\big)\,,
$$
then formula \eqref{cnrest} yields
\begin{eqnarray}\label{developc}
&&\sigma_\Sigma(s)=\gamma(s)-\gamma(0)
+\sum_{k=0}^{N-1}\sum_{n=2}^\iota\frac{(-1)^{n-1}}{n!}
\left[\frac{\gamma(t_k)+\gamma(t_{k+1})}{2},
\gamma(t_{k+1})-\gamma(t_k)\right]_{n-1}\\
&&+\sum_{k=0}^{N-1}\sum_{n=2}^\iota R_n\big(-\gamma(t_k),\gamma(t_{k+1})\big)
\,.\nonumber
\end{eqnarray}
We define $\mu_k=\big(\gamma(t_k)+\gamma(t_{k+1})\big)/2$,
then the absolute continuity of $\gamma$ yields
\[
\sum_{k=0}^{N-1}\left[\mu_k,\gamma(t_{k+1})-\gamma(t_k)\right]_{n-1}=
\sum_{k=0}^{N-1}\left[\mu_k,\int_{t_k}^{t_{k+1}}\dot\gamma(l)\,dl\right]_{n-1}=\sum_{k=0}^{N-1}
\int_{t_k}^{t_{k+1}}\left[\mu_k,\dot\gamma(l)\right]_{n-1}\,dl\,.
\]
Now we notice that
\begin{eqnarray*}
&&\sum_{k=0}^{N-1}
\int_{t_k}^{t_{k+1}}\left[\mu_k,\dot\gamma(l)\right]_{n-1}\,dl=
\int_0^s[\gamma(l),\dot\gamma(l)]_{n-1}\,dl\\
&&+\sum_{p=0}^{n-2}\sum_{k=0}^{N-1}\int_{t_k}^{t_{k+1}}
[\mu_k,[\mu_k-\gamma(l),[\gamma(l),\dot\gamma(l)]_p]]_{n-p-2}\,dl,
\end{eqnarray*}
hence the estimate
\begin{eqnarray*}
&&\left\|\;\sum_{p=0}^{n-2}\sum_{k=0}^{N-1}\int_{t_k}^{t_{k+1}}
[\mu_k,[\mu_k-\gamma(l),[\gamma(l),\dot\gamma(l)]_p]]_{n-p-2}\,dl\;\right\|\\
&&\leq\,(p-1)\,\beta^{n-2}\,\eta(\|\Sigma\|)\,
\max_{t\in[0,s]}\|\gamma(t)\|^{n-2}\,
\int_{0}^{s}\|\dot\gamma(l)\|\,dl
\end{eqnarray*}
immediately proves the validity of the following
\begin{eqnarray}\label{airebalaye}
\lim_{\|\Sigma\|\ra0^+}
\sum_{k=0}^{N-1}\left[\frac{\gamma(t_k)+\gamma(t_{k+1})}{2},
\gamma(t_{k+1})-\gamma(t_k)\right]_{n-1}=
\int_0^s[\gamma(l),\dot\gamma(l)]_{n-1}\,dl\,,
\end{eqnarray}
where we have set
$$
\eta(\tau)=\sup_{|\beta-\alpha|\leq\tau}
\left|\,\int_\alpha^\beta\|\dot\gamma(l)\|\,dl\right|
$$
which goes to zero as $\tau\ra0^+$, due to the absolute continuity
of $\gamma$.
Taking into account \eqref{estimsmall}, \eqref{developc} and \eqref{airebalaye},
then limit \eqref{aireb} follows. $\Box$
\begin{Rem}{\rm
Incidentally, the limit \eqref{airebalaye} exactly coincides with the 
``$(n$$-$$1)$-i\`eme aire balay\'ee'' of Definition~4.6 in \cite{Pan2}.}
\end{Rem}
\begin{Def}{\rm
Let $\Gamma:[a,b]\lra\M$ be a continuous curve. Then we define
\[
\mbox{\rm Var}_a^t\Gamma=\sup_{\substack{t_0=a<t_1<\cdots<t_N=t\\ N\in\N}}\sum_{k=1}^N
\rho\left(\Gamma(t_k),\Gamma(t_{k-1})\right)\,.
\]
}\end{Def}
\begin{The}\label{hdiffcurve}
Let $\Gamma:[a,b]\lra\M$ be a curve and define $\gamma=\exp^{-1}\circ\Gamma=\sum_{i=1}^\upsilon\gamma_i$,
where $\gamma_i$ takes values in $V_i$.
Then the following statements are equivalent:
\begin{enumerate}
	\item $\Gamma$ is absolutely continuous,
	\item $\gamma$ is absolutely continuous and
	the differential equation
	\begin{eqnarray}\label{eqhdiffcurve1}
	\dot\gamma_i(t)=\sum_{n=2}^\upsilon\frac{(-1)^n}{n!}\;
	\pi_i\left([\gamma(t),\dot\gamma(t)]_{n-1}\right)
	\end{eqnarray}
	is satisfied a.e. for every $i\geq2$.
\end{enumerate}
If one of the previous conditions holds, then
there exists a constant $C>0$ only depending
on $\rho$ and $\|\cdot\|$,
such that for any $\tau_1<\tau_2$, we have
\begin{eqnarray}\label{absestim}
\rho\big(\Gamma(\tau_1),\Gamma(\tau_2)\big)\leq C\,\int_{\tau_1}^{\tau_2}\|\dot\gamma_1(t)\|\,dt\,.
\end{eqnarray}
\end{The}
{\sc Proof.}
We first assume that $\Gamma$ is absolutely continuous
with respect to a homogeneous distance $\rho$ fixed in $\M$.
Observing that the image of $\Gamma$ is bounded and applying Lemma~\ref{estimi} with $i=1$, we immediately conclude that
$\gamma$ is also absolutely continuous with respect
to the norm $\|\cdot\|$ fixed in $\cM$; then $\gamma$ is a.e. differentiable. 
Let $a<\tau<b$ be both an approximate continuity point
of $\dot\gamma$ and a differentiability point of the total variation
\begin{eqnarray}\label{veegamma}
[a,b]\ni l\lra\mbox{\rm Var}_a^l\Gamma
\end{eqnarray}
We fix $\ep>0$ suitable small.
We fix a bounded open neighbourhood $U$ of $e$ and choose
$\delta>0$ such that $\exp\big(\gamma(\tau+t_{i-1})^{-1}\gope\gamma(\tau+t_i)\big)\in U$
whenever $\Sigma=(t_0,t_1,\ldots,t_N)$ is an arbitrary
partition of $[0,s]$ with $\|\Sigma\|\leq\delta$ and $0<s\leq\ep$.
For every $i=1,\ldots,\upsilon$, there exists a constant $k_i>0$ such that
\begin{eqnarray}
&&\|\pi_i\big(\sigma_\Sigma(\gamma(\tau+\cdot))(s)\big)\|\leq k_i\;
\|\pi^i\big(\sigma_\Sigma(\gamma(\tau+\cdot))(s)\big)\|\\
&&\leq k_i\, K_U\,\sum_{j=1}^N \rho\big(\Gamma(\tau+t_j),\Gamma(\tau+t_{j-1})\big)^i\leq
k_i\, K_U\,
\bigg(\mbox{\rm Var}_{\tau}^{\tau+s}\Gamma\bigg)^i\,.
\end{eqnarray}
Taking the limit in the previous inequality as $\|\Sigma\|\lra0^+$
and applying Lemma~\ref{laireb}, we obtain
\begin{eqnarray*}
&&\Big\|\gamma_i(\tau+s)-\gamma_i(\tau)
+\sum_{n=2}^\upsilon\frac{(-1)^{n-1}}{n!}
\int_0^s\,\pi_i\left([\gamma(\tau+l),\dot \gamma(\tau+l)]_{n-1}\right)\,dl\Big\| \\
&&\leq\,k_i\, K_U\,\bigg(\mbox{\rm Var}_{\tau}^{\tau+s}\Gamma\bigg)^{i-1}
\bigg(\mbox{\rm Var}_a^{\tau+s}\Gamma-\mbox{\rm Var}_a^\tau\Gamma\bigg)
\end{eqnarray*}
By our assumptions, $\tau$ is a differentiability point of both $\gamma$
and of the total variation \eqref{veegamma} and it is an approximate
continuity point of $\dot\gamma$. Thus, dividing by $s$ the last
inequality and taking the limit as $s\ra0^+$ for every $i\geq2$,
equations \eqref{eqhdiffcurve1} follow.

Now we assume the $\gamma$ is absolutely continuous and that
\eqref{eqhdiffcurve1} a.e. holds for every $i\geq2$. 
Let us fix $a<\tau_1<\tau_2<b$ and consider
$$
h\lra\Theta(h)=\Gamma(\tau_1)^{-1}\Gamma(\tau_1+h)=\exp\theta(h)=
\exp\big(\theta_1(h)+\theta_2(h)+\cdots+\theta_\upsilon(h)\big)\,,
$$
where $\theta_i(h)\in V_i$.
In view of Proposition~\ref{horprop} and observing that $\Theta$,
as left translated of $\Gamma$, is a.e. horizontal, we have that
\begin{eqnarray}\label{eqhdiffcurvetheta}
\dot\theta_i(h)=\sum_{n=2}^\upsilon\frac{(-1)^n}{n!}\;
\pi_i\left([\theta(h),\dot\theta(h)]_{n-1}\right)
\end{eqnarray}
a.e. holds and $\theta$ is absolutely continuous.
Now, we wish to prove by induction that there exists a constant
$C_i>0$, only depending on the norm constant $\beta$, such that
\begin{eqnarray}\label{piabs}
\int_0^{h_0}\|\dot\theta_i(t)\|\,dt\leq C_i\left(\int_0^{h_0}\|\dot\theta_1\|\right)^i
\end{eqnarray}
for every $i\geq1$, with $h_0=\tau_2-\tau_1$.
To prove \eqref{piabs},
we first write \eqref{eqhdiffcurvetheta} as follows
\begin{equation}\label{confshapetheta}
\dot\theta_i(h)=\sum_{n=2}^\upsilon\frac{(-1)^n}{n!}
\sum_{\substack{1\leq l_1,\ldots,l_n\leq\upsilon\\l_1+\cdots+l_n=i}}
[\theta_{l_1}(h),[\theta_{l_2}(h),[\cdots,[\theta_{l_{n-1}}(h),\dot\theta_{l_n}(h)],],\ldots,]\,,
\end{equation}
since $\M$ is a graded group.
For $i=1$, inequality \eqref{piabs} is trivial, choosing $C_1\geq1$.
Let us assume that \eqref{piabs} holds for every $i\leq j$.
Due to \eqref{confshapetheta} for $i=j+1$ and our
induction hypothesis, we have
\begin{eqnarray*}
&&\int_0^{h_0}\|\dot\theta_{j+1}(t)\|\,dt\leq	
\sum_{n=2}^\upsilon\frac{\beta^{n-1}}{n!}\!\!\!\!\!
\sum_{\substack{1\leq l_1,\ldots,l_n\leq\upsilon\\l_1+\cdots+l_n=j+1}}
\!\!\!\int_0^{h_0}\|\theta_{l_1}(t)\|\,\|\theta_{l_2}(t)\|\cdots
\|\theta_{l_{n-1}}(t)\|\,\|\dot\theta_{l_n}(t)\|\,dt\\
&&\leq\sum_{n=2}^\upsilon\frac{\beta^{n-1}}{n!}
\sum_{\substack{1\leq l_1,\ldots,l_n\leq\upsilon\\l_1+\cdots+l_n=j+1}}
C_{l_1}\cdots C_{l_n}
\Big(\int_0^{h_0}\|\dot\theta_1(t)\|\,dt\Big)^{j+1}\\
&&=\left(\int_0^{h_0}\|\dot\theta_1(t)\|\,dt\right)^{j+1}\,
\Big(\sum_{n=2}^\upsilon\frac{\beta^{n-1}}{n!}
\sum_{\substack{1\leq l_1,\ldots,l_n\leq\upsilon\\l_1+\cdots+l_n=j+1}}
C_{l_1}\cdots C_{l_n}\Big)\,.
\end{eqnarray*}
Thus, we have proved that 
\begin{eqnarray}\label{piabsfinal}
\|\theta_i(h_0)\|
=\left\|\pi_i\big(-\gamma(\tau_1)\gope\gamma(\tau_2)\big)\right\|\leq C_i\left(\int_{\tau_1}^{\tau_2}\|\dot\gamma_1\|\right)^i\,,
\end{eqnarray}
for every $i=1,\ldots,\upsilon$. These estimates show the validity
of \eqref{absestim}, hence $\Gamma$ is
absolutely continuous. $\Box$
\begin{Cor}\label{hdifflip}
Let $\Gamma:[a,b]\lra\M$ be a curve and define $\gamma=\exp^{-1}\circ\Gamma=\sum_{i=1}^\upsilon\gamma_i$,
where $\gamma_i$ takes values in $V_i$.
Then the following statements are equivalent:
\begin{enumerate}
	\item $\Gamma$ is Lipschitz,
	\item $\gamma$ is Lipschitz and the differential equation
	\begin{eqnarray}\label{eqhdiffcurve2}
	\dot\gamma_i(t)=\sum_{n=2}^\upsilon\frac{(-1)^n}{n!}\;
	\pi_i\left([\gamma(t),\dot\gamma(t)]_{n-1}\right)
	\end{eqnarray}
	is satisfied a.e. for every $i\geq2$.
\end{enumerate}
If one of the previous conditions holds, then
there exists a constant $C>0$ only depending on $\rho$ 
and $\|\cdot\|$, such that
\begin{eqnarray}\label{lipestim}
C^{-1}\;\mbox{\rm Lip}(\gamma_1)
\leq \mbox{\rm Lip}(\Gamma)\leq C\;\mbox{\rm Lip}(\gamma_1)
\end{eqnarray}
\end{Cor}
{\sc Proof.}
In view of Theorem~\ref{hdiffcurve}, if $\Gamma$ is Lipschitz,
then $\gamma$ is absolutely continuous and \eqref{eqhdiffcurve2} a.e. hold.
By the general estimate \eqref{lip1}, $\gamma_1$ is Lipschitz and we have a constant $C_1>0$, only depending on $\rho$ and $\|\cdot\|$, such that
\[
C_1^{-1}\;\Lip(\gamma_1)\leq \Lip(\Gamma)\,.
\]
Thus, from \eqref{absestim} the estimate \eqref{lipestim} follows.
Applying recursively the equations \eqref{eqhdiffcurve2}, the Lipschitz
property of $\gamma_1$ implies the Lipschitz property of
all $\gamma_j$'s. Conversely, if $\gamma$ is Lipschitz and
\eqref{eqhdiffcurve2} a.e. hold, then Theorem~\ref{hdiffcurve}
implies that $\Gamma$ is absolutely continuous
and the estimate \eqref{absestim} yields the Lipschitz property of $\Gamma$.
$\Box$
\begin{Cor}
Let $\Gamma:[a,b]\lra\M$ be an absolutely continuous curve.
Then the following formula holds
\[
\mbox{\rm Var}_a^b\Gamma=\int_a^b\rho\Big(\exp\big(\dot\gamma_1(t)\big)\Big)\,dt\,.
\]
where $\mbox{\rm Var}_a^b\Gamma$ is computed with respect to the
homogeneous distance $\rho$.
\end{Cor}
{\sc Proof.}
Our claim is a consequence of a.e. P-differentiability of $\Gamma$
proved in Theorem~\ref{pansu} and the general formula
\begin{eqnarray}\label{VarGamma}
\mbox{\rm Var}_a^b\Gamma=\int_a^b |\dot\Gamma(t)|
\end{eqnarray}
that holds in metric spaces, where $|\dot\Gamma(t)|$ is the metric 
derivative, that in this case coincides with $\rho\big(\exp\big(\dot\gamma_1(t)\big)\big)$,
by definition of P-differentiability. 
Formula \eqref{VarGamma} can be found in \cite{AmbTil} for Lipschitz curves
in metric spaces. Its extension to absolutely continuous curves
is straightforard and can be proved in a similar manner,
taking into account formula \eqref{absestim}. $\Box$
%
%
%
%
%
%
%
%
%
%
%
%
%
%
%
%
%
\subsection{Lipschitz property of P-differentiable mappings}\label{LpPD}
In the sequel, $\G$ will denote a stratified group with Lie algebra
$\cG$ and $\Omega\subset\G$ will be an open set.
The mapping $f:\Omega\lra\M$ will be represented in the algebra
by $F=\exp^{-1}\circ f:\Omega\lra\cM$ and $F_j=\pi_j\circ F$,
where $\pi_j:\cM\lra W_j$ are the canonical projections onto
the layers.

\begin{Lem}\label{PlipV}
Let $V$ be a finite dimensional normed vector space and
let $F:\Omega\lra V$ be continuously P-differentiable.
There exist constants $\mu_0>1$ and $C_0>0$, only depending
on $(\G,d)$, such that 
\[
\|F(x)-F(y)\|\leq C_0\max_{z\in D_{\xi,\mu_0r}}\|dF(z)\|\;\,d(x,y)
\]
for every $x,y\in D_{\xi,r}$, where $D_{\xi,\mu_0r}\subset\Omega$.
\end{Lem}
{\sc Proof.}
We consider $c_0=c(\G,d)$ as in \eqref{defca} and
fix $\mu_0=2(1+c_0N)$. If $D_{\xi,\mu_0r}\subset\Omega$ for a suitable small $r>0$, then one can check that
\[
A_{\xi,r}=\{y\in\G\mid \dist(y,B_{\xi,r})\leq 2\,c_0\,N\,r\}\subset B_{\xi,\mu_0r}\,.
\]
Let $x,y\in B_{\xi,r}$ and observe that Remark~\ref{rmkPN} implies
the existence of $a\in (P^N)^{-1}(D_1)$ such that $y=x\delta_{d(x,y)}P^N(a)$.
We fix $\lambda=d(x,y)$, then taking into account that $x\in B_{\xi,r}$
one easily realizes that 
\[
[0,\lambda\, a_s]\ni l\lra x\delta_\lambda P^{s-1}(a)\delta_lh_{i_s}\in A_{\xi,r}\subset
B_{\xi,\mu_0r}\subset\Omega
\]
for every $s=1,\ldots,N$.
Now, for every horizontal direction $h=\exp X\in\exp (V_1)$
and every $z\in\Omega$ such that
\[
[0,1]\ni t\lra \Gamma(t)=z\exp tX=z\delta_th\in\Omega
\] we have the
integral formula 
\[
F(xh)-F(x)=\int_0^1 DF(x\exp t X)(h)\,dt\,,
\]
since $t\lra F\circ \Gamma(t)$ is continuously differentiable with
$\big(F\circ\Gamma\big)'(t)=DF(x\exp tX)(h)$.
It follows that
\begin{eqnarray}\label{estFh}
\|F(xh)-F(x)\| \leq\|X\|\,\int_0^1\|dF(x\delta_th)\|\,dt\,.
\end{eqnarray}
For every $s=1,\ldots,N$, we apply \eqref{estFh} to the curve
\[
\Gamma(t)=x\delta_\lambda P^{s-1}(a)\delta_{t\lambda a_s}h_{i_s}\in D_{\xi,\mu_0r}\quad\mbox{for every}\quad t\in[0,1]
\]
obtaining the final estimates
\begin{eqnarray}
&&\|F(y)-F(x)\|\leq\sum_{s=1}^N
\|\exp^{-1}\big(\delta_{\lambda a_s}h_{i_s}\big)\|\;
\int_0^1\|dF(x\delta_\lambda P^{s-1}(a)\delta_{t\lambda a_s}h_{i_s})\|\,dt \\
&&\leq c_0\,d(x,y)\,\max_{1\leq s\leq N}\|\exp^{-1}\big(h_{i_s}\big)\|\,
\max_{y\in D_{\xi,\mu_0r}}\|dF(y)\|=C_0\,d(x,y)\,
\max_{y\in D_{\xi,\mu_0r}}\|dF(y)\|\,,\nonumber
\end{eqnarray}
that conclude the proof. $\Box$
\begin{Pro}\label{lipFcurve}
Let $c:[a,b]\lra\Omega$ be a Lipschitz curve
and let $f:\Omega\lra\M$ be a continuously P-differentiable mapping.
Then there exists a constant $C>0$, only depending on $d$,
$\rho$ and the norms on $\cG$ and $\cM$, such that
$\Gamma=f\circ c$ is Lipschitz and the following formulae hold
\begin{eqnarray}
&&\rho\big(\Gamma(t),\Gamma(\tau)\big)\leq C\;\mbox{\rm Lip}(c)
\max_{x\in c([a,b])} \|dF_1(x)\|\;|t-\tau|\,, \label{Liplip} \\
&& \dot\gamma_1(t)=dF_1\big(c(t)\big)\circ dc(t)=dF_1\big(c(t)\big)
\circ\dot\alpha_1(t) \label{dF_1}\,,
\end{eqnarray}
where we have set $\gamma_j=\pi_j\circ\exp^{-1}\circ\Gamma$, 
$\alpha_j=\pi_j\circ\exp^{-1}\circ c$,
\[
\Gamma=\exp\big(\gamma_1+\cdots+\gamma_\upsilon\big)\qquad
\mbox{and}\qquad c=\exp(\alpha_1+\cdots+\alpha_\iota)\,.
\]
\end{Pro}
{\sc Proof.}
Let $c:[a,b]\lra\Omega$ be a horizontal curve, where $f:\Omega\lra\M$ is 
continuously P-differentiable.
Due to Theorem~\ref{PdifPdifj}, $F_j$
are continuously P-differentiable, hence Lemma~\ref{PlipV}
yields their local Lipschitz property. 
Thus, the composition \[F\circ c:[a,b]\lra\cM\] is
locally Lipschitz, hence it is absolutely continuous.
By Corollary~\ref{hdifflip}, $c$ is horizontal, then
Theorem~\ref{pansu} yields the a.e. P-differentiability of $c$.
In view of the chain rule stated in Proposition~\ref{chain},
it follows that $\Gamma$, equal to $f\circ c$, is a.e. P-differentiable and
\[
d\Gamma(t)=df\big(c(t)\big)\big(dc(t)\big)\,.
\]
for a.e. $t\in[a,b]$.
Applying Corollary~\ref{Pansuimpl} and the characterization of Proposition~\ref{horprop}, it follows that $\Gamma$ is horizontal.
In view of \eqref{eqF_1}, the formula above yields the a.e.
validity of \eqref{dF_1}, since $c$ is horizontal.
Due to \eqref{dF_1} and \eqref{lipestim}, it follows that
$\Gamma$ is Lipschitz and formula \eqref{Liplip} follows. $\Box$
\begin{Cor}\label{corlip}
Every continuously P-differentiable mapping $f:\Omega\lra\M$ is locally Lipschitz.
\end{Cor}
{\sc Proof.}
It suffices to observe that $\Omega$ is locally connected by
geodesics with respect to a Carnot-Carath\'eodory distance
fixed on $\G$ and then apply Proposition~\ref{lipFcurve}. $\Box$
%
%
%
%
%
%
%
%
%
%
%
%
%
%
%
%
%
%
%
\section{Mean value inequality}\label{MVI}
Throughout this section, $\G$ and $\cG$ denote a stratified group
along with its Lie algebra, $\Omega\subset\G$ is a fixed open set, $d$ is a homogeneous distance on $\G$. $\M$ is a graded group equipped with a homogeneous distance $\rho$. When the norm $\|\cdot\|$ is applied to
a vector of either $T\G$ or $T\M$, it is understood that it represents
the left invariant Finsler norm generated by the fixed norm in the
Lie algebra of either $\G$ or $\M$. The mapping $f:\Omega\lra\M$ will be 
represented in the algebra by $F=\exp^{-1}\circ f:\Omega\lra\cM$ and 
$F_j=\pi_j\circ F$, where $\pi_j:\cM\lra W_j$ are the canonical projections onto
the layers.
\begin{The}\label{unifestimth}
Let $f:\Omega\lra\M$ be a continuously P-differentiable mapping and let  $\Omega'$ be an open subset compactly contained in $\Omega$. Let $c:[a,b]\lra\Omega'$ be a Lipschitz curve with $\|\dot c(t)\|\leq1$ a.e.
Then there exists a constant $C$, depending on the modulus of continuity of $dF_1$ on $\overline{\Omega'}$ and on $L=\max_{x\in\overline{\Omega'}}\|dF_1(x)\|$,
such that at approximate continuity points $t$ of $\dot c$
the following estimate holds
\begin{eqnarray}\label{unifestimch}
&&\quad\rho\left(Df\lls c(t)\rls\big(\exp \dot c(t)\,h\big),f\lls c(t)\rls^{-1}f\lls c(t+h)\rls\right)\\
&&\leq C\Bigg[L^{1/\iota}
\sup_{s\in I_{t,t+h}}\left|\,\medint_t^{t+h}\|\dot c(l)-\dot c(t)\|\,dl\,\right|^{1/\iota}
\!\!\!+\!\max_{s\in I_{t,t+h}}\|dF_1\lls c(s)\rls-dF_1\lls c(t)\rls\|^{1/\iota}
\Bigg]\!|h|,
\nonumber
\end{eqnarray}
where we have set
$I_{t,t+h}=[\min\{t,t+h\},\max\{t,t+h\}]$.
\end{The}
{\sc Proof.}
By Proposition~\ref{lipFcurve}, it follows that
$\Gamma=f\circ c=\exp\circ \gamma$ is Lipschitz and 
\begin{eqnarray}\label{gammalpha}
\dot\gamma_1(t)=dF_1(c(t))\big(\dot\alpha_1(t)\big)\,,
\end{eqnarray}
where $df=\exp^{-1}\circ Df\circ \exp$ and
$c=\exp\big(\alpha_1+\cdots+\alpha_\iota)$.
Observing that $\|\dot\alpha_1\|\leq 1$ a.e., we then
apply Theorem~\ref{pansu} to $\Gamma$.
As a consequence, taking into account \eqref{gammalpha},
the estimate \eqref{pansuestim} applied to $\Gamma$ yields
\begin{eqnarray}\label{pansuestim1}
&&\left\|\delta_{1/h}
\Big[-h\,dF_1\lls c(t)\rls\,\big(\dot\alpha_1(t)\big)\gope\big(-F\lls c(t)\rls\gope F\lls 
c(t+h)\rls\Big]\right\|\\
&&\leq\Upsilon(L)\;\sup_{s\in I_{t,t+h}}\left|\,\medint_t^s
\|dF_1\big(c(l)\big)\big(\dot\alpha_1(l)\big)-dF_1\big(c(t)\big)
\big(\dot\alpha_1(t)\big)\|dl\right|
\nonumber \\
&&\leq\Upsilon(L)\;
\Bigg[L\sup_{s\in I_{t,t+h}}\left|\,\medint_t^s
\|\dot\alpha_1(l)-\dot\alpha_1(t)\|\,dl\,\right|
+\max_{s\in I_{t,t+h}}\|dF_1\lls c(s)\rls-dF_1\lls c(t)\rls\|\Bigg]\nonumber
\end{eqnarray}
where $\Upsilon$ is the a nondecreasing function defined in
Theorem~\ref{pansu}. 
Now, applying \eqref{rhoestnorm} to \eqref{pansuestim1},
we obtain a constant $C>0$ such that \eqref{unifestimch}
holds. It is not difficult to notice that $C$ only depends on the
quantities stated in our claim. $\Box$
\begin{Cor}
Let $f:\Omega\lra\M$ be a continuously P-differentiable
mapping and let $\Omega_1$, $\Omega_2$ be open subsets
such that $\Omega_2$ is compactly contained in $\Omega$
and for every $x,y\in\Omega_1$ such that $x^{-1}y\in\exp(V_1)$
we have $x\delta_l(x^{-1}y)\in\overline{\Omega}_2$ whenever $0\leq l\leq1$.
Then there exists a constant $C$, depending on 
the modulus of continuity of $\omega_{\overline{\Omega}_2,dF_1}$
and on $\max_{x\in\overline{\Omega_2}}\|dF_1(x)\|$ such that
\begin{equation}\label{unifestimh}
\rho\left(Df\lls x\rls\big(x^{-1}y\big),
f\lls x\rls^{-1}f\lls y\rls\right)\leq C\;
\omega_{\overline{\Omega}_2,dF_1}\big(d(x,y)\big)^{1/\iota}\,d(x,y)\,.
\end{equation}
whenever $x,y\in\Omega_1$ and $x^{-1}y\in\exp V_1$.
\end{Cor}
{\sc Proof.}
It suffices to apply Theorem~\ref{unifestimth} to the curve
$c(l)=x\delta_lu$, with $[a,b]=[0,\rho(x,y)]$, $t=0$, $h=d(x,y)$
and $u=\delta_{1/\rho(x,y)}(x^{-1}y)$. Then
\begin{equation*}
\rho\left(Df\lls x\rls\big(x^{-1}y\big),
f\lls x\rls^{-1}f\lls y\rls\right)\leq C\;\Bigg[
\max_{l\in I_{t,t+h}}\|dF_1\lls c(l)\rls-dF_1\lls x\rls\|^{1/\iota}
\Bigg]\,d(x,y)\,.
\end{equation*}
Observing that $d(c(l),x)\leq d(x,y)$,
estimate \eqref{unifestimh} follows. $\Box$
\vskip.2truecm
{\sc Proof of Theorem~\ref{unifestim}.}
Let $x,y\in\Omega_1$ and let $a\in(P^N)^{-1}(D_1)$ be such that
$P^N(a)=\delta_{1/d(x,y)}\big(x^{-1}y\big)$.
Thus, due to Definition~\ref{PN}, we have
\[
d(P^s(a))\leq N\,c(\G,d)
\]
where $c(\G,d)$ and $N$ are defined in \eqref{defca}
and in Lemma~\ref{140FS}, respectively.
For the sequel, we set $c_0=c(\G,d)$.
Analogously, we have
\[
d\big(\delta_tP^{s-1}(a)\delta_lh_{i_s}\big)\leq c_0\,N\,\diam(\Omega_1)
\]
if $l\in [0,ta_s]$ and $0\leq t\leq\diam(\Omega_1)$.
We fix $t_0=d(x,y)\leq\diam(\Omega_1)$, hence
from our assumptions it follows that each curve
\[
[0,t_0a_s]\ni l\lra x\delta_{t_0}P^{s-1}(a)\delta_lh_{i_s}
\]
is contained in $\Omega_2$.
Applying \eqref{unifestimch} to this curve, there exists a constant
$C_1$, only depending on both $\omega_{\overline\Omega_2,dF_1}$
and $\max_{x\in\overline{\Omega_2}}\|dF_1(x)\|$, such that
\begin{eqnarray}
&&	\rho\Big(f(x\delta_{t_0}P^{s-1}(a))^{-1}f\big(x\delta_{t_0}P^s(a)\big),
	Df(x\delta_{t_0}P^{s-1}(a))\big(\delta_{a_st_0}h_{i_s}\big)\Big) \\
&&\leq  C_1\;|a_s\,t_0|\;
	\omega_{\overline{\Omega}_2,dF_1}\big(c_0\,t_0\big)^{1/\iota}.\nonumber
\end{eqnarray}
Then we obtain the key uniform estimate
\begin{eqnarray}\label{keyunif}
&&	\rho\Big(
	\delta_{1/t_0}\lLs 
	f(x\delta_{t_0}P^{s-1}(a))^{-1}f\big(x\delta_{t_0}P^s(a)\big)
	\rLs,Df(x\delta_{t_0}P^{s-1}(a))\big(\delta_{a_s}h_{i_s}\big)\Big) \\
&&  \leq  C_1\,c_0\;
	\omega_{\overline{\Omega}_2,dF_1}\big(c_0\,t_0\big)^{1/\iota}.\nonumber
	\end{eqnarray}
We fix $\nu=\max_{x\in\overline{\Omega_2}}\|dF_1(x)\|$, so that 
\eqref{rhonormiota} yields a constant $K_\nu>0$ such that
$$
\rho\big(Df(x)\big)\leq k_\nu\|dF_1(x)\|^{1/\iota}
$$
where we have set $\rho(Df(x))=\max_{d(h)=1}\rho\big(Df(x)(h)\big)$.
We have the estimate
\begin{eqnarray}\label{diffest}
	\rho\big(Df(x)(\delta_{a_s}h_{i_s}\cdots\delta_{a_N}h_{i_N})\big)\leq
	c_0\,k_\nu\,N\,\nu^{1/\iota} 
\end{eqnarray}
In view of \eqref{keyunif} and \eqref{diffest}, we use \eqref{estimprod}, getting a constant $C_2$, only depending on $\G$, $\nu$ and $\omega_{\overline{\Omega},dF_1}$ such that
\begin{eqnarray}\label{firstmult}
&&\rho\Big(\delta_{1/t_0}\lLs f(x)^{-1}f\big(x\delta_{t_0}P^N(a)\big)\rLs,
	Df(x)\big(P^N(a)\big)\Big)\\
&&\leq C_2 \sum_{s=1}^N	
  \rho\Big(\delta_{1/t_0}
	\lLs f(x\delta_{t_0}P^{s-1}(a))^{-1}f\big(x\delta_{t_0}P^s(a)\big)\rLs,
	Df(x)\big(\delta_{a_s}h_{i_s}\big)\Big)^{1/\iota}\,,\nonumber
\end{eqnarray}
where we have used the decomposition 
\begin{eqnarray}\label{hordec}
\delta_{1/t_0}\lLs f(x)^{-1}f\big(x\delta_{t_0}P^N(a)\big)\rLs=	\prod_{s=1}^N\delta_{1/t_0}\lLs
f\big(x\delta_{t_0}P^{s-1}(a)\big)^{-1}f\big(x\delta_{t_0}P^s(a)\big)\rLs.
\end{eqnarray}
By \eqref{keyunif} we have
\begin{eqnarray}\label{previntest}
&&\rho\Big(\delta_{1/t_0}
	\lLs f\big(x\delta_{t_0}P^{s-1}(a)\big)^{-1}f\big(x\delta_{t_0}P^s(a)\big)\rLs,
	Df(x)\big(\delta_{a_s}h_{i_s}\big)\Big)^{1/\iota}\\
&&\leq \big(C_1\,c_0\big)^{1/\iota}\;
	 \omega_{\overline{\Omega}_2,dF_1}\big(c_0\,t_0\big)^{1/\iota^2}
	 +\rho\Big(Df(x\delta_{t_0}P^{s-1}(a)(\delta_{a_s}h_{i_s}),
	 Df(x)\big(\delta_{a_s}h_{i_s}\big)\Big)^{1/\iota}\nonumber
\end{eqnarray}
By \eqref{rhoestnorm} we have a constant $C_3>0$ depending on 
$c$ and $\nu$, such that
\begin{eqnarray}\label{intermest}
	\rho\Big(Df(x\delta_{t_0}P^{s-1}(a))(\delta_{a_s}h_{i_s}),
	 Df(x)\big(\delta_{a_s}h_{i_s}\big)\Big)\leq c_0\;C_3\;
	 \omega_{\overline{\Omega_2},dF_1}(N\,c_0\,t_0)^{1/\iota}
\end{eqnarray}
Joining \eqref{firstmult}, \eqref{previntest} and \eqref{intermest},
we achieve
\begin{eqnarray*}
&&\rho\Big(\delta_{1/t_0}\lLs f(x)^{-1}f\big(x\delta_{t_0}P^N(a)\big)\rLs,
	Df(x)\big(P^N(a)\big)\Big)\\
&&\leq N\,C_2\Big(c_0^{1/\iota}\,C_1^{1/\iota}
\omega_{\overline{\Omega}_2,dF_1}(c_0t_0)^{1/\iota^2}
  +c_0^{1/\iota}\,C_3^{1/\iota}\omega_{\overline\Omega_2dF_1}
  (N\,c_0\,t_0)^{1/\iota^2}\Big)\,.\nonumber
\end{eqnarray*}
Taking into account that $y=x\delta_{t_0}P^N(a)$,
then the previous estimate can be precisely written as follows
\begin{eqnarray*}
&&\frac{\rho\Big(f(x)^{-1}f(y),Df(x)(x^{-1}y)\Big)}{d(x,y)}\\
&&\leq N\,C_2\,c_0^{1/\iota}
\Big(C_1^{1/\iota}\omega_{\overline\Omega_2,dF_1}(c_0\,d(x,y))^{1/\iota^2}
  +C_3^{1/\iota}\omega_{\overline\Omega_2dF_1}
  (N\,c_0\,d(x,y))^{1/\iota^2}\Big)\,.\nonumber
\end{eqnarray*}
The estimate obtained in $\Omega_1$ easily extends to its closure.
This concludes the proof. $\Box$

\subsection{Inverse mapping theorem}\label{invmapthe}
As an immediate consequence of the mean value inequality we achieve
the inverse mapping theorem for P-differentiable mappings. 
\begin{The}\label{inverse}
Let $\Omega$ be an open subset of $\G$ and let $f:\Omega\lra\G$ be a continuously 
P-differentiable mapping and assume that $Df(x):\G\lra\G$ is invertible at some
$\overline{x}\in\Omega$. Then there exists a neighbourhood $U$
of $\overline{x}$ such that the restriction $f_{|U}$ has an inverse mapping
$g$ that is also P-differentiable and for every $y\in f(U)$ we have
$$
Dg(y)=Df(g(y))^{-1}.
$$
\end{The}
{\sc Proof.}
By continuity of $x\lra Df(x)$, there exists a compact
neighbourhood $U'$ of $\overline{x}$ and a number $\mu>0$
such that
\[
\min_{x\in U'}\min_{d(u)=1}\rho\big(Df(x)(u)\big)=\mu.
\]
Then triangle inequality and estimates \eqref{keyunifest}
give us a constant $\beta>0$ and an open neighbourhood
$U$ contained in $U'$ and depending on the modulus
of continuity $\omega_{df,U'}$, such that
\[
d\big(f(x),f(y)\big)\geq \beta\;d(x,y)
\]
whenever $x,y\in U$. Then $f_{|U}$ is invertible onto its image
$f(U)$ and by the domain invariance theorem,
see for instance Theorem~3.3.2 of \cite{LLoyd}, 
$f_{|U}$ is an open mapping. Then the inverse mapping
$g:f(U)\lra U$ is continuous and the chain rule of
Proposition~\ref{chain} concludes the proof. $\Box$

%
%
%
%
%
%
%
%
%
%
%
%
%
%
%
%
%
\section{Homogeneous subgroups}\label{homsubgr}
Throughout this section, we assume that all subgroups are closed,
connected and simply connected Lie subgroups, if not otherwise stated.
Here $\G$ denotes an arbitrary graded group, that is not necessarily stratified. Its Lie algebra is denoted by $\cG$.
From 5.2.4 of \cite{Stein}, we recall the notion of homogeneous subalgebra and the corresponding notion of homogeneous subgroup.
%
%
%
%
%
%
%
%
%
%
%
\begin{Def}[Homogeneous subalgebra]{\rm
Let $\cp\subset\cG$ be a Lie subalgebra. We say that $\cp$ is a
{\em homogeneous subalgebra} if $\delta_r\cp\subset\cp$ for every $r>0$.
}\end{Def}
\begin{Rem}{\rm
It is not difficult to find examples of subalgebras
which are not homogeneous. It suffices to consider
$\cL=\span\{X+Z\},$ which is a subalgebra of the Heisenberg
algebra $\ch^1$ of brackets $[X,Y]=Z$.
However, $\delta_2(X+Z)=2X+4Z\notin\cL$.
}\end{Rem}
\begin{Pro}\label{pfactrs}
Let $\ca$ be a homogeneous subalgebra of $\cG$, where $\cG$
is decomposed into the direct sum $V_1\oplus\cdots\oplus V_\iota$.
Then we have
$$
\ca=(V_1\cap\ca)\oplus (V_2\cap\ca)\oplus\cdots\oplus(V_\iota\cap\ca)\,.
$$
\end{Pro}
{\sc Proof.}
We have only to prove that each vector
$\xi\in\ca$ with the unique decomposition $\xi=\sum_{j=1}^\iota\xi_j$,
with $\xi_j\in V_j$, satisfies the property $\xi_j\in\ca$.
By hypothesis, $\delta_r\xi\in\ca$ whenever $r>0$, then
the closedness of $\ca$ implies that
$$
\lim_{r\ra0^+}\frac{1}{r}\;\delta_r\xi=\xi_1\in\ca\,.
$$
This implies that $\xi-\xi_1\in\ca$, then 
$$
\lim_{r\ra0^+}\frac{1}{r^2}\;\delta_r(\xi-\xi_1)=\xi_2\in\ca\,.
$$
Iterating this process by induction, our claim follows. $\Box$
\begin{Cor}\label{homsub}
Every homogeneous subalgebra $\ca\subset\cG$ is a graded algebra.
\end{Cor}
\begin{Rem}{\rm
Let $\ca=\ca_1\oplus\cdots\oplus\ca_\iota$ be a homogeneous
subalgebra. In general, some factor $\ca_j$ might be the null space.
Consider for instance the homogeneous subalgebra $\ca=V_2\oplus V_4\oplus\cdots\oplus V_{2[\iota/2]}$, where 
$\cG=V_1\oplus V_2\oplus\cdots\oplus V_\iota$.
}\end{Rem}
\begin{Def}[Homogeneous subgroup]{\rm
Let $P\subset\G$ be a Lie subgroup. We say that $P$ is a
{\em homogeneous subgroup} if $\delta_rP\subset P$ for every $r>0$.
}\end{Def}
It is clear that all properties of homogeneous subalgebras are exactly
translated to homogeneous subgroups.
Thus, in the sequel we will equivalently work with either
homogeneous subalgebras or homogeneous subgroups.

\subsection{Complementary subgroups}\label{complemsbsect}
This subsection is devoted to the notion of complementary
subgroup, that plays an important role in the algebraic side of this paper.

When  $A$ and $B$ are subsets of an abstract group $G$ we will use the notation
\[AB=\{ab\,|\, a\in A, b\in B\}.\]
\begin{Def}[Complementary subgroup]\label{complem}{\rm
Let $P$ be a homogeneous subgroup of $\G$.
If there exists a homogeneous subgroup $H$ of $\G$
satisfying the properties $PH=\G$ and $P\cap H=\{e\}$,
then we say that $H$ is a {\em complementary subgroup} to $P$.
}\end{Def}
\begin{Rem}{\rm
Recall that for abstract subgroups $A,B$ of an abstract group $G$,
the subset $AB$ is an abstract subgroup if and only if $AB=BA$,
see for instance \cite{Her}. As a consequence,
$H$ is complementary to $P$
if and only if $P$ is complementary to $H$.
}\end{Rem}
Due to the previous remark, being complementary is a symmetric relation
and we can say that two subgroups are
{\em complementary}.
\begin{Rem}\label{uniqued}
{\rm Let $H$ and $P$ be complementary subgroups of $\G$
and let $g\in\G$. Then it is immediate to check that
there exist unique elements $h,h'\in H$
and $p,p'\in P$ such that $g=hp=p'h'$.
}\end{Rem}
The following proposition characterizes complementary subgroups
through their Lie subalgebras.
\begin{Pro}\label{homdec}
Let $\cp$ and $\ch$ be homogeneous subalgebras of $\cG$
and let $P$ and $H$ denote their corresponding homogeneous subgroups,
respectively. Then the condition $\cp\oplus\ch=\cG$ is equivalent
to require that $P$ and $H$ are complementary subgroups.
Furthermore, if one of these conditions hold, then the mapping
\begin{eqnarray}\label{phidecomp}
\phi:\cp\times\ch\lra\G,
\qquad \phi(W,Y)=\exp W\;\exp Y
\end{eqnarray}
is a diffeomophism.
\end{Pro}
{\sc Proof.}
Assume that $\cp\oplus\ch=\cG$.
By Theorem~2.10.1 of \cite{Vara} the exist bounded
neighbourhoods $\Omega_1$ and $\Omega_2$ of the origin
in $\cp$ ad in $\ch$, respectively, such that the mapping
restricted to $\Omega_1\times\Omega_2$
is a diffeomorphism onto some neighbourhood $U$ of $e$.
Observing that
\[
\delta_r\big(\phi(W,Y)\big)=\phi(\delta_rW,\delta_rY)
\]
for every $r>0$, it follows that $\phi$ is invertible.
This shows that $P$ and $H$ are complementary subgroups.
Viceversa, suppose that $P$ and $H$ are complementary subgroups
and assume by contradiction that there exists $X\in\cp\cap\ch\sm\{0\}$.
This implies that $\exp X\in H\cap P=\{e\}$, that is absurd.
Now, we wish to show that $\cp\oplus\ch=\cG$.
By contradiction, we assume that $\dim\big(\cp\oplus\ch\big)<\dim\cG$.
We notice that the mapping $\phi$ defined in \eqref{phidecomp}
has everywhere injective differential, since $\ch\cap\cp=\{0\}$
and both $\exp_{|\cp}$ and $\exp_{|\ch}$ have everywhere injective
differential.
Thus, $\phi(\cp\times\ch)\subset\G$ is an open subset,
whose dimension is less than $\dim\cG$.
As a consequence, $\phi$ cannot be surjective and this conflicts
with hypothesis $PH=\G$. $\Box$
\begin{Rem}{\rm
Joining Proposition~\ref{pfactrs} and Proposition~\ref{homdec},
we have the following property.
If $K$ and $P$ are complementary subgroups, then 
\begin{eqnarray}\label{QKP}
Q=\mbox{$\mcH$-$\dim(K)$}+\mbox{$\mcH$-$\dim(P)$},
\end{eqnarray}
where $\mbox{$\mcH$-$\dim$}$ denotes the Hausdorff dimension
with respect to a fixed homogeneous distance.
One can interpret this fact as a ``proper splitting" of the ambient
group even with respect to the metric point of view.
}\end{Rem}
Through Proposition~\ref{homdec},
it is easy to observe that not any subgroup admits a complementary
subgroup.
\begin{Exa}\label{nnexcmpl}{\rm
Let us consider the second layer $\cn=\span\{Z\}$
of the Heisenberg algebra $\ch^1$, with bracket relations $[X,Y]=Z$.
Then the normal subgroup $N=\exp(\cn)$ does not possess any
complementary subgroup, see also \cite{FSSC6}.
In fact, let $\ca$ be any 2-dimensional homogeneous subalgebra of $\ch^1$.
Then Proposition~\ref{pfactrs} gives the decomposition $\ca=\ca_1\oplus\ca_2$, 
where $\ca_j$ are contained in the $j$-th layer, $j=1,2$.
If $\ca_2=\{0\}$, then $\ca=\ca_1=\cv$ is not a subalgebra of
$\ch^1$. Thus, $\ca$ must contain $\cn$ and this conflicts
with existence of a 2-dimensional subalgebra complementary to $\cn$.
}\end{Exa}
\begin{Rem}{\rm More generally there do not
exist 2-dimensional subalgebras $\ca$ of $\ch^1$ such that $\ca\oplus\cn=\ch^1$,
even if we do not require the homogeneity of $\ca$.
In fact, let $\ca=\span\{U+\alpha Z,U'+\alpha'Z\}$ be a 2-dimensional
subalgebra of $\ch^1$, where $U,U'$ belong to the first layer $\cv$.
If $[U,U']=0$ then they are proportional and $\ca\oplus\cn$ is 2-dimensional. It follows that $[U,U']=\gamma Z$, with $\gamma\neq0$.
On the other hand, $\ca$ is a subalgebra, hence $Z\in\ca$
and this would imply $\ca=\ch^1$.
This is also a contradiction, then such a complementary subalgebra
cannot exist.
}\end{Rem}
We also point out that a complementary subgroup, when it exists
need not be unique, as we show in the next
\begin{Exa}{\rm
All the homogeneous subalgebras
$$
\ca_\lambda=\span\{X+\lambda Y\}\qquad \lambda\in\R,
$$
of $\H^1$ yield complementary subgroups $A_\lambda=\exp\ca_\lambda$
of $S=\exp\big(\span\{Y,Z\}\big).$
}\end{Exa}
Next, we wish to see how the notion of complementary subgroup
is translated for the corresponding subalgebras, when
$\G$ is a graded group.
First of all, we notice that a decomposition of $\cG$ as a direct sum of
two subspaces does not ensure a corresponding decomposition
of $\G$ as the product of their images through the exponential mapping.
This simple fact is shown in the next example.
\begin{Exa}{\rm
We consider the Heisenberg algebra $\ch^1$ with basis
$(X,Y,Z)$ and bracket relation $[X,Y]=Z$.
We consider the following subspaces of $\ch^1$
\[
\left\{\begin{array}{l}
\cu=\span\{X,Y\} \\
\cw=\span\{Y+\frac{1}{2}Z\}
\end{array}\right..
\]
It is immediate to check that $\cu\oplus\cw=\ch^1$,
but $\exp\cu\;\exp\cw\neq\H^1$, since
\[
\exp(-X+\lambda Z)\notin \exp\cu\;\exp\cw
\]
for every $\lambda\neq0$.
}\end{Exa}
Even the decomposition of $\cG$ as the direct sum
of two subalgebras does not imply that their corresponding subgroups
factorize $\G$, as we show in the next
\begin{Exa}{\rm
We consider the Heisenberg algebra $\ch^2$ with basis
$(X_1,Y_1,X_2,Y_2,Z)$ and bracket relations
$[X_1,Y_1]=[X_2,Y_2]=Z$.
We consider the following subalgebras of $\ch^2$
\[
\left\{\begin{array}{l}
\ca=\span\{X_1,X_2,Z+Y_1\} \\
\cb=\span\{Y_1,Y_2\}
\end{array}\right..
\]
It is immediate to check that $\ca\oplus\cb=\ch^2$,
but
\[
\exp(2X_1+Z)\notin AB,
\]
where $A=\exp\ca$ and $B=\exp\cb$ are the corresponding Lie subgroups.
}\end{Exa}
\subsection{h-epimorphisms and h-monomorphisms}\label{hepimono}
In this subsection we use complementary subgroups to 
characterize both h-epimorphisms and h-monomorphisms.
These characterizations will be used in the proof of both
the implicit function theorem and the rank theorem.
\begin{Pro}[Characterization of h-epimorphisms]\label{GlinH}
Let	$L:\G\lra\M$ be a surjective h-homomorphism and let $N$
be its kernel. The conditions
\begin{enumerate}
	\item there exists a subgroup $H$ complementary to $N$
	\item $L$ is an h-epimorphism
\end{enumerate}
are equivalent and if one of them holds, then the restriction $L_{|H}$ is an h-isomorphism.
\end{Pro}
{\sc Proof.}
We first assume the validity of the first condition.
Then we consider the restriction $T=L_{|H}:H\lra\M$. 
Let $m\in\M$ and let $g\in\G$ such that $L(g)=m$. 
Then the property $NH=\G$ implies that $g=nh$, where $(n,h)\in N\times H$.
As a consequence, we have
$$
L(nh)=L(h)=T(h)=m,
$$
hence $T$ is surjective.
If we have $T(h)=e$, then
$$
h\in H\cap N=\{e\}.
$$
We have shown that $L_{|H}$ is an h-homomorphism.
Clearly, $T^{-1}:\M\lra H$ is an h-homomorphism and
satisfies $L\circ T^{-1}=\Id_\M$, hence $L$ is an h-epimorphism.
Conversely, if $L$ is an h-epimorphism, then there exists
a right inverse $R:\M\lra\G$ that is also an h-homomorphism.
We set $H=R(\M)$ and easily observe that $H\cap N=\{e\}$.
Let $g\in\G$ and consider
\[
m=L(g)=L\big(R(m)\big),
\]
therefore $g^{-1}R(m)\in N$ and this implies that $g\in HN$.
It follows that $\G=HN=NH$ and this concludes the proof. $\Box$
\begin{Pro}[Characterization of h-monomorphisms]\label{GlinHmono}
Let	$T:\G\lra\M$ be an injective h-homomorphism and let $H$
be its image. The conditions
\begin{enumerate}
	\item there exists a normal subgroup $N$ complementary to $H$,
	\item there exists an h-epimorphism $p:\M\lra H$ such that
	$p_{|H}=\Id_H$,
	\item $T$ is an h-monomorphism
\end{enumerate}
are equivalent.
\end{Pro}
{\sc Proof.}
We show that (1) implies (2).
We define the projection $p:\M\lra H$ that associates
to any element $m=nh\in\M$, with $(n,h)\in N\times H$
the element $h\in H$. This definition is well posed, since 
$N$ and $H$ are complementary.
The fact that $N$ is normal and the uniqueness of representation
of the product $nh$ give
\[
p\big((n_1h_1)(n_2h_2)\big)=p\big((n_1h_1n_2h_1^{-1})(h_1h_2)\big)=h_1h_2
=p(n_1h_1)\,p(n_2h_2)\,.
\]
It is trivial to observe that $p$ is homogeneous and that
its restriction to $H$ is exactly $\Id_H$.
Then $p$ is a surjective h-homomorphism. Furthermore, $H$ is
a complementary subgroup of the kernel $N$, then Proposition~\ref{GlinH}
implies that $p$ is an h-epimorphism.

To prove that (2) implies (3), we define the mapping
$L=J\circ p$, where the injectivity of $T$ allows us to define
the h-isomorphism $J:H\lra\G$ such that $J\big(T(g)\big)=g$
for every $g\in\G$. Then $L$ is an h-homomorphism as composition
of h-homomorphisms. In addition, one can easily verify that $L\circ T=\Id_\G$.

We are left to show that (3) implies (1).
Bt definition of h-monomorphism we have an h-homomorphism
$L:\M\lra\G$ that is a left inverse of $T$. Let $N$ be its kernel, that is a normal homogeneous subgroup of $\M$.
We have to show that $N$ and $H$ are complementary subgroups.
Let $m\in\M$ and consider $L(m)=g\in\G$. We have
\[
L\big(m\,T(g^{-1})\big)=L(m)g^{-1}=e\,,
\]
then $m\,T(g^{-1})=n\in N$, namely, $m=n\,T(g)\in NH$.
Let $x\in N\cap H$ and let $g\in\G$ such that $x=T(g)$.
Thus, we get
\[
e_\G=L(x)=L\circ T(g)=g\,,
\]
that implies $g=e_\G$ and $x=T(e_\G)=e_\M$, namely, $N\cap H=\{e_\M\}$. $\Box$
%
%
%
%
%
%
%
%
%
%
%
%
%
%
%
%
\section{Quotients of graded groups}\label{hquot}

In this section we show that the group quotient of graded groups is still
graded with a natural left invariant and homogeneous distance.
We first recall some elementary facts of Lie groups Theory 
in order to study the relationship between quotients of
Lie algebras and quotients of Lie groups.

Let $G$ be a real Lie group and let $H$ be a Lie subgroup of $G$.
Recall that the quotient $G/H$ has a unique manifold structure
that makes the projection $\pi:G\lra G/H$ a smooth mapping.
$G/H$ is called homogeneous manifold, see Theorem~3.58 of \cite{Warn}.

If we consider a normal Lie subgroup $N$, then $G/N$ is in addition
a Lie group, according to Theorem~3.64 of \cite{Warn} and in this case 
\[
\pi:G\lra G/N
\]
is clearly is a Lie group homomorphism.
As a result, by Theorem~3.14 of \cite{Warn}
\[
d\pi:\cG\lra\overline\cG
\]
is a Lie algebra homomorphism, where 
$\cG$ and $\overline\cG$ are the Lie algebras of $G$ and
$G/N$, respectively.
\begin{Rem}{\rm The mapping $d\pi$ is also surjective.
In fact, Theorem~3.58 ensures the existence of a neighbourhood $W\subset G/N$
of the unit element $\overline{e}$ of $G/N$ and a smooth mapping
$\tau:W\lra G$ such that $\pi\circ\tau=\Id_W$.
Let $\Gamma:]-\ep,\ep[\lra W$ be a smooth curve with
$\dot\Gamma(0)=\overline{X}_{\overline e}\in T_{\overline e}\big(G/N\big)$.
Then the curve $\gamma=\tau\circ\Gamma:]-\ep,\ep[\lra G$ satisfies
$\pi\circ\gamma=\Gamma$ and $d\pi(X_e)=\overline{X}_{\overline e}$,
where $X_e=\dot\gamma(0)$. This shows surjectivity of $d\pi$.
}\end{Rem}
\begin{Pro}\label{algquot}
The Lie algebra $\overline{\cG}$ of $G/N$ is isomophic
to the quotient algebra $\cG/\cN$, where $\cG$ is the
Lie algebra of $G$ and $\cN$ is the ideal corresponding
to the normal subgroup $N$. Furthermore, we can represent
the exponential mapping of $G/N$ on $\cG/\cN$ by setting
$\Exp:\cG/\cN\lra G/N$ and
\[
\Exp(X+\cN)=\pi\big(\exp(X)\big)
\]
where $\exp:\cG\lra G$ is the exponential mapping of $G$,
$\pi:G\lra G/N$ is the canonical projection and the diagram
\[
\xymatrix@C=30pt@R=30pt{\cG/\cN \ar_{{\rm Exp}}[rd] \ar^{\overline{d\pi}}[r]
& \overline\cG  \ar^{\overline{\exp}}[d]  \\
  & G/N   }
\]
commutes, where $\overline{\exp}:\overline{\cG}\lra G/N$ is
the canonical exponential mapping of $G/N$.
\end{Pro}
{\sc Proof.}
We have seen above that $d\pi:\cG\lra\overline{\cG}$ is 
surjective, hence
\[
\dim\big(\ker(d\pi)\big)=\dim(\cN).
\]
For every $U\in\cN$ the curve $\pi\big(\exp(tX)\big)$ is constantly
equal to the unit element $\overline e$ of $G/N$,
then $d\pi(X)=0$ and clearly $\cN\subset\ker d\pi$.
It follows that $\cN=\ker d\pi$ and it is well defined
the algebra isomorphism $\overline{d\pi}:\cG/\cN\lra\overline{\cG}$,
where $p:\cG\lra\cG/\cN$ is the canonical projection and
\[
[p(X),p(Y)]=p\big([X,Y]\big)
\]
for every $X,Y\in\cG$. If we consider the diagram
\begin{eqnarray}\label{comdiag}
\xymatrix@C=30pt@R=30pt{\cG \ar_{p}[d] \ar^{d\pi}[rd] &   \\
\cG/\cN \ar_{{\rm Exp}}[rd] \ar^{\overline{d\pi}}[r] & \overline\cG
 \ar^{\overline{\exp}}[d]  \\
  & G/N   }
\end{eqnarray}
then we are left to show that its lower part commutes.
In fact, by definition of $\overline{d\pi}$,
for every $X\in\cG$ we have
\begin{eqnarray}
\overline{\exp}\circ\overline{d\pi}\big(p(X)\big)=
\overline{\exp}\circ d\pi(X).
\end{eqnarray}
On the other hand, due to Theorem~3.32 of \cite{Warn}, the diagram
\[
\xymatrix@C=30pt@R=30pt{G  \ar^{\pi}[r] & G/N \\
\cG \ar^{\exp}[u] \ar^{d\pi}[r] & \overline{\cG} \ar_{\overline{\exp}}[u]  }
\]
commutes, namely $\overline{\exp}\circ d\pi(X)=\pi\circ\exp(X)$
for every $X\in\cG$. It follows that
\[
\overline{\exp}\circ\overline{d\pi}=\pi\circ\exp=\Exp.
\]
This concludes the proof. $\Box$
\begin{Pro}\label{quotgrad}
Let $\M$ be a graded group of algebra 
$\cM=W_1\oplus\cdots\oplus W_\upsilon$
and let $N$ be a normal homogeneous subgroup of corresponding ideal 
$\cN=\cN_1\oplus\cdots\oplus\cN_\upsilon$.
Then the surjective Lie algebra homomorphism $d\pi:\cM\lra\overline\cM$
induces a grading
\[
\overline W_1\oplus\cdots\oplus \overline W_\upsilon
\]
on the Lie algebra $\overline\cM$ of $\M/N$, that is 
h-isomorphic to
\[
W_1/\cN_1\oplus\cdots\oplus W_\upsilon/\cN_\upsilon\,,
\]
where $d\pi(W_i)=\overline W_i$ is linearly isomorphic to $W_i/\cN_i$.
Furthermore, a one parameter group of dilations $\overline\delta_r$
can be defined on $\M/N$ such that $\pi:\M\lra\M/N$ is an h-homomorphism.
Finally, if $\M$ is stratified then so is $\M/N$.
\end{Pro}
{\sc Proof.}
The decomposition of $\cN$ into factors $\cN_i$ follows
from Corollary~\ref{homsub}, where some of $\cN_i$'s are possibly vanishing.
By Proposition~\ref{algquot} the Lie algebra $\overline\cM$ of the quotient
$\M/N$ is isomorphic to $\cM/\cN$ and
\[
p:\cM\lra\cM/\cN\qquad\mbox{is equivalent to}\qquad
d\pi:\cM\lra\overline\cM\,,
\]
by commutativity of diagram \eqref{comdiag},
hence we consider the grading induced by $p$.
Let us check that $p(W_i)\cap p(W_j)=\{0\}$ whenever $i\neq j$.
By contradiction, let $\overline X\in p(W_i)\cap p(W_j)$ be nonvanishing.
Then there exist $X_i\in W_i\sm\cN$ and $X_j\in W_j\sm\cN$
such that $X_i=X_j+Z$, where $Z\in\cN$.
Since $\cN$ is homogeneous, we must have $Z=Z_i+Z_j$, where
$Z_i\in W_i\cap\cN$ and $Z_j\in W_j\cap\cN$, hence 
$Z_i=X_i$ and $Z_j=-X_j$ this conflicts with our initial assumption
on $X_i$ and $X_j$. This shows that
$\overline\cM=\overline W_1\oplus\cdots\oplus\overline W_\upsilon$ is a grading,
where $\overline W_i=d\pi(W_i)$. This immediately implies that $\overline \cM$
is stratified when so is $\cM$. If we set $\overline\delta_r(\overline X_i)=r^i\overline X_i$ for every $\overline X_i\in\overline W_i$, then one can
easily check that $d\pi\circ\delta_r=\overline\delta_r\circ d\pi$ and
$\pi$ is an h-homomorphism. 
Finally, we consider the following commutative diagram
\[
\xymatrix@C=30pt@R=30pt{\cM \ar_{p_0}[d] \ar^{p}[r]
& \cM/\cN   \\
W_1/\cN_1\oplus\cdots\oplus W_\upsilon/\cN_\upsilon \ar^{J}[ru] & }
\]
where we have set
\[
p_0(w_1+\cdots+w_\upsilon)=\big(w_1+\cN_1,\ldots,w_\upsilon+\cN_\upsilon\big)
\]
and the following definition
\[
J\big(w_1+\cN_1,\ldots,w_\upsilon+\cN_\upsilon\big)=
w_1+\cdots+ w_\upsilon+\cN
\]
is well posed. It is immediate to observe that $\ker J=\{0\}$,
then $J$ is a linear isomorphism. The commutativity of the diagram above
implies that
\[
p(W_i)=J\big(W_i/\cN_i\big)\,.
\]
This concludes the proof. $\Box$
\begin{Rem}{\rm 
The induced grading on $\cM/\cN$ considered on Proposition~\ref{quotgrad} is 
given by the subspaces $\overline W_j=d\pi(W_j)$, that satisfy
$\cM/\cN=\overline W_1\oplus\cdots\oplus\overline W_\upsilon$, where some of 
these factors are possibly vanishing.
On the other hand, when $\cN$ is not homogeneous $\cM/\cN$
still can have a grading with dilations, but these ones
do not commute with $\pi$, namely, $\pi$ is no longer an
h-homomorphism.
}\end{Rem}
\begin{Exa}{\rm
Let $W_1=\span\{X_1,X_2,X_3\}$ and $W_2=\span\{Z\}$ where the
only nontrivial bracket relation is $[X_1,X_2]=Z$.
$\cM=W_1\oplus W_2$ is a 2-step stratified algebra and
$\cN=\span\{X_3-Z\}$ is an ideal of $\cM$ that is clearly not homogeneous.
Elements of $\cM/\cN$ are $\overline X_i=X_i+\cN$, where $i=1,2,3$
and $\overline X_3=Z+\cN$. Clearly, the only nontrivial bracket
relation is $[\overline X_1,\overline X_2]=\overline X_3$ and
$\cM/\cN$ is isomorphic to the 3-dimensional Heisenberg algebra.
On the other hand, there are no dilations $\overline\delta_r$ on
$\cM/\cN$ such that $d\pi\circ\delta_r=\overline\delta_r\circ d\pi$, since
\[
d\pi\big(\delta_rX_3\big)=r\,d\pi(X_3)=r\,\overline X_3\qquad
\mbox{and}\qquad
d\pi\big(\delta_rZ\big)=r^2\,d\pi(Z)=r^2\,\overline X_3.
\]
We also notice that in this case $d\pi$ cannot induce any grading
on $\overline\cM$, since $d\pi(W_1)\cap d\pi(W_2)=\{\overline X_3\}$.
}\end{Exa}
%
%
%
%
%
%
%
%
%
%
%
%
%
%
%
%
\section{Implicit Function Theorem and Rank Theorem}\label{MProofs}
In this section we prove both the implicit function theorem
and the rank theorem stated in the introduction.
We will use the following notation
$D_{\overline{n},r}^N=D_{\overline{n},r}\cap N$ and
$D_{\overline{h},s}^H=D_{\overline{h},s}\cap H$.
%
%
%
%
%
%
%
%
%
%
%
%
\vskip.2truecm
\noindent
{\sc Proof of Theorem~\ref{implth}.}
Our arguments are divided into two main steps.

{\bf Step 1}: {\em Existence.}

Due to Proposition~\ref{GlinH}, the restriction
$Df(\overline{x})(h):H\lra\G$ is invertible, then
$$
\min_{\substack{d(h)=1\\ h\in H}}
\rho\lls Df(\overline{x})(h)\rls>0.
$$
By continuity of
$
\Omega\times H\ni (x,h)\lra\rho\lls Df(x)(h)\rls,
$
there exists $R>0$ such that 
\begin{equation}\label{mumu}
\mu=\min_{x\in D_{\overline{x},R}}\min_{\substack{d(h)=1\\ h\in H}}
\rho\lls Df(x)(h)\rls>0\quad\mbox{and}\quad
\overline{\mu}=\max_{x\in D_{\overline{x},R}}\max_{\substack{d(h)=1\\ h\in H}}
\rho\lls Df(x)(h)\rls>0
\end{equation}
with $D_{\overline{x},R}\subset\Omega$.  
Now we set $R_0=R/(2+2\,c(\G,d)\,N)$,
where $c(\G,d)$ and $N$ are as \eqref{defca} and
Lemma~\ref{140FS}, respectively.
Then condition \eqref{inclomega1} is satisfied and 
Theorem~\ref{unifestim} yields the estimate
\begin{equation}\label{applunifestim}
	\max_{x,y\in D_{\overline{x},R_0}}
	\rho\lls f(x)^{-1}f(y),Df(x)(x^{-1}y)\rls\leq
	C\;\;\big[\omega_{D_{\overline{x},R},df}\big(N\,c\,d(x,y)\big)
\big]^{1/\iota^2}\;d(x,y)\,.
\end{equation}
Thus, by definition of $\overline{\mu}$ in \eqref{mumu} we have
\begin{eqnarray}\label{frabove}
\rho\big(f(x),f(y)\big)\leq
\left(C\big[\omega_{D_{\overline{x},R},df}\big(N\,c\,d(x,y)\big)
\big]^{1/\iota^2}+\overline{\mu}\right)\;d(x,y)
\end{eqnarray}
for every $x,y\in D_{\overline{x},R_0}$.
Let $(\overline{n},\overline{h})$ be the unique couple of $N\times H$ such that $\overline{n}\overline{h}=\overline{x}$ and let
$r,s>0$ be such that
$D^N_{\overline{n},r}D^H_{\overline{h},s}\subset D_{\overline{x},R_0}$.
Then, using definition of $\mu$ in \eqref{mumu},
for every $n\in D^N_{\overline{n},r}$
and every $h,h'\in D^H_{\overline{h},s}$ there holds
\begin{eqnarray}\label{frbelow}
\rho\big(f(nh),f(nh')\big))\geq
\left(\mu-C\big[\omega_{D_{\overline{x},R},df}\big(N\,c\,d(h,h')\big)
\big]^{1/\iota^2}\right)\;d(h,h')
\end{eqnarray}
observing that $d(nh,nh')=d(h,h')$.
Using \eqref{frabove} and \eqref{frbelow} and possibly taking
a smaller $s>0$ depending on $C, N, c>0$ we get a constant $\beta>0$
such that
\begin{eqnarray}\label{biLip}
\beta^{-1}\;d(h,h')\leq\rho\big(F_n(h),F_n(h')\big)\leq\beta\;d(h,h')
\end{eqnarray}
where for each $n\in D^N_{\overline{n},r}$ we have defined
$$
F_n:D^H_{\overline{h},s}\lra \M,\qquad F_n(h)=f(nh).
$$
Observing that $F_{\overline{n}}(\overline{h})=f(\overline{x})$
and taking into account the biLipschitz estimate \eqref{biLip}
we have
$$
F_{\overline{n}}(h)\neq f(\overline{x})\quad\mbox{for every}\quad
h\in\partial D^H_{\overline{h},s}\,.
$$
Injectivity of
$F_{\overline{n}}$ implies that 
$$
\deg\big(F_{\overline{n}},B^H_{\overline{h},s},f(\overline{x})\big)
\in\{-1,1\},
$$
see for instance Theorem 3.3.3 of \cite{LLoyd}.
By continuity of 
$$
n\lra\max_{h\in\partial D^H_{\overline{h},s}}
\rho\big(f(nh),f(\overline{n}h)\big)\,,
$$
up to choosing a smaller $r>0$, we can assume that
\begin{eqnarray}
\max_{h\in\partial D^H_{\overline{h},s}}
\rho\big(f(nh),f(\overline{n}h)\big)<\frac{s}{2\beta}
\end{eqnarray}
for every $n\in D^N_{\overline{n},r}$.
As a consequence, applying \eqref{biLip},
it follows that
$$
\rho\big(F_n(h),f(\overline{x})\big)\geq
\rho\big(F_{\overline{n}}(h),f(\overline{x})\big)
-\rho\big(F_{\overline{n}}(h),F_n(h)\big)>\frac{s}{2\beta}
$$
for every $h\in\partial D^H_{\overline{h},s}$
and every $n\in D^N_{\overline{n},r}$.
For an arbitrary $n\in D^N_{\overline{n},r}\sm\{\overline{n}\}$
one can consider the continuous curve
$\gamma:[0,1]\lra D^N_{\overline{n},r}$,
defined by $\gamma(t)=\overline{n}\delta_t\big((\overline{n})^{-1}n\big)$.
Notice that $\gamma$ has image in $D^N_{\overline{n},r}$,
since $N$ is a homogeneous subgroup of $\G$.
By previous estimates, the mapping $\Phi:[0,1]\times D_{\overline{h},s}\lra\M$ defined by
$$
\Phi(t,h)=f\big(\gamma(t)h\big)
$$
is a homotopy between $F_{\overline{n}}$ and $F_n$
such that $\Phi(t,h)\neq f(\overline{x})$
for every $t\in[0,1]$ and every $h\in\partial D^H_{\overline{h},s}$.
Thus, homotopy invariance of topological degree, \cite{JTSch}, implies
$$
\deg\big(F_n,B^H_{\overline{h},s},f(\overline{x})\big)=
\deg\big(F_{\overline{n}},B^H_{\overline{h},s},f(\overline{x})\big)\neq0,
$$
hence there exists at least one element $h'\in D^H_{\overline{h},s}$,
depending on $n$, such that $F_n(h')=f(\overline{x})$.
Injectivity of $F_n$ gives uniqueness of $h'$, hence
there exists $\ph: D^N_{\overline{n},r}\lra D^H_{\overline{h},s}$,
uniquely defined, such that
$$
F_n\big(\ph(n)\big)=f\big(n\,\ph(n)\big)=f(\overline{x}).
$$
This concludes the proof of the first step.

{\bf Step 2}: {\em Regularity.} 

We keep the same notation of the previous step.
By definition of $\mu>0$ in \eqref{mumu}
the restriction $L(x)=\big(Df\big(x)\big)_{|H}:H\lra\M$ is invertible
for every $x\in D_{\overline{x},R_0}$ and setting $T(x)=L(x)^{-1}$,
we have
\begin{eqnarray}\label{unifinvest}
d\big(T(x)(m)\big)\leq \mu^{-1}\;\rho(m)
\end{eqnarray}
for every $m\in\M$. Due to \eqref{applunifestim},
we can choose a possibly smaller $R>0$, hence a smaller
$R_0=R/(2+2cN)$ such that
\begin{eqnarray*}
\max_{x,y\in D_{\overline{x},R_0}}
	\rho\lls f(x)^{-1}f(y),Df(x)(x^{-1}y)\rls\leq\frac{\mu\;d(x,y)}{2}\,.
\end{eqnarray*}
It follows that the remainder $E(n,n')$
in the expression 
\begin{eqnarray*}
f\big(n\ph(n)\big)^{-1}f(n\ph(n'))
=L\big(n\ph(n)\big)\big(\ph(n)^{-1}\ph(n')\big)E(n,n')
\end{eqnarray*}
satisfies the uniform estimate
\begin{eqnarray}\label{appl1unifestim}
\rho\big(E(n,n')\big)\leq\frac{\mu}{2}\;d\big(\ph(n),\ph(n')\big)
\end{eqnarray}
for every $n,n'\in D^N_{\overline{n},r}$.
The implicit mapping $\ph$ satisfies 
$$
f\big(n\ph(n)\big)^{-1}f(n\ph(n')\big)
=\left(f\big(n\ph(n')\big)^{-1}f\big(n'\ph(n')\big)\right)^{-1},
$$
therefore one easily gets 
$$
\ph(n)^{-1}\ph(n')=T\big(n\ph(n)\big)\left(f\big(n\ph(n')\big)^{-1}f\big(n'\ph(n')\big)\right)^{-1}\,T\big(n\ph(n)\big)\big(E(n,n')\big)^{-1}.
$$
As a result, in view of \eqref{unifinvest} and \eqref{appl1unifestim},
we obtain
\begin{eqnarray*}
d\big(\ph(n),\ph(n')\big)\leq \frac{2}{\mu}\;
\rho\left(f\big(n\ph(n')\big),f\big(n'\ph(n')\big)\right).
\end{eqnarray*}
By \eqref{frabove}, up to choosing a possibly smaller $R>0$,
we can suppose that 
\[
\rho\big(f(x),f(y)\big)\leq 2\,\overline{\mu}\,\,d(x,y)
\]
for every $x,y\in D_{\overline{x},R_0}$, hence
\begin{eqnarray}
d\big(\ph(n),\ph(n')\big)\leq \frac{4\,\overline{\mu}}{\mu}\;
d\big(\ph(n')^{-1}n^{-1}n'\ph(n')\big).
\end{eqnarray}
Finally, formulae \eqref{conjest1} and \eqref{leftinveucl} lead us
to the conclusion. $\Box$
\begin{Rem}{\rm
Under hypotheses of Theorem~\ref{implth}, one immediately gets
\begin{eqnarray}\label{holdergraph}
d\big(n\ph(n),n'\ph(n')\big)\leq(1+\kappa)\;
d\big(\ph(n')^{-1}n^{-1}n'\ph(n')\big)\,,
\end{eqnarray}
where $\kappa$ is as in the statement of the above mentioned theorem.
}\end{Rem}
%
%
%
%
%
%
%
%
%
%
%
%
%
%
\vskip.2truecm
\noindent
{\sc Proof of Theorem~\ref{Pembed}.}
By Proposition~\ref{GlinHmono} we have a projection $p:\M\lra H$
that is an h-epimorphism and $p_{|H}=\Id_H$. We also know that
the kernel $N$ of $p$ is complementary to $H$.
As a first consequence, the mapping $p\,\circ f:\Omega\lra H$ is 
continuously P-differentiable and
$D(p\,\circ f)(\overline{x}):\G\lra H$ is invertible.
Due to Theorem~\ref{inverse}, there exist open neighbourhoods
$V\Subset\Omega$ and $W\Subset H$ of $\overline{x}$ and of $p\big(f(\overline{x})\big)$,
respectively, such that the restriction
\[
(p\circ f)_{|\overline{V}}:\overline{V}\lra\overline{W}
\]
is invertible. We denote by $\psi:\overline{W}\lra\overline{V}$ its inverse function.
By Remark~\ref{uniqued}, there exists a unique function
$f_N:\Omega\lra N$ such that $f(x)=p\big(f(x)\big)f_N(x)$
for every $x\in\Omega$, hence we consider $g:\overline{W}\lra\M$
defined by
\[
g(h)=f\big(\psi(h)\big)=h\,f_N\big(\psi(h)\big)
\]
for every $h\in\overline{W}$.
Uniqueness of factorization implies that
\[
g(W')=g(W')\cap W'N,
\]
whenever $W'\subset\overline{W}$. If $W'$ is open, then so is
$W'N$, since the mapping \eqref{phidecomp} is open.
As a consequence, $g$ is an open mapping on its image,
hence so is the restriction $f_{|V}$.
This implies that $f_{|V}$ is a topological embedding.
If we set
\[
\ph=f_N\circ\psi\quad\mbox{and}\quad\Psi=I\circ J^{-1}\circ\psi\circ p_{|f(V)}
\]
then $\ph(h)=h^{-1}g(h)$ and \eqref{imagerep} holds, where
$I$ and $J$ are defined in the statement of Theorem~\ref{Pembed}.
We notice that $g$ is continuously P-differentiable on $\overline{W}$,
then in particular it is Lipschitz on $\overline{W}$ up to a suitable choice of $W$,
for instance a closed ball.
Our next computations are performed in the Lie algebra $\cG$
of $\G$. We set $\ph=\exp\circ F$ and $g=\exp\circ G$, therefore
Lemma~\ref{estimi} gives a constant $C>0$,
depending on both $\overline{W}$ and the Lipschitz constant of $g$ such that
\begin{eqnarray}\label{lipesth}
\|G(h)-G(h')\|\leq C\; d(h,h')\quad\mbox{and}\quad
\|\xi-\xi'\|\leq C\;d(h,h'),
\end{eqnarray}
where $h=\exp\xi$ and $h'=\exp\xi'$, with $h,h'\in\overline{W}.$
We have
\begin{eqnarray*}
&&F(h)-F(h')=(-\xi)\gope G(h)-\big((-\xi')\gope G(h')\big)\\
&&=\xi'-\xi+G(h)-G(h')
+\sum_{n=2}^\iota \left[c_n\big(-\xi,G(h)\big)-c_n\big(-\xi',G(h')\big)\right],
\end{eqnarray*}
due to the formula \eqref{absBCH}. From \eqref{bchdiffprod}, it follows that
\begin{equation}
\|F(h)-F(h')\|\leq\|\xi'-\xi\|+\|G(h)-G(h')\|
+\sum_{n=2}^\iota \gamma_n\,\nu^{n-1}
\big(\|\xi-\xi'\|+\|G(h)-G(h')\|\big),\nonumber
\end{equation}
where we have set
\[
\nu=2\,\Big(\max_{\eta\in\overline{W}}\|G(\eta)\|+\sup_{\eta\in\overline{W}}\|\eta\|\Big).
\]
By virtue of \eqref{lipesth}, we get a constant $C_1>0$, depending
on $\overline{W}$, $\nu$ and the Lipschitz constant of $g$ such that
\[
\|F(h)-F(h')\|\leq C_1\,d(h,h').
\]
This concludes the proof. $\Box$
%
%
%
%
%
%
%
%
%
%
%
%
%
%
%
%
%
\section{Intrinsic graphs and $(\G,\M)$-regular sets}\label{GMregsets}

The notion of complementary subgroup allows us to set properly
the well established notion of intrinsic graph in those graded groups
that admit a factorization by two complementary subgroups.
\begin{Def}\label{intrgraph}{\rm
Let $P$ and $H$ be complementary subgroups of $\G$
and let $S\subset\G$. We say that $S$ is an
{\em intrinsic graph} with respect to $(P,H)$
if there exists a subset $A\subset P$ and a mapping
$\ph:A\lra H$ such that
\[
S=\{p\,\ph(p)\mid p\in A\}.
\]
}\end{Def}
\begin{Rem}[Translations of intrinsic graphs]\label{trintrgr}{\rm
As already observed in \cite{FSSC6}, intrinsic graphs
are preserved under left translations.
If $P$ and $H$ are complementary subgroups of $\G$
and $S=\{p\,\ph(p)\mid p\in A\}$ is an intrinsic graph
with respect to $(P,H)$, then the left translated
$xp\,\ph(p)$, for some $x\in\G$,
is of the form $\tilde p\,\psi(\tilde p)$ where
$\tilde p\in\tilde A\subset P$ and $\psi:\tilde A\lra H$.
In fact, by unique decomposition of elements in $\G$
with respect to $P$ and $H$, see Remark~\ref{uniqued},
there exist functions $\psi_1:A\lra P$
and $\psi_2:A\lra H$ such that
\[
xp\,\ph(p)=\psi_1(p)\psi_2(p)\,.
\]
To conclude, it suffices to show that $\psi_1$ is injective.
If $\psi_1(p)=\psi_1(p')$, then
\[
xp\,\ph(p)\psi_2(p)^{-1}=xp'\ph(p')\psi_2(p')^{-1},
\]
where $\ph(p)\psi_2(p)^{-1},\ph(p')\psi_2(p')^{-1}\in H$.
Again, uniqueness yields $p=p'$. 
}\end{Rem}
\begin{Def}[($\G,\M$)-regular set of $\G$]\label{gmsing}{\rm
Let $\G$ be a stratified group and let $\M$ be a graded group
such that $\M$ is an h-quotient of $\G$.
We say that a subset $S\subset\G$ is {\em $(\G,\M)$-regular}
if for every point $\overline{x}\in S$, there exists an open neighbourhood
$U$ of $\overline{x}$ and a continuously P-differentiable mapping
$f:U\lra\M$ such that $S\cap U=f^{-1}(e)$ and 
$Df(x):\G\lra\M$ is an h-epimorphism for every $x\in U$.
}\end{Def}
\begin{Def}[($\G,\M$)-regular set of $\M$]\label{gmsinm}{\rm
Let $\G$ be a stratified group and let $\M$ be a graded group
such that $\G$ h-embeds into $\M$.
A subset $S\subset\M$ is {\em $(\G,\M)$-regular}
if for every point $\overline{x}\in S$, there exist open neighbourhoods
$U\subset\M$ of $\overline{x}$ and $V\subset\G$ of $e\in\G$,
along with a continuously P-differentiable topological
embedding $f:V\lra\M$, such that $S\cap U=f(V)$ and 
$Df(y):\G\lra\M$ is an h-monomorphism for every $y\in V$.
}\end{Def}
\begin{Def}[($\G,\M$)-regular set]{\rm
When a subset $S$ is either $(\G,\M)$-regular in $\G$ or
$(\G,\M)$-regular in $\M$, we simply say that it is $(\G,\M)$-regular
without further specification, or we can say that it is an
{\em intrinsically regular set}.
}\end{Def}
\subsection{Tangent cone to $(\G,\M)$-regular sets}\label{intrblw}
The notion of ``tangent cone" given in 3.1.21 of \cite{Fed}
can be easily extended to graded groups.
\begin{Def}\label{tancone}{\rm
Let $\G$ be a graded group and let $S\subset\G$ with $a\in\G$.
The {\em homogeneous tangent cone} of $S$ at $a$ is the homogeneous subset
\begin{eqnarray}
&&\mbox{Tan}(S,a)=\Big\{v\in\G\mid v=\lim_{k\ra\infty}\delta_{r_k}(a^{-1}s_k),\,
\mbox{for some sequences}\; (r_k)\subset\R^+,\\
&&(s_k)\subset S,\;\mbox{where}\;s_k\ra a\Big\}. \nonumber
\end{eqnarray}
}\end{Def}
\begin{Rem}{\rm
If $a\in\overline{S}$, then $\mbox{Tan}(S,a)\neq\emptyset$.}
\end{Rem}
{\sc Proof of Theorem~\ref{blwlevelimage}.}
Replacing $f$ with $f\circ l_x$, where $l_x(y)=xy$, 
it is not restrictive to assume that $\overline{x}=e$,
since
\[
(\overline{x})^{-1}S=(f\circ l_{\overline{x}})^{-1}(f\circ l_{\overline{x}})(e)
\]
and $Df(\overline{x})=D(f\circ l_{\overline{x}})(e)$.
Now we use the notation of Theorem~\ref{implth} applying it in the
case $\overline{x}=e$.
From P-differentiability of $f$ at $e$ and representation \eqref{impleq}, we have
\[
f(e)=f(n\ph(n))=f(e)
Df(e)\big(\ph(n)\big)o\big(n\ph(n)\big).
\]
Setting $L=Df(e)$, we obtain
\begin{eqnarray}
d\big(\ph(n)\big)=
d\left(L_{|H}^{-1}\big(o(n\ph(n))\big)\right)\leq \ep \big(d(n)+d(\ph(n))\big)
\end{eqnarray}
for $n$ suitable close to $e$ and $\ep>0$ small, arbitrarily fixed.
From this, we easily conclude that $\ph$ is P-differentiable at $e_N$ and
$D\ph(e_N)$ is the null h-homomorphism.
Now, we arbitrarily fix $R>0$ and consider $\lambda>0$ sufficiently small,
such that
\[
D_R\cap\delta_{1/\lambda}S=D_R\cap
\{\delta_{1/\lambda}n\delta_{1/\lambda}\ph(n)\mid n\in D_r^N\}.
\]
Let $(\lambda_k)$ be an infinitesimal sequence of positive numbers.
Using Proposition 4.5.5 of \cite{AmbTil}, we have to prove 
the following conditions:
\begin{itemize}
\item[(i)] if $x=\lim_{k\to\infty} x_k$ for some sequence $(x_k)$
such that $x_k\in D_R\cap\delta_{1/\lambda_k}S$, then we have
$x\in D_R\cap N$;
\item[(ii)] if $x\in D_R\cap N$, then there exists a sequence
$(x_k)$ such that $x_k\in D_R\cap\delta_{1/\lambda_k}S$
and $x_k\to x$.
\end{itemize}
To prove (i), we write $x_k=\big(\delta_{1/\lambda_k}n_k\big)\big(\delta_{1/\lambda_k}\ph(n_k)\big)$
and use Proposition~\ref{homdec} with respect to our complementary subgroups
$N$ and $H$. Through the global diffeomorphism $\phi$ in \eqref{phidecomp}, 
convergence of the sequence $(x_k)$ implies convergence of both $(\delta_{1/\lambda_k}n_k)$ and $\big(\delta_{1/\lambda_k}\ph(n_k)\big)$.
We notice that the limit $n_0$ of $(\delta_{1/\lambda_k}n_k)$
belongs to $N$. P-differentibility of $\ph$ at $e_N$ gives
\[
\lim_{k\ra\infty}\delta_{1/\lambda_k}\ph(\delta_{\lambda_k}n'_k\big)
=D\ph(e_N)(n_0)=e_H
\]
where we have set $\delta_{1/\lambda_k}n_k=n'_k$. Then $x_k\ra x=n_0\in D_R\cap N$
and our claim is achieved. Now we choose $x\in N\cap D_R$ and consider
a sequence $(x_l)$ contained in $N\cap B_R$ converging to $x$.
We observe that for every fixed $l$, there exists a sufficiently large
$k_l$ such that $\delta_{\lambda_{k_l}}x_l\in D^N_r$ and
\[
\delta_{1/\lambda_{k_l}}
\big(\delta_{\lambda_{k_l}} x_l\ph(\delta_{\lambda_{k_l}}x_l)\big)
\in B_R\cap\delta_{1/\lambda_k}S,
\]
since P-differentiability of $\ph$ at $e_N$ implies that the 
sequence $\delta_{1/\lambda_{k}}\big(\delta_{\lambda_k}x_l
\ph(\delta_{\lambda_k}x_l)\big)$ converges to $x_l$ as $k\ra\infty$.
The fact that $D\ph(e_N)$ is the null h-homomorphism also implies that
\[
\delta_{1/\lambda_{k_l}}
\big(\delta_{\lambda_{k_l}} x_l\ph(\delta_{\lambda_{k_l}}x_l)\big)\lra
xD\ph(e_N)(x)=x\quad\mbox{as}\quad l\ra\infty.
\]
Finally, we notice that conditions (i) and (ii) exactly prove that
the two inclusions between $N$ and $\mbox{Tan}(S,e)$.
This concludes the proof of the first part of Theorem~\ref{blwlevelimage}.

Now, we keep the same notation used in the proof of Theorem~\ref{Pembed}
and observe that it is not restrictive assuming that $f(\overline{x})=e_\M$.
In fact, in the general case it suffices to replace $f$ with
the left translated $x\ra f(\overline{x})^{-1}f(x)$,
that holds the same assumptions. 
We know that the mapping $p\circ f_{|V}$
is invertible with P-differentiable inverse $\psi:W\lra V$ and
$g(h)=f\circ \psi(h)=h\ph(h)$ for every $h\in W$.
By Proposition~\ref{chain}, we have
that $Dg(e_H):H\hookrightarrow\M$ identically embeds $H$ into $\M$, then
\[
g(h)=h\,o(h)=h\,\ph(h)
\]
so that $\ph(h)=o(h)$. This implies that there exists
$D\ph(e_H)$ and it is constantly equal to the unit element $e_N\in N$.
To show the remaining claims, one can argue exactly as in the first part
of this proof, using the fact that $\ph$ is P-differentiable at $e_H$
with null P-differential at this point and replacing $N$ with $H$. $\Box$
\begin{Cor}\label{GMGTan}
Let $S$ be a $(\G,\M)$-regular set of $\G$. Then for every
$x\in S$, we have
\[
\mcH\mbox{\rm -dim}\big(\mbox{\rm Tan}(S,x)\big)
=\mcH\mbox{\rm -dim}\G-\mcH\mbox{\rm -dim}\M\,.
\]
\end{Cor}
{\sc Proof.}
By Theorem~\ref{blwlevelimage}, $\mbox{\rm Tan}(S,x)$ is h-isomorphic to
$N=\ker Df(x)$, where $f$ is the continuously P-differentiable mapping
defining $S$ around $x$ as a level set. Due to Proposition~\ref{GlinH},
the factorizing property of h-epimorphisms yields a complementary subgroup
$H$ that is h-isomorphic to $\M$. Due to formula \eqref{QKP}, our claim
follows. $\Box$
\begin{Exa}\label{nnhisotan}{\rm
The fact that homogeneous tangent cones to a $(\G,\M)$-regular set of $\G$
have a fixed Hausdorff dimension does not mean that they are all
algebraically h-isomorphic. Let us consider the mapping
\[
f:\H^2\lra\R^2,\qquad
(x_1,x_2,x_3,x_4,x_5)\lra\Big(\sqrt{x_2^2+x_3^2},x_4\Big)
\]
where $p=\exp\big(\sum_{j=1}^5x_jX_j)\in\H^2$ and 
$[X_1,X_2]=[X_3,X_4]=X_5$ are the nontrivial bracket relations of $\H^2$.
The Pansu differential $Df(x)$ is represented by the matrix
\[
\left(\begin{array}{ccccc}
0 & \frac{x_2}{\sqrt{x_2^2+x_3^2}} & \frac{x_3}{\sqrt{x_2^2+x_3^2}} & 0   & 0 \\
0 & 0       & 0   & x_4 & 0 
\end{array}\right)\,.
\]
Then, defining the connected open set
\[
\Omega=\left\{(x_1,x_2,x_3,x_4,x_5)\in\H^2\mid x_4>0,\;x_2^2+x_3^2>0 \right\}
\]
we notice that $S=\Omega\cap f^{-1}\big((1,1)\big)$ is $(\H^2,\R^2)$-regular in $\H^2$.
Now we fix $\xi=(0,1,0,1,0)$ and $\eta=(0,0,1,1,0)$, observing that
$\xi,\eta\in S$. Thus, it is easy to check that
$\mbox{Tan}(S,\xi)$ is commutative and $\mbox{Tan}(S,\eta)$ is
h-isomorphic to $\H^1$.
}\end{Exa}
%
%
%
%
%
%
%
%
%
%
%
%
%
%
%
%
\section{Factorizing groups}\label{sectfactor}
In this section we investigate the algebraic conditions
under which either surjective or injective h-homomorphisms
are h-monomorphisms or h-epimorphisms, respectively.
Let $\G$ and $\M$ be arbitrary graded groups
with algebras $\cG$ and $\cM$.
\begin{Def}[h-quotients and h-embeddings]\label{hquotembed}{\rm
We say that $\M$ is an {\em h-quotient} of $\G$
if there exists a normal homogeneous subgroup $N\subset\G$
such that $\G/N$ is h-isomorphic to $\M$.
Analogously, $\G$ {\em h-embeds} into $\M$ if and only
if there exists a homogeneous subgroup $H$ of $\M$ which is
h-isomorphic to $\G$.
}\end{Def}
\begin{Pro}\label{normhsu}
There exists a normal subgroup $N$ of $\G$ such that $\G/N$
is h-isomorphic to $\M$ if and only if there exists a surjective
h-homomorphism $L:\G\lra\M$.
\end{Pro}
{\sc Proof.}
If the first condition holds, then we denote by $T:\G/N\lra\M$
the canonical h-isomorphism and notice that $L=T\circ \pi:\G\lra\M$ 
is a surjective h-homomrphism, where
$\pi:\G\lra\G/N$ is the canonical projection, that is also an h-homomorphism, according to Proposition~\ref{quotgrad}.
The converse is trivial. $\Box$
\begin{Rem}\label{hhomosub}{\rm
The previous proposition shows that $\M$ is an {\em h-quotient}
of $\G$ if and only if there exists a surjective h-homomorphism $L:\G\lra\M$.
In view of Proposition~\ref{quotgrad}, we also notice that $\M$ is stratified
when so is $\G$.
On the other side, it is easy to check that $\G$ h-embeds into $\M$
if and only if there exists an injective h-homomorphism $T:\G\lra\M$.
}\end{Rem}
\begin{Lem}\label{hepim}
If $m$ is the dimension of the first layer of $\cG$,
then $\R^k$ is an h-quotient of $\G$ if and only if $k\leq m$.
\end{Lem}
{\sc Proof.}
Let $V_1\oplus V_2\oplus\cdots\oplus V_\iota$ be the direct decomposition
of $\cG$ into its layers. Let $\cu$ be an $(m-k)$-dimensional
subspace of $V_1$. Then $\cN=\cu\oplus V_2\oplus\cdots\oplus V_\iota$
is an ideal and clearly $\cG/\cN$ is h-isomorphic to $\R^k$.
By Proposition~\ref{algquot}, $\G/N$ is also h-isomorphic to
$\R^k$, where $N=\exp\cN$. Conversely, assume that
$\R^k$ is an h-quotient of some $\G$. From Proposition~\ref{quotgrad},
it follows that there exists $N$ of ideal $\cN=\cN_1\oplus\cdots\oplus\cN_\iota$ such that $\R^k$ is linearly isomorphic to $V_1/\cN_1$,
then $k\leq\dim V_1=m$. $\Box$
\begin{Exa}{\rm
The Heisenberg group $\H^k$ is not an h-quotient of $\H^n$,
whenever $n>k$. By contradiction, assume that there exists
a normal subgroup $N=\exp\cn$ such that $\H^n/N$ is h-isomorphic
to $\H^k$. In terms of algebras, we have a $2(n-k)$-dimensional
ideal $\cn$ of $\ch^n$ such that $\ch^n/\cn$ is h-isomorphic to $\ch^k$. 
This implies that there exists $\cn+Z'\in\ch^n/\cn$ that corresponds
to a nonvanishing element of the second layer of $\ch^k$.
It follows that $\cn\subset\cv$ and $[\cn,\cn]=\{0\}$,
since $\cv$ is the first layer of $\ch^n$.
If we pick an element $X\in\cv$ such that $[X,\cn]\neq\{0\}$,
then we meet a contradiction. In fact, $X\notin\cn$,
since $\cn$ is commutative and this conflicts
with the fact that $\cn$ is an ideal.
}\end{Exa}
\begin{Rem}{\rm
In view of Remark~\ref{hhomosub}, the previous example
shows the nonexistence of surjective h-homomorphisms
from $\H^n$ to $\H^k$ whenever $n>k$, as it was proved
in Theorem~2.8 of \cite{Mag1}, by explicit representation
of h-homomorphisms from $\H^n$ to $\H^k$.
}\end{Rem}
\begin{Def}[Factorizing group as a quotient]\label{factorquot}{\rm
We say that $\M$ {\em factorizes $\G$ as a quotient} if it is an h-quotient
of $\G$ and every normal subgroup $N$ of $\G$ such that $\G/N$
is h-isomorphic to $\M$ has a complementary subgroup $H$.
}\end{Def}
\begin{Def}[Factorizing group as a subgroup]\label{factorsubg}{\rm
We say that $\G$ {\em factorizes $\M$ as a subgroup} if it h-embeds into
$\M$ and every subgroup $H$ of $\M$ which is h-isomorphic to $\G$
has a complementary normal subgroup $N$.
}\end{Def}
As a corollary of both Proposition~\ref{GlinH} and
Proposition~\ref{GlinHmono}, we have the following 
\begin{Pro}\label{factorpro}
$\M$ factorizes $\G$ as a quotient if and only if every surjective 
h-homomorphism is an h-epimorphism and $\G$ factorizes 
$\M$ as a subgroup if and only if every injective h-homomorphism is an
h-monomorphism.
\end{Pro}
\begin{Pro}\label{rksub}
If $\R^k$ h-embeds into $\G$, then it factorizes $\G$ as a subgroup.
\end{Pro}
{\sc Proof.}
By hypothesis, the class of $k$-dimensional subgroups
$H=\exp \ch$ of $\G$ which are h-isomorphic to $\R^k$ is nonempty
and it correspond to $k$-dimensional commutative subgroups contained
in $\exp V_1$. Let $\cN_1$ be a subspace of $V_1$ such that
$\ch\oplus\cN_1=V_1$. Then $\cN=\cN_1\oplus V_2\oplus\cdots\oplus V_\iota$
is an ideal and $N=\exp \cN$ is complementary to $H$. $\Box$
\begin{Rem}\label{rksubrm}{\rm
If $\R^k$ h-embeds into $\M$, then every injective h-homomorphism
from $\R^k$ to $\M$ is an h-monomorphism. 
This immdiately follows joining
Proposition~\ref{factorpro} and Proposition~\ref{rksub}. 
}\end{Rem}
\begin{Def}[Factorizing group]\label{factor}{\rm
We say that $\G$ {\em factorizes $\M$} if it factorizes $\M$
both as a subgroup and as a quotient.
}\end{Def}
\begin{Rem}
{\rm It might be interesting from an algebraic viewpoint
to investigate whether factorizing groups as quotient are also
factorizing as a subgroup and viceversa.
}\end{Rem}
\begin{Pro}\label{rfactg}
$\R$ factorizes every stratified group.
\end{Pro}
{\sc Proof.}
By Lemma~\ref{hepim}, $\R$ is an h-quotient of an arbitrary stratified group $\G$. Let $N=\exp\big(\cN_1\oplus\cdots\oplus\cN_\iota\big)$
be a homogeneous normal subgroup such that $\G/N$ is h-isomorphic to $\R$.
By Proposition~\ref{quotgrad}, we must have that
$V_1/\cN_1$ is one-dimensional.
The Lie algebra $\cG$ of $\G$ is the direct sum
$V_1\oplus\cdots\oplus V_\iota$.
We choose $X\in V_1\sm\cN_1$ and define \[\ch=\{tX\mid t\in\R\}.\]
It follows that $\cn$ and $\ch$ are homogeneous subalgebras
such that $\cn\oplus\ch=\cG$. By virtue of Proposition~\ref{homdec},
$\R$ factorizes $\G$ as a quotient.
Finally, by Proposition~\ref{rksub} $\R$ factorizes $\G$
as a subgroup, hence our proof is complete. $\Box$ 
\begin{Exa}\label{cntrex}{\rm
One can easily find $\G$ and $\M$ such that there exist
h-epimorphisms $L:\G\lra\M$, although $\M$ does not factorize $\G$
as a quotient. In fact,
let us consider the 2-step stratified group $\G$
with Lie algebra $\cG=V_1\oplus V_2$, where
\[
V_1=\span\{X_1,X_2,X_3,X_4\},\qquad V_2=\span\{Z_{23},Z_{24},Z_{34}\}
\]
and the only nontrivial bracket relations are the following
$$
[X_2,X_3]=Z_{23},\quad[X_2,X_4]=Z_{24},
\quad\mbox{and}\quad [X_3,X_4]=Z_{34}\,.
$$ 
Let $L_1:\G\lra\R^2$ be the h-epimorphism defined by
\[
L_1\Big(\exp
\big(\sum_{i=1}^4x_iX_i+z_{23}Z_{23}+z_{24}Z_{24}+z_{34}Z_{34}\big)\Big)
=(x_1,x_2)
\]
and let $\ker L_1=\exp\cn=\exp(\cn_1\oplus\cn_2)=N$.
We have $\cn_2=V_2$ and $\cn_1=\span\{X_3,X_4\}$,
then $H=\exp\{X_1,X_2\}$ is complementary to $N$,
as it easily follows from Proposition~\ref{homdec}.
However, if we consider the h-epimorphism
$L_2:\G\lra\R^2$ defined by
\[
L_2\Big(\exp
\big(\sum_{i=1}^4x_iX_i+z_{23}Z_{23}+z_{24}Z_{24}+z_{34}Z_{34}\big)\Big)
=(x_3,x_4),
\]
then $\ker L_2=\exp\cn$, with
$\cn=\span\{X_1,X_2\}\oplus V_2$ and we will check
that $N=\exp\cn$ does not admit any complementary subgroup.
By contradiction, if $H$ is complementary to $N$,
then Proposition~\ref{GlinH} shows that 
the restriction $T:H\lra\R^2$ of $L_2$ is a group isomorphism.
In particular, $H$ is commutative.
In addition, Proposition~\ref{homdec} shows that 
$\ch\oplus\cn=\cG$, where $\ch$ denotes the Lie algebra of $H$.
Then $\ch$ is a 2-dimensional commutative subalgebra of $\cG$.
We consider a basis $(v,w)$ of $\ch$, given by
$$
v=\sum_{j=1}^4\alpha_{j}X_j+Z,\quad\mbox{and}\quad
w=\sum_{j=1}^4\beta_jX_j+T,
$$
where $T,Z\in V_2$. The decomposition $\ch\oplus\cn=\cG$
implies that $(v,w,X_1,X_2,Z_{23},Z_{24},Z_{34})$ is a basis of
$\cG$, hence we must have $\alpha_3\beta_4-\beta_3\alpha_4\neq0$.
As a consequence,
$$
[v,w]=(\alpha_2\beta_3-\alpha_3\beta_2)Z_{23}+
(\alpha_2\beta_4-\alpha_4\beta_2)Z_{24}+
(\alpha_3\beta_4-\alpha_4\beta_3)Z_{34}\neq0\,.
$$
This conflicts with the fact that $\ch$ is commutative.
}\end{Exa}

\subsection{Factorizations in some H-type groups}\label{htypeg}
Our aim here is to present some factorizing properties of $\R^k$
in some important H-type groups. These groups were introduced in \cite{Kap}.
\begin{Def}\label{htype}{\rm
Let $\cg$ be a Lie algebra equipped with an inner product $\lan\cdot,\cdot\ran$
and let $\cz$ be a nontrivial subspace of $\cg$ such that $[\cg,\cz]=0$ and
$[\cg,\cg]\subset\cz$. Let $\cv$ be the orthogonal space of $\cz$
and define $J:\cz\lra\mbox{End}(\cv)$ by the formula
\begin{eqnarray}
\lan J_ZX,Y\ran=\lan Z,[X,Y]\ran
\end{eqnarray}
for every $Z\in\cz$ and $X,Y\in\cv$. If $J$ satisfies the condition
$$
|J_ZX|=|Z|\,|X|\,,
$$
then we say that $\cg$ is an {\em H-type algebra}.
}\end{Def}
From the previous definition, it follows that
\begin{eqnarray}\label{keyhtype}
[X,J_ZX]=|X|^2\;Z\qquad\mbox{for every}\qquad (X,Z)\in \cv\times\cz\,.
\end{eqnarray}
Then $\cg=\cv\oplus\cz$ is a 2-step stratified algebra of first
layer $\cv$ and second layer $\cz$.
In the sequel, we will utilize another well known formula
\begin{eqnarray}\label{polH}
J_ZJ_W+J_WJ_Z=-2\,\lan Z,W\ran\,\mbox{Id}_\cv.
\end{eqnarray}
\begin{Exa}\label{heishtype}
{\rm Let $(X,Y,Z)$ be an orthonormal basis of the Heisenberg algebra
$\ch^1$, where $[X,Y]=Z$. Then setting $\cz=\span\{Z\}$, $\cv=\span\{X,Y\}$
and
\[
J_ZX=Y,\qquad J_ZY=-X
\]
extended by linearity, it follows that $J$ makes $\ch^1$ an H-type group.
The higher dimensional Heisenberg algebras $\ch^n$
can be seen as direct product of the irreducible Heisenberg algebras
isomorphic to $\ch^1$ as follows
\[
\ch^n=\cv_1\oplus\cdots\oplus\cv_n\oplus\cz,
\]
where $(X_i,Y_i)$ is an orthonormal basis of $\cv_i$ and
$J_ZX_i=Y_i$, $J_ZY_i=-X_i$.
Here we notice that setting $\omega(U_1,U_2)=\lan J_Z U_1,U_2\ran$,
we define a symplectic form on the $2n$-dimensional first layer
$\cv=\cv_1\oplus\cdots\oplus\cv_n$.
}\end{Exa}
\begin{Exa}\label{h12k}{\rm
It is easy to check that $\R^2$ does not factorize $\ch^1$ as a quotient,
although it is an h-quotient of $\ch^1$.
In fact, the only 1-dimensional ideal $\cn$ of 
$\ch^1$ is the second layer $\cz$ and we have already shown in 
Example~\ref{nnexcmpl} that there do not exist 2-dimensional subalgebras 
complementary to $\cz$.
}\end{Exa}
However, as first observed in \cite{FSSC6}, $\R^k$ factorizes
the Heisenberg group $\H^n$ as a quotient, whenever $k\leq n$.
In view of Proposition~\ref{GlinH}, the next statement
translates Proposition~3.24 of \cite{FSSC6} into our setting,
with a different proof.
\begin{Pro}\label{kfacth}
If $1\leq k\leq n$, then $\R^k$ factorizes $\H^n$.
\end{Pro}
{\sc Proof.}
By Lemma~\ref{hepim}, $\R^k$ is an h-quotient $\H^n$.
Let $N$ be a normal homogeneous subgroup of $\H^n$ such that
$\H^n/N$ is h-isomorphic to $\R^k$. Let $\cn=\cn_1\oplus\cn_2$
be the Lie algebra of $N$, which is also an ideal of $\ch^n$.
It is clear that $\cn_2=\cz$, where $\cv$ and $\cz$ denote
the first and the second layers of $\ch^n$, respectively.
If we prove that there exists a $k$-dimensional commutative subalgebra $\cs$
contained in $\cv$ such that $\cs\oplus\cn_1=\cv$,
then Proposition~\ref{homdec} concludes the proof.

To prove this, we fix a unit vector $Z\in\cz$ and consider the mapping
$J_Z:\cv\lra\cv$ that makes $\ch^n$ an H-type algebra,
see Example~\ref{heishtype}. We set $\cn_1\cap J_Z(\cn_1)=\cw$
and notice that $J_Z(\cw)=\cw$. If $\cw\neq\{0\}$,
then $\cw$ is a symplectic subspace with respect to the
symplectic form $\omega(X,Y)=\lan J_ZX,Y\ran$.
In fact, one can easily construct the following orthogonal basis 
$(e_1,\ldots,e_l,J_Z(e_1),\ldots,J_Z(e_l))$ of $\cw$.
We proceed as follows, choosing a unit vector $e_1\in\cw$ and then
considering the orthogonal $J_Z(e_1)\in\cw$. The procedure continues
selecting a unit vector $e_2\in\cw\cap\span\{e_1,J_Z(e_1)\}^\bot$ and checking
that $J_Z(e_2)\in\cw\cap\span\{e_1,J_Z(e_1)\}^\bot$ and it is
repeated up to reaching the desired basis.
Notice that this basis is orthonormal by construction.
If $p=\dim(\cn_1)>2l$,
then we choose $u_1\in \cn_1\cap\cw^\bot$ and observe that
$J_Z(u_1)\notin\cn_1$, otherwise $u_1,J_Z(u_1)$ would belong to
$\cw$, that is a contradiction. One can iterate this argument up
to reaching $u_1,\ldots,u_{p-2l}$, where for every $j=1,\ldots,p-2l$
we have
\[
J_Z(u_j)\notin\cn_1.
\]
By construction $u_1,\ldots,u_{p-2l},J_Z(u_1),\ldots,J_Z(u_{p-2l})$
are orthogonal vectors and 
\[
(e_1,\ldots,e_l,J_Z(e_1),\ldots,J_Z(e_l),u_1,\ldots,u_{p-2l})
\]
is an orthonormal basis of $\cn_1$.
If $(p-l)<n$, then repeating the previous argument we can
complete 
\begin{eqnarray}\label{symplbasis}
\left(e_1,\ldots,e_l,J_Z(e_1),\ldots,J_Z(e_l),u_1,\ldots,u_{p-2l},
J_Z(u_1),\ldots,J_Z(u_{p-2l})\right)
\end{eqnarray}
to the following orthonormal basis of $\cv$
\begin{eqnarray*}
&&\big(e_1,\ldots,e_l,J_Z(e_1),\ldots,J_Z(e_l),u_1,\ldots,u_{p-2l},
J_Z(u_1),\ldots,J_Z(u_{p-2l}),\\
&&w_1,\ldots,w_r,J_Z(w_1),\ldots,J_Z(w_r)\big)\,.
\end{eqnarray*}
This is precisely a symplectic basis with respect to $\omega(\cdot,\cdot)$, since
\[
\lan Z,[w_s,J_Zw_s]\ran=\lan J_Zw_s,J_Zw_s\ran=1
=\lan Z,[e_i,J_Ze_i]\ran=\lan Z,[u_j,J_Zu_j]\ran
\]
implies that $[w_s,J_Zw_s]=[e_i,J_Ze_i]=[u_j,J_Zu_j]=Z$
and in the same way one can check that these are the only nontrivial brackets.
For instance, we have
\[
\lan Z,[J_Z(u_i),w_j]\ran=\lan J_Z^2(u_i),w_j\ran=-\lan u_i,w_j\ran=0
\]
that gives $[J_Z(u_i),w_j]=0$. 
Finally, the dimensional condition $2l+2(p-2l)+2r=2n$ gives $l-r=p-n\geq0$,
by our hypothesis on $p\geq n$. Thus, the basis
\[
\Big(J_Z(u_1),\ldots,J_Z(u_{p-2l}),
w_1+e_1,\ldots,w_r+e_r,J_Z(w_1)-J_Z(e_1),\ldots,J_Z(w_r)-J_Z(e_r)\Big)
\]
defines a $k$-dimensional subspace $\cs$ that is commutative
and by construction satisfies $\cs\cap\cn_1=\{0\}$.
This shows that $\R^k$ factorizes $\H^n$ as a quotient.
Due to Proposition~\ref{rksub}, $\R^k$ also factorizes $\H^n$
as a subgroup. This concludes the proof. $\Box$
\begin{Rem}{\rm
Joining Proposition~\ref{factorpro} and Proposition~\ref{rfactg},
it follows that every surjective h-homomorphism $L:\G\lra\R$ and every injective h-homomorphism $L:\R\lra\G$ are an h-epimorphism and an h-monomoprhism, respectively.
Joining Proposition~\ref{factorpro} and Proposition~\ref{kfacth},
for every $1\leq k\leq n$, it follows that every surjective h-homomorphism $L:\H^n\lra \R^k$ and every injective h-homomorphism $L:\R^k\lra\H^n$ are an h-epimorphism and an h-monomoprhism, respectively.
On the other hand, Example~\ref{h12k} shows that 
$\R^k$ is not an h-quotient of $\H^n$ whenever $k>n$.
}\end{Rem}
Another example of H-type group is the complexified Heisenberg group $\H_2^1$, 
where the center $\cz$ of its Lie algebra $\ch_2^1$ is 2-dimensional and the 
first layer $\cv$ has dimension four. More information on this group
can be found in \cite{ReiRic}.
\begin{Lem}\label{baslemcompheis}
Let $\ch_2^1=\cv\oplus\cz$ be the 6-dimensional real Lie algebra of the complexified Heisenberg group $\H_2^1$ and 
let $X\in\cv$ with $|X|=1$. Let $(Z_1,Z_2)$ be an orthonormal
basis of $\cz$. Then the following vectors
\begin{eqnarray}\label{R_i}
R_0=X,\quad R_1=J_{Z_1}X,\quad R_2=J_{Z_2}X,\quad R_3=J_{Z_1}J_{Z_2}X,
\end{eqnarray}
form an orthonormal basis of $\cv$ and the only nontrivial bracket
relations are given by
\begin{eqnarray}
[R_0,R_1]=[R_2,R_3]=Z_1, \quad [R_0,R_2]=-[R_1,R_3]=Z_2.
\end{eqnarray}
\end{Lem}
{\sc Proof.}
Using just the properties of $J$ and formula
$[Y,J_ZY]=|Y|^2Z$ for every $Y\in\cv$ and $Z\in\cz$
it is easy to check that $R_i$ form an orthonormal basis of $\cv$.
The previous formula also yields
$[R_0,R_1]=Z_1$ and $[R_0,R_2]=Z_2$. In addition, we have
\begin{eqnarray*}
&&\lan Z_1,[R_1,R_2]\ran=\lan J_{Z_1}^2X,J_{Z_2}X\ran=-\lan X,J_{Z_2}X\ran=0, \\
&&\lan Z_2,[R_1,R_2]\ran=\lan J_{Z_2}J_{Z_1}X,J_{Z_2}X\ran
=\lan J_{Z_1}X,X\ran=0, 
\end{eqnarray*}
then $[R_1,R_2]=0$. To prove that $[R_0,R_3]=0$,
$[R_1,R_3]=-Z_2$ and $[R_2,R_3]=Z_1$,
one argues in the same way, using also formula \eqref{polH}. $\Box$
\begin{Pro}\label{h21r2}
$\R^2$ factorizes the complexified Heisenberg group $\H_2^1$ 
as a quotient.
\end{Pro}
{\sc Proof.}
Lemma~\ref{hepim} ensures that $\R^2$ is an h-quotient of $\H_2^1$.
Let $N$ be the kernel of an h-epimorphism $L:\H_2^1\lra\R^2$
and set $\cn=\exp^{-1}(N)=\cn_1\oplus\cn_2$, with
$\cn_1=\cn\cap V_1$ and $\cn_2=\cn\cap V_2$, due to
Proposition~\ref{pfactrs}.
Clearly $\cn_2=\cz$ and $\dim(\cn_1)=2$. 
We have either $\dim([\cn_1,\cn_1])=0$ or $\dim([\cn_1,\cn_1])=1$.
In the first case we choose an orthonormal basis
$(X,Y)$ of $\cn_1$ and represent
$Y$ as a linear combination of $R_i$, according to \eqref{R_i}.
Then the fact that $(X,Y)$ is an orthonormal basis and
$[X,Y]=0$ imply that $Y=J_{Z_1}J_{Z_2}X$ for a
fixed orthonormal basis $(Z_1,Z_2)$ of the center $\cz$.
Now we simply notice that the commutative subalgebra
$\ch=\span\{J_{Z_1}X,J_{Z_2}X\}$ satisfies
$\ch\oplus\cn=\ch_2^1$.
Let us consider the remaining case and take an orthonormal basis
$(X,Y)$ of $\cn_1$. We have $[X,Y]=Z\neq0$.
Let $(T_1,T_2)$ be an orthonormal basis of $\cz$
such that $T_1=Z/|Z|$.
Replacing $(Z_1,Z_2)$ in Lemma~\ref{baslemcompheis} with the orthonormal
basis $(T_1,T_2)$, we get $Y=\alpha_1 J_{T_1}X+\alpha_3J_{T_1}J_{T_2}X$,
with $|\alpha_1|=|Z|>0$.
By direct computation, one can check that the commutative subalgebra
$\ch=\span\{X-\lambda J_{T_2}X,\lambda J_{T_1}X+J_{T_1}J_{T_2}X\}$
satisfies the condition $\ch\oplus\cn=\ch_2^1$ if we fix
$\lambda\neq0$ and $\lambda^{-1}\neq\alpha_3/\alpha_1$.
This concludes the proof. $\Box$
\begin{Pro}\label{nnexcomm}
Let $\ch_2^1=\cv\oplus\cz$ be the complexified Heisenberg algebra.
Then there do not exist commutative subalgebras of $\cv$
with dimension greater than two.
\end{Pro}
{\sc Proof.}
By the general formula \eqref{keyhtype}, for every $X\in\cv\sm\{0\}$
the mapping 
\[
\mbox{ad X}:\cv\lra\cz\,,\qquad Y\lra [X,Y]
\]
is surjective, then its kernel is 2-dimensional. By contradiction,
the existence of a commutative subalgebra of $\cv$ with dimension
greater than two would conflict with the dimension of the kernel. $\Box$
\begin{Cor}
For $k=3,4$ we have that $\R^k$ does not h-embeds into $\H_2^1$
and does not factorize $\H^1_2$ as a quotient.
\end{Cor}
{\sc Proof.}
The first assertion immediately follows from
both Definition~\ref{hquotembed} and Proposition~\ref{nnexcomm}.
Concerning the proof of the second assertion, in view of Proposition~\ref{factorpro} we consider a surjective h-homomorphism
$L:\ch_2^1\lra\R^3$. By contradiction, if $L$ is an h-epimorphism,
then we get a 3-dimensional subalgebra of $\cv$ that is h-isomorphic
to $\R^3$. This conflicts with Proposition~\ref{nnexcomm} and
concludes the proof. $\Box$
\subsection{Factorizations in some free stratified groups}\label{factfree}

We denote by $\cg_{p,\upsilon}$ the free $\upsilon$-step stratified algebra on 
$p$ generators. The corresponding simply connected Lie group
will be denoted by $\G_{p,\upsilon}$.
\begin{Rem}\label{rpupsilon}{\rm
We notice that $\cg_{p,\upsilon}$ is an h-quotient of
$\cg_{r,\upsilon}\oplus\ca$ for every $p\leq r$, 
where $\ca$ is a stratified algebra and $\oplus$ denotes the
direct product of Lie algebras. It suffices to consider the basis
\begin{eqnarray*}
&&\cB_{r,\upsilon}=\big\{[X_{j_1},[X_{j_2},[\cdots,[X_{j_{s-1}},X_{j_s}],],\ldots,] \mid 1\leq s\leq \upsilon,\\
&& (j_1,\ldots,j_s)\in \cA^s_{r,\upsilon}\subset\{1,\ldots,r\}^s\big\}
\end{eqnarray*}
of $\cg_{r,\upsilon}$ and the basis
\begin{eqnarray*}
&&\cB_{p,\upsilon}'=\big\{[X_{j_1}',[X_{j_2}',[\cdots,[X_{j_{s-1}}',
X_{j_s}'],],\ldots,] \mid 1\leq s\leq \upsilon,\\
&& (j_1,\ldots,j_s)\in \cA'^s_{p,\upsilon}\subset\{1,\ldots,p\}^s\big\}
\end{eqnarray*}
of $\cg_{p,\upsilon}$, setting
\begin{eqnarray*}
L\big([X_{j_1},[X_{j_2},[\cdots,[X_{j_{s-1}},X_{j_s}],],\ldots,]\big):=
[X_{j_1}',[X_{j_2}',[\cdots,[X_{j_{s-1}}',X_{j_s}'],],\ldots,]\,,
\end{eqnarray*}
for every $1\leq s\leq \upsilon$ and every
$(j_1,\ldots,j_s)\in \cA'^s_{p,\upsilon}\subset\cA^s_{r,\upsilon}$,
\begin{eqnarray*}
L\big([X_{j_1},[X_{j_2},[\cdots,[X_{j_{s-1}},X_{j_s}],],\ldots,]\big):=0
\end{eqnarray*}
for every $1\leq s\leq \upsilon$ and every 
$(j_1,\ldots,j_s)\in\cA^s_{r,\upsilon}\sm\cA'^s_{p,\upsilon}$ and $L(\ca)=\{0\}$.
It is easy to check that $L:\cg_{r,\upsilon}\oplus\ca\lra\cg_{p,\upsilon}$
is a surjective h-homomorphism.
}\end{Rem}

\begin{Pro}\label{Gpupsilon}
Let $\P$ be a stratified group such that $\G_{p,\upsilon}$
is an h-quotient of $\P$. Then $\G_{p,\upsilon}$ factorizes
$\P$ as a quotient.
\end{Pro}
{\sc Proof.}
We apply Proposition~\ref{factorpro}. 
By hypothesis, we have $L:\cP\lra\cg_{p,\upsilon}$ that is
a surjective h-homomorphism,
where $\cP=W_1\oplus\cdots\oplus W_\iota$ is the stratified algebra of $\P$.
Let $X_1,\ldots,X_p$ be generators of $\cg_{p,\upsilon}$.
Then there exist $U_1,\ldots, U_p\in W_1$ such that $L(U_i)=X_i$.
By hypothesis, we consider the basis  
\[
\cB=\{[X_{j_1},[X_{j_2},[\cdots,[X_{j_{s-1}},X_{j_s}],],\ldots,]\mid
1\leq s\leq \upsilon,\;(j_1,\ldots,j_s)\in \cA^s\subset\{1,\ldots,p\}^s\}
\]
of $\cg_{p,\upsilon}$.
The homomorphism property of $L$ implies that 
\[
L\big([U_{j_1},[U_{j_2},[\cdots,[U_{j_{s-1}},U_{j_s}],],\ldots,]\big)=
[X_{j_1},[X_{j_2},[\cdots,[X_{j_{s-1}},X_{j_s}],],\ldots,]\,,
\]
then the family 
\[
\cB'=\{[U_{j_1},[U_{j_2},[\cdots,[U_{j_{s-1}},U_{j_s}],],\ldots,]\mid
1\leq s\leq \upsilon,\;(j_1,\ldots,j_s)\in \cA'^s\subset\{1,\ldots,p\}^s\}
\]
is a basis of $\ch=\mbox{Lie-}\span\{U_1,\ldots,U_p\}$, then
$L$ maps $\ch$ h-isomorphically onto $\cg_{p,\upsilon}$. $\Box$
\begin{Cor}
$\H^1$ factorizes $\G_{r,2}\times \G'$ as a quotient
for every $r\geq2$.
\end{Cor}
{\sc Proof.}
Observe that $\ch^1=\cg_{2,2}$, take into account Remark~\ref{rpupsilon}
and apply Proposition~\ref{Gpupsilon}. $\Box$
\vskip.2truecm
The following corollaries are straightforward.
\begin{Cor}
$\H^1$ factorizes $\G_{r,2}\times \G_{r,2}\times\cdots\times\G_{r,2}$
as a quotient for every $r\geq2$.
\end{Cor}
\begin{Cor}
$\H^1$ factorizes $\H^1\times \H^1\times\cdots\times\H^1$ as a quotient.
\end{Cor}
\begin{Rem}{\rm
Although $\H^1$ factorizes $\H^1\times \H^1\times\cdots\times\H^1$
as a quotient, one notice that there are a few nontrivial
h-epimorphisms between these groups. In fact, they are all of the form
\[
L:\overbrace{\H^1\times \H^1\times\cdots\times\H^1}^{n-times}\lra\H^1\,\quad
L(a_1,\ldots,a_n)=J(a_k)
\]
for some fixed $k\in\{1,\ldots,n\}$ and some h-isomorphism
$J:\H^1\lra\H^1$.
}\end{Rem}
%
%
%
%
%
%
%
%
%
%
%
%
%
\section{Examples of $(\G,\M)$-regular sets}\label{SecEx}

Existence of different types of intrinsically regular sets
in a given graded group $\P$ depends on the corresponding
algebraic factorizations.
In correspondence to the fact that $\R^k$ factorizes a graded
group $\M$ as a subgroup, we have the following
\begin{The}\label{legalg}
Let $\M$ be a graded group and let $n$ be the maximum over all dimensions
of commutative subalgebras contained in the first layer.
Then the family of $k$-dimensional Legendrian
$C^1$ smooth submanifolds is nonempty if and only if $1\leq k\leq n$
and it coincides with that of $(\R^k,\M)$-regular sets of $\M$.
\end{The}
{\sc Proof.}
Every $k$-dimensional commutative subalgebra of the first layer
also represents a trivial example of $k$-dimensional Legendrian submanifold.
This shows, by definition of $n$, that this family is nonempty for every
$1\leq k\leq n$. Recall that a $k$-dimensional Legendrian submanifold
is locally parametrized by a $C^1$ contact mapping defined on an open subset
of $\R^k$. In other words, it is locally parametrized by
a continuously h-differentiable contact mapping with injective differential.
Then Theorem~\ref{PdifContact} shows that this mapping is continuously 
P-differentiable. Furthermore, since in this case classical differentiability coincides with h-differentiability, then formula \eqref{hpdif} shows that the classical differential coincides with the P-differential.
Then the P-differential of the local parametrization is an injective h-homomorphism. As a result,
by Remark~\ref{rksubrm} the P-differential is an h-monomorphism.
Taking into account Definition~\ref{gmsinm},
we have shown that every $C^1$ smooth Legendrian submanifold contained
in $\M$ is an $(\R^k,\M)$-regular set of $\M$.
Clearly, we have also proved that there do not exist $k$-dimensional
Legendrian submanifold when $k>n$, since the injective P-differential
of the local parametrization would imply the existence 
of a commutative subalgebra in the first layer of $\cM$ with
dimension greater than $k$, that is a contradiction.

Now, let $S$ be an $(\R^k,\M)$-regular set, where clearly $1\leq k\leq n$.
Again Theorem~\ref{PdifContact} shows that continuously P-differentiable
mappings on $\R^k$ exactly correspond to contact mappings of class $C^1$,
then $S$ is a $C^1$ smooth Legendrian submanifold. $\Box$
\vskip.2truecm
It is well known that in the 
Heisenberg group $\H^n$ the only intrinsically regular sets
are $(\H^n,\R^k)$-regular sets and $(\R^k,\H^n$-regular sets
with $1\leq k\leq n$, \cite{FSSC6}.
\begin{Def}{\rm
We say that a homogeneous subgroup of $\H^n$ is {\em horizontal} if its 
Lie algebra is contained in the first layer $\cv$ and
{\em vertical} if it contains the second layer $\cz$.}
\end{Def}
\begin{Rem}{\rm
Notice that horizontal subgroups are always commutative, whereas
homogeneous normal subgroups exactly characterize
vertical subgroups.
}\end{Rem}
\begin{Pro}\label{anyfact}
Every couple of complementary subgroups of $\H^n$
is formed by a horizontal and a vertical subgroup.
\end{Pro}
{\sc Proof.}
Let $\ca$ and $\cb$ be the corresponding homogeneous subalgebras
of $A$ and $B$, respectively.
From our hypothesis and Proposition~\ref{homdec}, it follows that
$\ca\oplus\cb=\ch^n$. 
In particular, $Z\in\ca\oplus\cb$ hence we can suppose that
$\ca$ is the subalgebra containing a vector $W=W_0+\lambda Z$ with $\lambda\neq0$. Homogeneity of $\ca$ gives
$W_0+r\lambda Z\in\ca$ for every $r>0$, therefore $W_0\in\ca$
and $Z\in\ca$. This shows that $\ca$ is an ideal of $\ch^n$.
Clearly, the same argument shows that $\cb$ cannot have
vectors of the form $U_0+\mu Z$ with $\mu\neq0$, otherwise
$Z$ would belong to $\cb$, that conflicts with the decomposition
$\ca\oplus\cb=\ch^n$. This shows that $\cb$ is contained
in the first layer of $\ch^n$ and concludes the proof. $\Box$
\vskip.2truecm
The special factorization of $\H^n$ cannot be extended to all H-type groups.
\begin{Exa}\label{nnorhtype}{\rm
Let us consider the complexified Heisenberg group $\H^1_2$
along with the bases $(R_0,R_1,R_2,R_3)$ and $(Z_1,Z_2)$,
defined in Lemma~\ref{baslemcompheis}.
We define the homogeneous commutative subalgebras
$\ca=\span\{R_0,R_3,Z_1\}$ and $\cb=\span\{R_1,R_2,Z_2\}$.
Clearly, $\ca\oplus\cb=\ch_2^1$, hence from Proposition~\ref{homdec}
it follows that $A=\exp\ca$ and $B=\exp\cb$ are complementary subgroups of $\H_2^1$.
On the other hand, $[R_0,R_2]=Z_2\notin\ca$ and $[R_0,R_1]=Z_1\notin\cb$
impliy that both $\ca$ and $\cb$ are not ideals of $\ch_2^1$.}
\end{Exa}
On the other hand, we have the following
\begin{Pro}\label{clsh21}
Let $N$ and $H$ be complementary subgroups of $\H_2^1$,
where $N$ is a normal. Then $N$ contains the center
of the group, $\mbox{\rm top-}\dim(N)\in\{4,5\}$
and $H$ is commutative and horizontal, namely, $H\subset\exp\cv$.
\end{Pro}
{\sc Proof}
Let $N=\exp \cn$ and $H=\exp \ch$. If $\dim(\cn)\geq3$, then
$\cn\supset\cz$ by formula \eqref{keyhtype}. Then in this case
$\dim(\ch)\leq 3$ and $\ch\subset\cv$. By Proposition~\ref{nnexcomm},
it follows that $\dim(\ch)\leq 2$. Then the case $\dim(\cn)=3$ cannot occur.
If $\dim(\cn)\leq2$, then $\cn\subset\cz$, then Proposition~\ref{nnexcomm}
prevents the existence of $\ch\subset\cv$, that should be commutative.
As a result, the only allowed possibilities are
$\dim(\cn)=4$ and $\dim(\cn)=5$. $\Box$
\begin{The}\label{intregh21}
The only intrinsically regular sets of $\H_2^1$ are 
contained in the following list
\begin{enumerate}
\item
$(\H_2^1,\R)$-regular sets, 
\item
$(\H_2^1,\R^2)$-regular sets,
\item
$(\R^2,\H_2^1)$-regular sets,
\item
$(\R,\H_2^1)$-regular sets.
\end{enumerate}
\end{The}
{\sc Proof.}
By Proposition~\ref{GlinH} and Proposition~\ref{clsh21},
we are allowed to consider only
level sets defined through continuously P-differentiable mappings
with values in $\R^k$, $k=1,2$ and defined on an open set of $\H_2^1$.
Every $C^1$ mapping $f:\Omega\subset\H_2^1\lra\R^2$ with surjective differential 
defines a $(\H_2^1,\R^2)$-regular set, since it is continuously
P-differential and its P-differential is an h-epimorphism
due to Proposition~\ref{h21r2}. This allows us to apply Theorem~\ref{implth}.
We argue in the same way for $C^1$ mappings $f:\Omega\subset\H_2^1\lra\R$,
due to Proposition~\ref{rfactg}.

Analogously, Proposition~\ref{GlinHmono} and Proposition~\ref{clsh21}
allow us to consider only image sets from open subsets of $\R^k$,
$k=1,2$ with values in $\H_2^1$. By Theorem~\ref{legalg},
$(\R,\H_2^1)$-regular sets and $(\R^2,\H_2^1)$-regular sets of
$\H_2^1$ exist and correspond to one dimensional and two dimensional
$C^1$ smooth Legendrian submanifolds of $\H_2^1$. $\Box$

\end{document}